%% file: main_arxiv.tex
\documentclass[a4paper,10pt
]{amsart}
\usepackage{mathptmx}      
\usepackage{mathrsfs}
\usepackage{latexsym}
\usepackage{amssymb}
\usepackage{bm}
\usepackage[usenames,dvipsnames]{xcolor}
\usepackage{graphicx,pifont}
\usepackage{grffile}
\usepackage{multicol,multirow}
\usepackage{tikz}
\usetikzlibrary[patterns]
\usepackage{comment}
\usepackage{array}

\usepackage{pgfplots}
\pgfplotsset{compat=1.17}
\usepgfplotslibrary{fillbetween}
\usepackage{mwe}
\usepackage{subcaption,float}
\newenvironment{nouppercase}{%
  \renewcommand{\uppercasenonmath}[1]{}}{}

\usepackage{txfonts}
\usepackage{mathtools}
\usepackage[shortlabels]{enumitem}
\usepackage{hyperref}
\usepackage{cleveref}
\usepackage{fancyvrb}
\usepackage{float}
\usepackage{booktabs}
\usepackage[abs]{overpic}
\crefname{equation}{}{}
\crefname{enumi}{}{}
\crefname{figure}{Figure}{Figure}
\crefname{subsection}{subsection}{subsections}
\crefname{lemma}{Lemma}{Lemma}
\crefname{proposition}{Proposition}{Proposition}

\newtheorem{theorem}{Theorem}[section]

\theoremstyle{definition}

\theoremstyle{remark}
\newtheorem{remark}[theorem]{Remark}
\numberwithin{equation}{section}
\numberwithin{figure}{section}

\newcommand{\bs}[1]{\boldsymbol{#1}}

\newcommand{\bb}[1]{{\mbox{\sffamily \bfseries #1}}}

\newcommand{\nrm}[1]{\left|\left|#1\right|\right|}
\newcommand{\eqdef}{\stackrel{\mathrm{def}}{=\joinrel=}}

\newcommand{\powth}{{\textrm{\scriptsize th}}}

\newcommand{\node}{\textrm{node}}
\newcommand{\pmc}{\textrm{pmc}}
\newcommand{\cfl}{{\textrm{cfl}}}
\newcommand{\freq}{{\textrm{freq}}}

\newcommand{\frob}{{\textrm{Frob}}}

\newcommand{\fom}{{\textrm{fom}}}
\newcommand{\rom}{{\textrm{rom}}}
\newcommand{\vms}{{\textrm{vms}}}

\newcommand{\ddiv}{{\textrm{div}}}
\definecolor{vargreen}{rgb}{0.0, 0.5, 0.0}

\usepackage[foot]{amsaddr}
\usepackage[margin=2.7cm]{geometry}
\usepackage{systeme}
\usepackage{graphicx} 
\usepackage[]{units}
\usepackage{color}
\usepackage{mathrsfs}
\usepackage{bm}
\usepackage[ansinew]{inputenc}
\usepackage{float}
\usepackage{amsmath}
\usepackage{isomath} 
\usepackage{amsthm}
\usepackage{amssymb}

\usepackage{hyperref}

\newcommand{\red}{\color{black}}

\usepackage{bm}
\usepackage{graphicx,color, amsfonts,amsmath}

\usepackage{mwe}

\usepackage{amssymb}
\providecommand{\norm}[1]{\lVert#1\rVert}

\DeclareMathAlphabet\mathbfcal{OMS}{cmsy}{b}{n}

\begin{document}

\title[Embedded domain ROMs for  the shallow water hyperbolic equations with the SBM]{Embedded domain Reduced Basis Models for the shallow water hyperbolic equations \\ with the Shifted Boundary Method}

\author{Xianyi Zeng\textsuperscript{1}
}
\email{xzeng@utep.edu}
\address{\textsuperscript{1}Department of Mathematical Sciences, University of Texas at El
Paso, United States.}
\author{Giovanni Stabile\textsuperscript{2}
}
\email{gstabile@sissa.it}

\author{Efthymios N. Karatzas\textsuperscript{3}
}
\address{\textsuperscript{2}SISSA, International School for Advanced Studies, Mathematics Area, mathLab, Trieste, Italy.}
\address{\textsuperscript{3}Department of Mathematics, School of Applied Mathematical and Physical Sciences, NTUA, Athens, Greece, and FORTH Institute of Applied and Computational Mathematics, Heraclion, Crete, Greece.}
\email{karmakis@math.ntua.gr}

\author{Guglielmo Scovazzi\textsuperscript{4}}
\email{guglielmo.scovazzi@duke.edu}
\address{\textsuperscript{3}Civil and Environmental Engineering Department, Duke University, Durham, NC 27708, United States.}

\author{Gianluigi Rozza\textsuperscript{2}}
\email{grozza@sissa.it}

\subjclass[2010]{78M34, 97N40, 35Q35}

\keywords{Shallow waters, 
          geometrical parametrization,
          embedded FEM,
          Shifted Boundary Method,
          reduced basis}
  
\date{{\red{\today}}}

\dedicatory{}
\input{abstract.tex}
\begin{nouppercase}
\maketitle
\end{nouppercase}
\input{section_intro.tex}
\input{section_sbm.tex}

\input{section_pmc.tex}

\input{section_rom.tex}
\input{section_num.tex}

\input{section_concl.tex}
\input{section_aknw.tex}
\bibliographystyle{amsplain_gio}
\bibliography{bibfile_sissa}
\end{document}

%% file: abstract.tex
\begin{abstract}
\normalsize
We consider fully discrete embedded finite element approximations for a shallow water hyperbolic  problem and its reduced-order model. Our approach is based on a fixed background mesh and an embedded reduced basis. The Shifted Boundary Method for spatial discretization is combined with an explicit predictor/multi-corrector time integration to integrate in time the numerical solutions to the shallow water equations, both for the full and reduced-order model. In order to improve the approximation of the solution manifold also for geometries that are untested during the offline stage, the snapshots have been pre-processed by means of an interpolation procedure that precedes the reduced basis computation. The methodology is tested on geometrically parametrized shapes with varying size and position.
\end{abstract}

%% file: section_intro.tex
\section{Introduction}\label{sec:intro}
The computational cost associated with the numerical solution of partial
differential equations might be in some cases prohibitive. This is happening,
for example, when the numerical solution is required in nearly real time or a
when large number of system configurations need to be tested. Shape optimization
problems are a typical example of the latter case, where a large number of
different geometrical configurations need to be analyzed to converge to an optimal
solution. Reduced order models demonstrated to be a viable approach to reduce
the computational burden and have been developed for a large variety of
different linear and nonlinear problems \cite{BeOhPaRoUr17,handbook}.

In recent times, immersed/embedded/unfitted methods have seen a great development from  the seminal ideas of Peskin~\cite{Peskin1977}. The key ideas in embedded methods is the use of grids that are not body-fitted, in which the geometry of the shapes to be simulated is immersed by way of computational geometry techniques.
In this work, we base the reduced order models on the Shifted Boundary Method (SBM), which is an embedded/unfitted finite element method originally proposed for the Poisson, Stokes and incompressible Navier-Stokes equations~\cite{MaSco17_2,MaSco17_3} and recently extended to wave equations and shallow water equations (SWE)~\cite{TSong:2018a}.
In the SBM, a surrogate boundary is introduced in proximity of the true immersed boundary, and the boundary conditions are imposed on the surrogate boundary, with appropriate corrections that rely on Taylor expansions~\cite{MaSco17_2,MaSco17_3}.
The SBM does not require complicated data structures and numerical quadratures to integrate the governing equations on cut element, typical of cutFEM/XFEM approaches.
Compared to other embedded finite element methods such as XFEM or cutFEM~\cite{JChessa:2003a,AGerstenberger:2008a,TPFries:2010a}, the SBM has also the advantage that the degrees of freedom (unknowns) stay the same for varying geometries, hence existing reduced order model methodologies are more easily adapted. Specifically, the total number of unknowns in SBM is determined by the background mesh and it is independent of the location of the embedded geometry.
In contrast, XFEM and other enriched finite element methods alike introduce new degrees of freedom such as the Heaviside functions within cut elements, hence the total number of unknowns typically varies with the embedded geometry locations and/or depends on a computationally expensive cutting of elements procedure.

In this article we focus our attention on projection-based reduced order models
specifically tailored to geometrically parametrized problems \cite{StabileZancanaroRozza2020,TezzeleDemoStabileMolaRozza2020}. The
idea is to combine the recently proposed Shifted Boundary Method
\cite{MaSco17_2,MaSco17_3,TheoreticalPoissonAtallahCanutoScovazzi2020} with Reduced Order Models based on the Proper Orthogonal Decomposition (POD) with Galerkin projection (SBM-ROM). This
combination, that has been recently proposed in previous works in a different
setting
\cite{KaratzasStabileAtallahScovazziRozza2018,KaratzasStabileNouveauScovazziRozza2018,KaratzasStabileNouveauScovazziRozzaNS2019},
allows to avoid the map of all the parametrized solutions to a common reference
geometry, see also \cite{KaRo21,KaKaTra21}. 
In this paper, the embedded methodology introduced in the
mentioned research works is extended to shallow water equations with explicit
time marching schemes. 
The idea of merging embedded approaches with reduced order models has been
proposed also in \cite{Nouy2011} where a fictitious domain method was
coupled with a Proper Generalized Decomposition approach to study uncertain
geometries. In \cite{BaFa2014} the authors proposed a projection based reduced
order model starting from an embedded full order simulation applied to evolving
interfaces.

With embedded simulations it is in fact easy to work with a common background
mesh also in the case of large geometrical changes.
By comparison, body-fitted meshes often require sophisticated re-meshing techniques, when complex geometrical deformations are present, and maintaining the topology of the underlying mesh is a difficult task.

In addition, we introduce a new approach to handle degrees
of freedom located in the \emph{``out of interest/ghost"} region which is based on a
radial basis function interpolation. We denote this new approach as SBM-iROM.
As shown in the numerical examples, this
approach allows to partially reduce the drawback associated with the slow decay
of the Kolmogorov N-width when dealing with embedded full order models.

The article is organized as follows: in section~\ref{sec:HF} we introduce the
mathematical formulation of the full order problem, the associated weak
formulation and the details concerning the specific discretization strategy.
Section~\ref{sec:ROM} describes in details the approach used for the
construction of the reduced order model with a focus on the relevant changes
required for the specific full order model formulation and introduces the employed interpolation preprocessing. In
section~\ref{sec:num_exp} we introduce three numerical examples to show the
properties and accuracy of the proposed methodology. Finally, in
section~\ref{sec:conclusions} we report some conclusions and outlooks for future
developments.
\section{The model problem and the full order approximation} \label{sec:HF}
{Before introducing the shallow water model, we briefly review the relevant literature:
for linearized shallow water equations arising from the equations of acoustics we refer to the work 
\cite{KaTsuYamKawFu92}, where only the generation of the second harmonic wave is considered (the higher order harmonics being neglected) under the assumption of weak non-linearity, while a set of uncoupled equations for the primary and secondary wave is discretized in space by a finite element method, and then solved by using the Newmark-$\beta$ integration scheme for time. 
In \cite{ScoCa12}, weak boundary conditions are considered for the hyperbolic structure of the wave equation based on stabilized methods and the variational multiscale analysis as well as we cite \cite{ScoSoZe17} for linear elastodynamics and \cite{TSong:2015a} for Nitsche and wave propagation problems. 
In \cite{HAUKE98} one may see the shallow water equations as a symmetric advective-diffusive systems with source terms solved using SUPG and GLS stabilized methods via a predictor multi-corrector algorithm. 
The work of \cite{Hauke02} treats two-dimensional shallow water equations applied to solve practical irrigation problems with large friction coefficients, dry bed conditions and singular infiltration terms. 
For a three-step shallow water flow explicit scheme, using parallel computing, and tidal flow in Tokyo Bay, we refer to \cite{KaItoBeTez95}, while a finite element method for the analysis of nearshore current, which is one of the principal currents in coastal seas analyzing two main characteristics of the wave, i.e. direction and height is introduced in \cite{KawKas84}. 
Surface wave motion handled by the Helmholtz equation is studied in \cite{KawKas85}. 
 A new combinative method of boundary-type finite elements and boundary solutions to study wave diffraction-refraction and harbour oscillation problems is presented in \cite{KasKawMu88} with model the mild-slope equation proved as an effective and accurate method for water surface wave problems. 
Large-scale computation of storm surges and tidal flows carried out with finite element methods are discussed in \cite{KaSaBeTez97}. 
 The work of \cite{KaMa94} presents an adaptive boundary-type finite element method for wave diffraction-refraction in harbors model based on the mild-slope equation and an arbitrary reflection condition and in \cite{KaSaHiKa88} boundary-type finite element method has been investigated and applied to the Helmholz and mild-slope equations useful for practical analysis. 
Numerical analysis of tsunamis applying the finite element method  based on the shallow water wave equation with the Lax-Wendroff finite difference method used also for the analysis of the Tokachi-oki Earthquake tsunami problem and compared with the tide gauge records is investigated in \cite{KaTaNoYoTa78}.
A stabilized Residual Distribution scheme for the simulation of shallow water flows with a nonlinear variant of a Lax-Friedrichs type discretization is proposed in  \cite{Ri09} while in \cite{TaKaTaTeTa10} a numerical modeling of seismic body wave propagation problems is introduced.
Finally for a space-time SUPG formulation of the shallow-water equations based on a proposed embedded with extension
of the proposed approach to discontinuous Galerkin methods or residual redistribution schemes  considering the complex
morphology of the ocean coastlines in real scenarios with the treatment of complex coastlines as reflective walls in the framework of large scale simulations of fine scale urban floods we refer to \cite{Smith75}. 
}
\subsection{Strong formulation of the shallow water problem}\label{par:strong form}
\begin{figure}[ht]\centering
  \includegraphics[width=.6\textwidth]{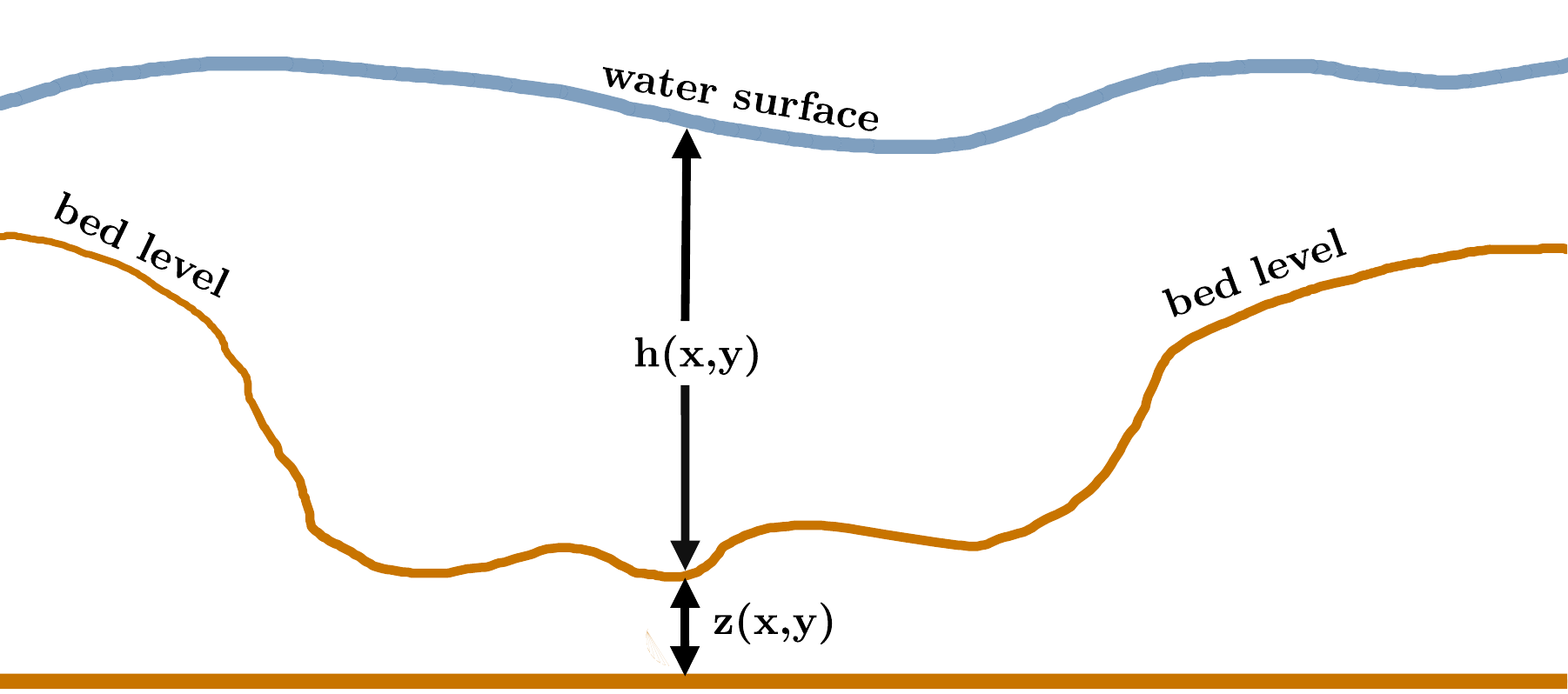}
  \caption{A sketch of the bathymetry $z(x,y)$ and the water height $h(x,y)$.}
  \label{fg:W_height_bathymetry}
\end{figure}

The shallow-water equations, also known as de Saint-Venant equations \cite{Barre1871}, are a system of hyperbolic partial differential equations simulating the behavior of a free surface of a fluid \cite{KouDou20, Vre86} when the depth of the fluid bed is shallow compared to the characteristic horizontal spatial length.
Such system of equations is derived from the Navier-Stokes equations after integrating through the depth and observing that, since the horizontal length scale is much greater than the vertical length scale, the vertical component of the fluid velocity field is small compared to the horizontal component and the vertical gradients of pressure are nearly hydrostatic. This allows to conclude that the horizontal velocity field is constant  over the entire depth of the fluid. Settings in fluid dynamics where the horizontal length scale is much greater than the vertical length scale are common and widely applicable, see e.g. atmospheric and oceanic modeling. As mentioned in the previous paragraph, the propagation of a tsunami wave can also be simulated efficiently with the shallow-water equations until it reaches the coast, see e.g. \cite{TSong:2018a} and references therein. 
The shallow-water equations we examine are simulating thin layers of fluid in hydrostatic balance (with constant density), with upper and lower bounds: a free surface and the bed/ground topography, respectively. These  systems can be formulated as
\begin{subequations}\label{eq:fom_swe}
\begin{eqnarray}
h_t + (hv_1)_x
+ (hv_2)_y &=& 0,
\\
(hv_1)_t
+ (hv_1^2 + \frac{1}{2} g h^2)_x
+ (hv_1v_2)_y
 &=& S_1 ,
 \\
(hv_2)_t  
+ (hv_1v_2)_x
 + (hv_2 ^2 + \frac{1}{2}gh^2)_y
&=& S_2,
\end{eqnarray}
\end{subequations}
where we indicate by $h$ the height of the water, $z$ the bathymetry of the water bed, and by 
$\eta = h + z$ the free surface measure,  as it is visualized in Figure \ref{fg:W_height_bathymetry}.
The two-dimensional position and velocity vectors in Cartesian coordinate axes are denoted by $x = (x , y)$ and
$v = (v_1 , v_2)$. The quantity $\partial/\partial t$ expresses derivation with respect to time, and
the source term $ \mathbf{S} = (S_1,S_2) = gh( S_{o_1} - S_{f_1}, S_{o_2} - S_{f_2})$, 
for $S_{o_i} = - \partial z/\partial x_i$, $i=1,2$   and $S_{f_i} = f^2 v_i (v^2_1 + v_2^2 )^{\frac{1}{2}} h^{-\frac{4}{3}}$ to be  the slope of the river bed or ocean floor, and the friction (in terms of the Manning's roughness coefficient $f$) respectively. We note that other types of force can be added to the source term if needed, like wind stress, Coriolis forces, etc.


\subsection{Weak formulation on body-fitted grids}
\label{sec:fom_weak}
The SWE is solved by a stabilized piecewise linear nodal finite element method (FEM) using two-dimensional triangular elements.
We first introduce the notation and the weak formulation assuming a body-fitted grid and the extension to embedded boundaries using the Shifted Boundary Method is given in the next subsection.
Let us denote the vector of conserved variables by $\bs{U}=[h,\;hv_1,\;hv_2]^T$ and the flux function by $\bs{F}(\bs{U})$ with components $\bs{F}_i,\;1\le i\le d=2$; then the equation~(\ref{eq:fom_swe}) can be written compactly as:
\begin{equation}\label{eq:fom_weak_swe}
  \partial_t\bs{U} + \nabla\cdot\bs{F}(\bs{U}) = \bs{Z}\;,\quad (\bs{x},\;t)\in\mathcal{D}\times[0,\;T]\;,
\end{equation}
where $\bs{Z}=[0,\;S_1,\;S_2]^T$.
Let $\mathcal{T}$ be a tessellation of the domain $\mathcal{D}$, we define the globally continuous piecewise linear solution and trial spaces:
\begin{equation}\label{eq:fom_weak_space}
  \mathcal{S}^h = \mathcal{V}^h = \left\{\bs{U}^h\in[C(\mathcal{D})]^{1+d}\,:\; \bs{U}^h\big|_K\in[\mathbb{P}_1(K)]^{1+d},\;\ \forall K\in\mathcal{T}\right\}\subset H^1(\mathcal{D})\times H^{\ddiv}(\mathcal{D})\;,
\end{equation}
where no boundary conditions are specified as they are enforced weakly.
The semi-discrete stabilized finite element method states as finding $\bs{U}^h:[0,\;T]\mapsto\mathcal{S}^h$ such that for all $0\le t\le T$ and $\bs{W}^h\in\mathcal{V}^h$:
\begin{equation}\label{eq:fom_weak_fem}
  (\bs{W}^h,\;\partial_t\bs{U}^h-\bs{Z})_{\mathcal{D}} - (\nabla\bs{W}^h,\;\bs{F})_{\mathcal{D}} + a_{\vms}(\bs{W}^h,\;\bs{U}^h) + b(\bs{W}^h,\;\bs{U}^h) = 0\;.
\end{equation}
Here $(\cdot,\;\cdot)_{\mathcal{D}}$
denotes the standard inner product for $L^2$ scalar, vectorial, and tensorial functions, $a_{\vms}$ contains a variational multiscale stabilization (VMS) term, and $b$ contains all boundary terms.

The purpose of the VMS term is to prevent spurious oscillations due to the fact that equal-order interpolation is used for the velocity and fluid height variables.
The spatial and full differential operators in linearized form are defined as:
\begin{equation}\label{eq:fom_weak_dop}
  \mathcal{L}\bs{U} = \sum_{i=1}^d\bs{A}_i\partial_{x_i}\bs{U}\;,\quad
  \mathcal{L}_t\bs{U} = \partial_t\bs{U} + \mathcal{L}\bs{U}\;,
\end{equation}
where $\bs{A}_i\eqdef\partial\bs{F}_i/\partial\bs{U}$.
Then we may compute the dual operator of $\mathcal{L}$ as $\mathcal{L}^\ast=\sum_{i=1}^d\bs{A}_i^T\partial_{x_i}$ and define the VMS term as:
\begin{equation}\label{eq:fom_weak_vms}
  a_{\vms}(\bs{W}^h,\;\bs{U}^h) = (\mathcal{L}^\ast\bs{W}^h,\;\tau_{\vms}\bs{A}_0^{-1}(\mathcal{L}_t\bs{U}^h-\bs{Z}))_{\mathcal{D}}\;,
\end{equation}
where $\tau_{\vms}$ 
is a parameter that scales with time and it is computed as $\tau_{\vms}=c_{\vms}\Delta t/2$ with $\Delta t$ being the time step size and $c_{\vms}=O(1)$ being a user defined parameter.
A fixed value $c_{\vms}=2.0$ is used in all computations in this work.
The matrix $\bs{A}_0$ is the Jacobian matrix converting the conservative variables $\bs{U}$ to primitive ones $\bs{Y}=[h,\;v_1,\;v_2]^T$:
\begin{equation}\label{eq:fom_weak_scaling}
  \bs{A}_0 = \begin{bmatrix}
    1 & 0 & 0 \\ v_1 & h & 0 \\ v_2 & 0 & h
  \end{bmatrix}\;.
\end{equation}
Other options for the scaling matrix $\bs{A}_0^{-1}$ include the Jacobian matrix between conservative and entropy variables, as suggested in~\cite{ETadmor:1984a}; we adopt~(\ref{eq:fom_weak_scaling}) for simplicity.

{\bf Remark}. In principle, a discontinuity capturing (artificial viscosity) term could also be included for improved stability, especially when strong shocks are present.
In this work, the problems considered do not involve strong shocks and we omit the discontinuity capturing operator. We point out however that the proposed methodology can be applied to the case in which shock capturing operators are used.

Lastly, the boundary condition can be classified into different types depending on the information available, such as incoming or outgoing flows, or subcritical or supercritical velocities.
We shall only consider the boundary terms that are relevant to the test problem here, namely the Neumann condition, the subcritical inflow condition, and the subcritical outflow condition, defined on the boundary portions $\Gamma_N$, $\Gamma_{I;sub}$, and $\Gamma_{O;sub}$, respectively; for a complete list of all boundary conditions the readers are referred to~\cite{SoMaScoRi17}.
Particularly, the three boundary conditions are given by:
\begin{align}
  \label{eq:fom_weak_nbc}
  &\bs{v}\cdot\bs{n} = v_N\;, &&\bs{x}\in\Gamma_N\;, \\
  \label{eq:fom_weak_isubbc}
  &h\bs{v}\cdot\bs{n} = m_{I;sub}\;,\quad \bs{v}\cdot\bs{\tau} = 0\;, &&\bs{x}\in\Gamma_{I;sub}, \\
  \label{eq:fom_weak_osubbc}
  &h\bs{v}\cdot\bs{n} = m_{O;sub}\;, &&\bs{x}\in\Gamma_{O;sub}\;.
\end{align}
Here $\bs{n}$ is the outer unit normal on the boundary, $\bs{\tau}$ is the tangent vector,  and $v_N$, $m_{I;sub}<0$, and $m_{O;sub}>0$ are prescribed normal velocity, inflow mass rate, and outflow mass rate, respectively.
Note that in principle (\ref{eq:fom_weak_isubbc}) and (\ref{eq:fom_weak_osubbc}) are only valid when the Fruode number is smaller than unity, otherwise the boundary condition needs to be switched to supercritical ones; the latter scenario, however, does not occur in the tests in this work.
While we shall describe the Neumann condition by assuming a general $v_N$, in all cases the value of $v_N$ is set to zero so that the boundary either represent a slippery wall or a symmetry plane.
To this end, the boundary term $b$ is given by:
\begin{equation}\label{eq:fom_weak_bterm}
  b(\bs{W}^h,\;\bs{U}^n) = \langle\bs{W}^h,\bs{H}^h_N\rangle_{\Gamma_N} + \langle\bs{W}^h,\bs{H}^h_{I;sub}\rangle_{\Gamma_{I;sub}} + \langle\bs{W}^h,\bs{H}^h_{O;sub}\rangle_{\Gamma_{O;sub}}\;,
\end{equation}
where the angle brackets are the inner product on a general boundary piece $\Gamma$ that is defined as $\langle\bs{W},\bs{H}\rangle_\Gamma=\int_{\Gamma}\bs{W}\cdot\bs{H}d\Gamma$.
Particularly, the vectors $\bs{H}^h_N$, $\bs{H}^h_{I;sub}$, and $\bs{H}^h_{O;sub}$ are given respectively by:
\begin{equation}\label{eq:fom_weak_h_bc}
  \bs{H}_N^h = v_N\begin{bmatrix}h \\ h\bs{v}\end{bmatrix} + \frac{1}{2}gh^2\begin{bmatrix}0\\\bs{n}\end{bmatrix}\;,\quad
  \bs{H}_{I;sub}^h = m_{I;sub}\begin{bmatrix}1\\(\bs{v}\cdot\bs{n})\bs{n}\end{bmatrix} + \frac{1}{2}gh^2\begin{bmatrix}0\\\bs{n}\end{bmatrix}\;,\quad
  \bs{H}_{O;sub}^h = m_{O;sub}\begin{bmatrix}1\\\bs{v}\end{bmatrix} + \frac{1}{2}gh^2\begin{bmatrix}0\\\bs{n}\end{bmatrix}\;.
\end{equation}
In the case of an embedded boundary, the boundary term (especially the Neumann one) needs to be modified, as described next.

%% file: section_sbm.tex
\subsection{Discretization: the Shifted Boundary Method}
In this subsection, we introduce the basic aspects of the Shifted Boundary Method adapted to the shallow water equations~\cite{MaSco17_2,MaSco17_3,TSong:2015a,TSong:2018a}. We consider a surrogate domain $\tilde {\mathcal{D}}$ and boundary $\tilde \Gamma$ together with the true computational domain $\mathcal{D}$ and its boundary $\Gamma$, as shown in Figure~\ref{SurrogateMesh} and Figure~\ref{fig:SBM}. 
We indicate by $\bm{\tilde{n}}$ the unit outward-pointing normal to the surrogate boundary $\tilde \Gamma$, which is distinct from the outward-pointing normal  $\bm n$  to ${\Gamma}$ as seen in Figure~\ref{fig:ntd}. 
$\tilde  \Gamma$ consists of the edges that are closest in some sense to the true boundary $\Gamma$, as shown in Figure~\ref{fig:ntd}. 
\begin{figure}
\centering
\begin{subfigure}[t]{.4\textwidth}\centering       
  \begin{overpic}[width=\textwidth,tics=10]{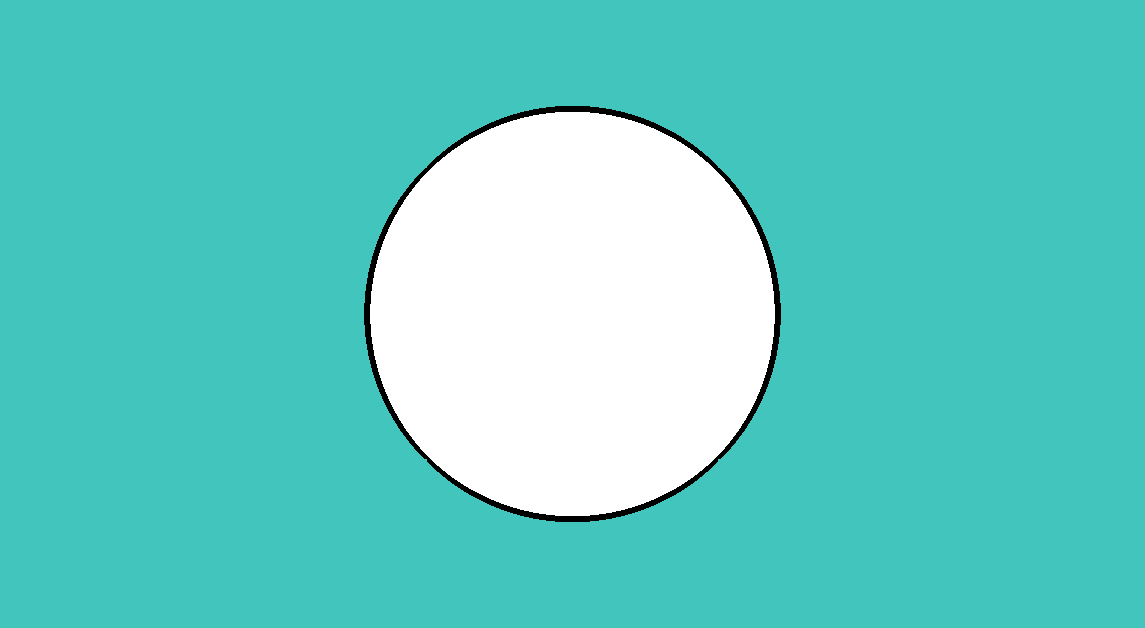}
    \put (30,75) {\Large$\mathcal{D}$}
    \put (70,25) {\Large$\Gamma$}
  \end{overpic}
  \caption{The geometry $\mathcal{D}$ surrounding a disk and its boundary $\Gamma$.}
  \label{SurrogateMesh_geom}
\end{subfigure}$\quad$
\begin{subfigure}[t]{.4\textwidth}\centering
  \begin{overpic}[width=\textwidth,tics=10]{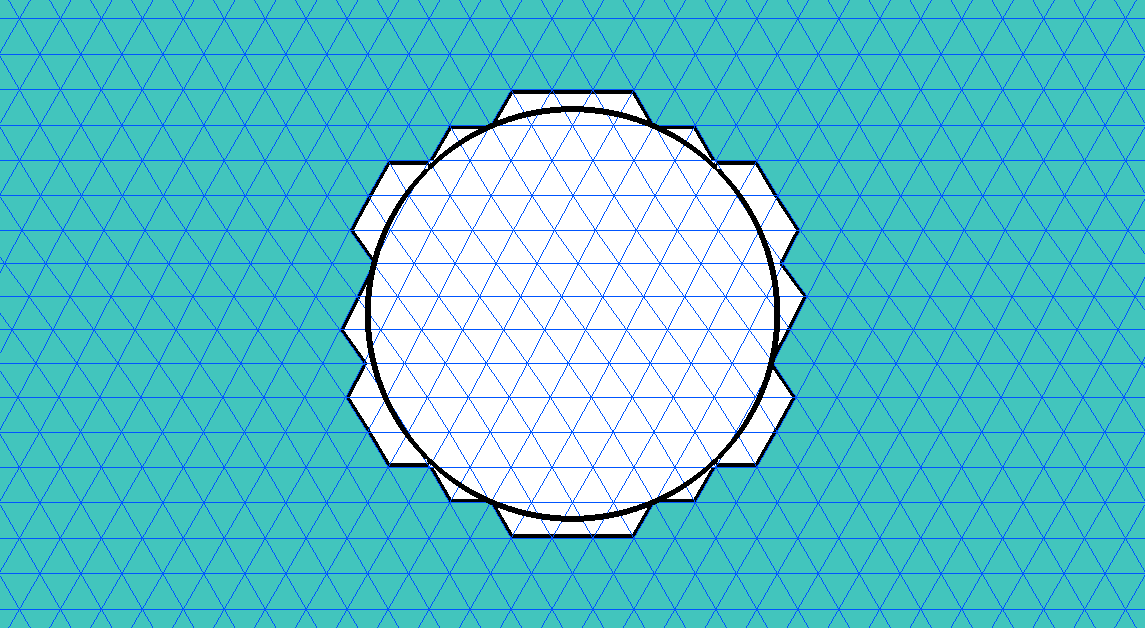}
    \put (30,75) {\Large$\tilde{\mathcal{D}}$}
    \put (70,25) {\Large$\Gamma$}
    \put (45,25) {\Large$\tilde{\Gamma}$}
    \put (61.5,45) {\large ``Ghost area''}
  \end{overpic}
  \caption{The SBM surrogate geometry $\tilde{\mathcal{D}}$, the surrogate boundary $\tilde{\Gamma}$, and ghost area.}
  \label{SurrogateMesh_sbm}
\end{subfigure}
\caption{Embedded geometry: (A) The geometry of a disk and (B) the SBM surrogate geometry attached with the ghost area.}
\label{SurrogateMesh}
\end{figure}
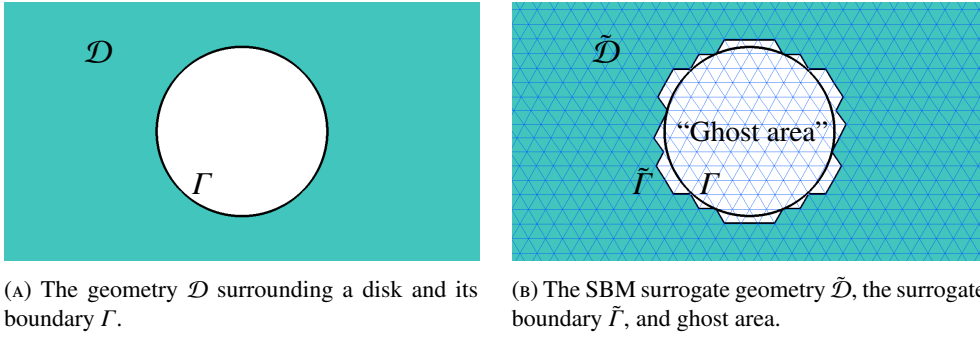
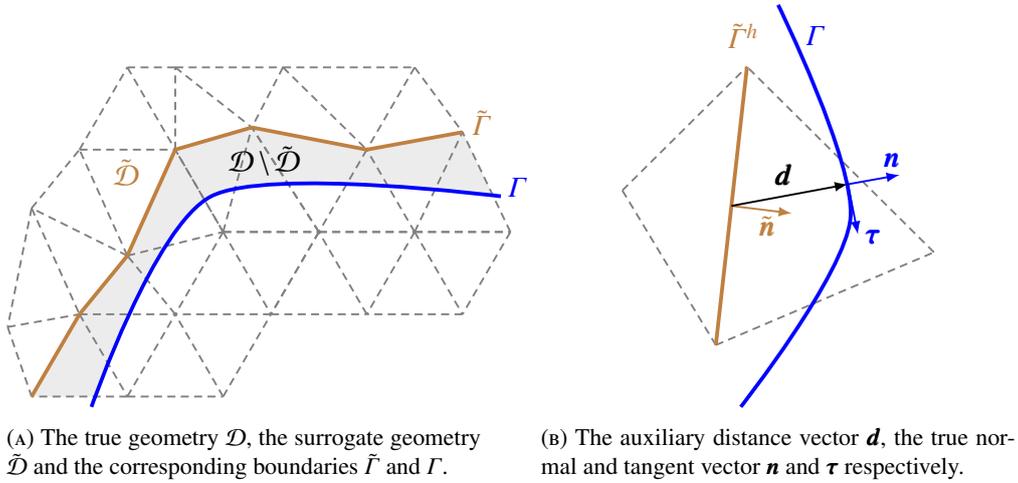
\begin{figure}
\centering
\begin{subfigure}[t]{.4\textwidth}
	\centering
\begin{tikzpicture}[scale=0.63]
\draw [black, line width = 0.95mm, draw=none,name path=surr] plot coordinates { (-2,-3.4641) (-1,-1.73205) (0,-0.5) (1,1.73205) (2.6,2.2) (5,1.73205) (7,2.1) (7.62,0.8) }; 
\draw [blue,  line width = 0.95mm, draw=none,name path=true] plot[smooth] coordinates {(-0.7,-3.4641) (1.75,0.75) (8.25,0.75)}; 
\tikzfillbetween[of=true and surr, split]{gray!15!};
\draw[line width = 0.25mm,densely dashed,gray] (-1,1.73205) -- (0,3.4641);
\draw[line width = 0.25mm,densely dashed,gray] (0,3.4641) -- (2,0);
\draw[line width = 0.25mm,densely dashed,gray] (1,1.73205) -- (1,3.4641);
\draw[line width = 0.25mm,densely dashed,gray] (1,3.4641) -- (0,3.4641);
\draw[line width = 0.25mm,densely dashed,gray] (1,3.4641) -- (2.6,2.2);
\draw[line width = 0.25mm,densely dashed,gray] (1,3.4641) -- (3.25,3.4641);
\draw[line width = 0.25mm,densely dashed,gray] (3.25,3.4641) -- (2.6,2.2);
\draw[line width = 0.25mm,densely dashed,gray] (3.25,3.4641) -- (5,1.73205);
\draw[line width = 0.25mm,densely dashed,gray] (3.25,3.4641) -- (6,3.4641);
\draw[line width = 0.25mm,densely dashed,gray] (6,3.4641) -- (5,1.73205);
\draw[line width = 0.25mm,densely dashed,gray] (6,3.4641) -- (7,2.1);
\draw[line width = 0.25mm,densely dashed,gray] (0,-0.5) -- (-2,0.5);
\draw[line width = 0.25mm,densely dashed,gray] (-2,0.5) -- (-1,1.73205);
\draw[line width = 0.25mm,densely dashed,gray] (-2,0.5) -- (-1,-1.73205);
\draw[line width = 0.25mm,densely dashed,gray] (-2,0.5) -- (-2.5,-2);
\draw[line width = 0.25mm,densely dashed,gray] (-2.5,-2) -- (-1,-1.73205);
\draw[line width = 0.25mm,densely dashed,gray] (-2.5,-2) -- (-2,-3.4641);
\draw[line width = 0.25mm,densely dashed,gray] (0,-0.5) -- (-1,1.73205);
\draw[line width = 0.25mm,densely dashed,gray] (-1,1.73205) -- (1,1.73205);
\draw[line width = 0.25mm,densely dashed,gray] (0,-0.5) -- (2,0);
\draw[line width = 0.25mm,densely dashed,gray] (2,0) -- (1,1.73205);
\draw[line width = 0.25mm,densely dashed,gray] (1,1.73205) -- (0,-0.5);
\draw[line width = 0.25mm,densely dashed,gray] (2,0) -- (2.6,2.2);
\draw[line width = 0.25mm,densely dashed,gray] (2.6,2.2) -- (1,1.73205);
\draw[line width = 0.25mm,densely dashed,gray] (2,0) -- (4,0);
\draw[line width = 0.25mm,densely dashed,gray] (4,0) -- (2.6,2.2);
\draw[line width = 0.25mm,densely dashed,gray] (4,0) -- (2.6,2.2);
\draw[line width = 0.25mm,densely dashed,gray] (2.6,2.2) -- (5,1.73205);
\draw[line width = 0.25mm,densely dashed,gray] (5,1.73205) -- (4,0);
\draw[line width = 0.25mm,densely dashed,gray] (4,0) -- (6,0);
\draw[line width = 0.25mm,densely dashed,gray] (6,0) -- (5,1.73205);
\draw[line width = 0.25mm,densely dashed,gray] (6,0) -- (7,2.1);
\draw[line width = 0.25mm,densely dashed,gray] (7,2.1) -- (5,1.73205);
\draw[line width = 0.25mm,densely dashed,gray] (6,0) -- (8,0);
\draw[line width = 0.25mm,densely dashed,gray] (8,0) -- (7,2.1);
\draw[line width = 0.25mm,densely dashed,gray] (0,-0.5) -- (-1,-1.73205);
\draw[line width = 0.25mm,densely dashed,gray] (-1,-1.73205) -- (1,-1.73205);
\draw[line width = 0.25mm,densely dashed,gray] (2,0) -- (1,-1.73205);
\draw[line width = 0.25mm,densely dashed,gray] (1,-1.73205) -- (0,-0.5);
\draw[line width = 0.25mm,densely dashed,gray] (2,0) -- (3,-1.73205);
\draw[line width = 0.25mm,densely dashed,gray] (3,-1.73205) -- (1,-1.73205);
\draw[line width = 0.25mm,densely dashed,gray] (4,0) -- (3,-1.73205);
\draw[line width = 0.25mm,densely dashed,gray] (2,0) -- (4,0);
\draw[line width = 0.25mm,densely dashed,gray] (4,0) -- (3,-1.73205);
\draw[line width = 0.25mm,densely dashed,gray] (3,-1.73205) -- (5,-1.73205);
\draw[line width = 0.25mm,densely dashed,gray] (5,-1.73205) -- (4,0);
\draw[line width = 0.25mm,densely dashed,gray] (4,0) -- (6,0);
\draw[line width = 0.25mm,densely dashed,gray] (6,0) -- (5,-1.73205);
\draw[line width = 0.25mm,densely dashed,gray] (6,0) -- (7,-1.73205);
\draw[line width = 0.25mm,densely dashed,gray] (7,-1.73205) -- (5,-1.73205);
\draw[line width = 0.25mm,densely dashed,gray] (6,0) -- (8,0);
\draw[line width = 0.25mm,densely dashed,gray] (8,0) -- (7,-1.73205);
\draw[line width = 0.25mm,densely dashed,gray] (0,-3.4641) -- (-2,-3.4641);
\draw[line width = 0.25mm,densely dashed,gray] (-2,-3.4641) -- (-1,-1.73205);
\draw[line width = 0.25mm,densely dashed,gray]  (-1,-1.73205) -- (0,-3.4641);
\draw[line width = 0.25mm,densely dashed,gray] (0,-3.4641) -- (1,-1.73205);
\draw[line width = 0.25mm,densely dashed,gray] (0,-3.4641) -- (2,-3.4641);
\draw[line width = 0.25mm,densely dashed,gray] (2,-3.4641) -- (1,-1.73205);
\draw[line width = 0.25mm,densely dashed,gray] (2,-3.4641) -- (3,-1.73205);
\draw [line width = 0.5mm,blue, name path=true] plot[smooth] coordinates {(-0.75,-3.681818) (1.75,0.75) (7.8,0.75)};
\draw[line width = 0.5mm,brown] (1,1.73205) -- (2.6,2.2);
\draw[line width = 0.5mm,brown] (2.6,2.2) -- (5,1.73205);
\draw[line width = 0.5mm,brown] (5,1.73205) --  (7,2.1);
\draw[line width = 0.5mm,brown] (1,1.73205) -- (0,-0.5);
\draw[line width = 0.5mm,brown] (0,-0.5) -- (-1,-1.73205);
\draw[line width = 0.5mm,brown] (-1,-1.73205) -- (-2,-3.4641);
\node[text width=0.5cm] at (7.6,2.3) {\large${\color{brown}\tilde{\Gamma}}$};
\node[text width=3cm] at (2.1,1.25) {\large${\color{brown}\tilde{\mathcal{D}}}$};
\node[text width=0.5cm] at (8.35,0.95) {\large${\color{blue}\Gamma}$};
\node[text width=3cm] at (4.5,1.5) {\large${\mathcal{D}} \setminus \tilde{\mathcal{D}} $};
\end{tikzpicture}
\caption{The true geometry ${\mathcal{D}}$, the surrogate geometry $\tilde{\mathcal{D}}$ and the corresponding boundaries $\tilde{\Gamma}$ and $\Gamma$.}
\label{fig:SBM}
\end{subfigure}
\qquad
\begin{subfigure}[t]{.4\textwidth}\centering
	\begin{tikzpicture}[scale=0.82]
\draw[line width = 0.25mm,densely dashed,gray] (0,0.5) -- (-1.5,3);
\draw[line width = 0.25mm,densely dashed,gray] (-1.5,3) -- (0.5,5);
\draw[line width = 0.25mm,densely dashed,gray] (0,0.5) -- (3.5,2);
\draw[line width = 0.25mm,densely dashed,gray] (3.5,2) -- (0.5,5);
\draw [line width = 0.5mm,blue, name path=true] plot[smooth] coordinates {(0.4,-0.5) (2.16,2.55) (1.0,6)};
\draw[line width = 0.5mm,brown] (0,0.5) -- (0.5,5);
\node[text width=0.5cm] at (0.5,5.5) {\large${\color{brown}\tilde{\Gamma}^h}$};
\node[text width=0.5cm] at (1.75,5.5) {\large${\color{blue}\Gamma}$};
\node[text width=0.5cm] at (1.25,3.25) {\large$\bm{d}$};
\node[text width=0.5cm] at (3,3.5) {\color{blue} \large$\bm{n}$};
\node[text width=0.5cm] at (1.0,2.43) {\color{brown} \large$\tilde{\bm{n}}$};
\node[text width=0.5cm] at (2.7,2.25) {\color{blue} \large$\bm{\tau}$};
\draw[->,brown, line width = 0.25mm,-latex] (0.25,2.75) -- (1.22,2.63);
\draw[->,line width = 0.25mm,-latex] (0.25,2.75) -- (2.12,3.1);
\draw[->, blue, line width = 0.25mm,-latex] (2.12,3.1) -- (2.28,2.29);
\draw[->, blue, line width = 0.25mm,-latex] (2.12,3.1) -- (2.95,3.25);
\end{tikzpicture}
    \caption{The auxiliary distance vector $\bm{d}$, the true normal and tangent vector $\bm{n}$ and  $\bm{\tau}$ respectively.}
    \label{fig:ntd}
\end{subfigure}
%
\caption{The true and the surrogate geometry and boundary, and the SBM related quantities $\bm{d}$, $\tau $, $\bm{n}$, $\tilde{\bm{n}}$.}
\label{fig:surrogates}
\end{figure}
The mapping 
\label{eq:defMmap}
\begin{align}
\bm{M}^{h}: \tilde{\Gamma} \to \Gamma,
\end{align}
is introduced, similarly to \cite[Section 2.1]{TheoreticalPoissonAtallahCanutoScovazzi2020}, which maps any point $\tilde{\bm{x}} \in \tilde{\Gamma}$ on the surrogate boundary, to a point $\bm{x} = \bm{M}^{h}(\tilde{\bm{x}})$ on the true physical boundary $\Gamma$.
Through $\bm{M}^{h}$, an auxiliary distance vector function $\bm{d}_{\bm{M}^{h}}$ is defined as 
\begin{align}
\label{eq:Mmap}
\bm{d}_{\bm{M}^{h}} (\tilde{\bm{x}})
\, = \, 
\bm{x}-\tilde{\bm{x}}
\, = \, 
[ \, \bm{M}^h-\bm{I} \, ] (\tilde{\bm{x}})
\; .
\end{align}
For brevity, and denoting by $\bm{\nu}$ the unit vector, we set $\bm{d} = \bm{d}_{\bm{M}^{h}} $ where $\bm{d} = \| \, \bm{d} \, \| \bm{\nu}$.
In the case of smooth surfaces with one type of boundary condition, the mapping $\bm{M}^{h}$ corresponds to the closest point projection and $\bm{\nu} = \bm{n}$, see e.g. Figure~\ref{fig:ntd}. 
%
The general construction and analysis of $\bm{M}^h$ are detailed in~\cite{TheoreticalPoissonAtallahCanutoScovazzi2020}, including the case when corners are present or the closures of the Dirichlet boundary $\Gamma_D$ and the Neumann boundary $\Gamma_N$ have non-empty intersections. Following the analysis therein, we may hypothesize that ${\bm{M}^{h}}$ is continuous and Lipschitz. This assumption makes sense since the true surface is smooth between edges and corners.
The mapping ${\bm{M}^{h}}$ can be used to extend the unit normal vector $\bm n$ from the boundary $\Gamma$ to the surrogate boundary $\tilde \Gamma$ as
$\bm \bar {\bm n}(\tilde {\bm x}) \equiv {\bm n}({\bm{M}^{h}}(\tilde {\bm x}))$.
In the next sections, we use the short-hand notation $\bm n(\tilde {\bm x})$, which means $\bm \bar {\bm n}(\tilde {\bm x})$ at a point ${ \tilde {\bm x}}\in \tilde\Gamma$.
In a similar way, we extend the boundary conditions on $\Gamma$ to the boundary $\tilde  \Gamma$ of the surrogate domain. 
\subsubsection{Semi-discrete shifted boundary weak formulation}
We can now set up the semi-discrete SBM weak formulation relying on the surrogate domain $\tilde{\mathcal{D}}$ whose boundary is composed of both body-fitted and embedded portions.
With a slight abuse of notation, we indicate with $\Gamma$ the portion of the boundary $\partial \tilde{\mathcal{D}}$ that is body-fitted, and by $\tilde{\Gamma}$ the portion of the boundary $\partial \tilde{\mathcal{D}}$ that is embedded (see, e.g., Figure~\ref{fg:num_cyl_setup}, where $\Gamma$ consists of the four exterior edges and $\tilde{\Gamma}$ is the surrogate boundary associated with the internal circle).
We discretize now $\tilde{\mathcal{D}}$ using a mesh triangulation ${\mathcal{{\tilde{D}}_T}}$ consisting of triangles $K$ that belong to a tessellation $\mathcal T$.
The weak SBM formulation is now given as: \\

Find $\bs{U}^h:\;[0,\;T]\mapsto\tilde{\mathcal{S}}^h$, such that for all $0\le t\le T$ and $\bs{W}^h\in\tilde{\mathcal{V}}^h$:
\begin{equation}\label{eq:sbm_weak}
  (\bs{W}^h,\;\partial_t\bs{U}^h-\bs{Z})_{\tilde{\mathcal{D}}} - (\nabla\bs{W}^h,\;\bs{F})_{\tilde{\mathcal{D}}}+a_{\vms}(\bs{W}^h,\;\bs{U}^h) + b(\bs{W}^h,\;\bs{U}^h) + \tilde{b}(\bs{W}^h,\;\bs{U}^h) = 0\;,
\end{equation}
where $\tilde{\mathcal{S}}^h$ and $\tilde{\mathcal{V}}^h$ are obtained by replacing $\mathcal{D}$ in (\ref{eq:fom_weak_space}) with the surrogate $\tilde{\mathcal{D}}$.
The VMS term remains the same as~(\ref{eq:fom_weak_vms}), except the inner product is evaluated on $\tilde{\mathcal{D}}$; the boundary term $b$ is given by~(\ref{eq:fom_weak_bterm} for all body-fitted boundaries; and the second boundary term $\tilde{b}$ is given below in the case of an embedded Neumann boundary ($\tilde{\Gamma}=\tilde{\Gamma}_N$):
\begin{equation}\label{eq:sbm_weak_nbc}
  \tilde{b}(\bs{W}^h,\;\bs{U}^h) = \langle\bs{W}^h,\;\tilde{\bs{H}}^h_N\rangle_{\tilde{\Gamma}_N}\;,\quad
  \tilde{\bs{H}}^h_N = \left((v_N-\bs{n}^T\nabla\bs{v}\bs{d})\bs{n}\cdot\tilde{\bs{n}}+(\bs{v}\cdot\bs{\tau})\bs{\tau}\cdot\tilde{\bs{n}}\right)\begin{bmatrix}h \\ h\bs{v}\end{bmatrix} + \frac{1}{2}gh^2\begin{bmatrix}0\\\tilde{\bs{n}}\end{bmatrix}\;.
\end{equation}
The first term in the definition of $\tilde{\bs{H}}^h_N$ attempts to enforce a normal velocity {\it shifted} to the surrogate interface by using Taylor series expansions:
\begin{align*}
  \bs{v}(\tilde{\bs{x}})\cdot\tilde{\bs{n}}(\tilde{\bs{x}}) 
  &= \left[(\bs{v}(\tilde{\bs{x}})\cdot\bs{n})\bs{n} + (\bs{v}(\tilde{\bs{x}})\cdot\bs{\tau})\bs{\tau}\right]\cdot\tilde{\bs{n}} \\
  &\approx \left\{\left[(\bs{v}(\bs{x})-\nabla\bs{v}(\tilde{\bs{x}})\bs{d})\bs{n}\right]\cdot\bs{n}+(\bs{v}(\tilde{\bs{x}})\cdot\bs{\tau})\bs{\tau}\right\}\cdot\tilde{\bs{n}}
  = (v_N-\bs{n}^T\nabla\bs{v}\bs{d})\bs{n}\cdot\tilde{\bs{n}}+(\bs{v}\cdot\bs{\tau})\bs{\tau}\cdot\tilde{\bs{n}}\;.
\end{align*}
Embedded boundary conditions of other types can be derived similarly but are omitted here, for the sake of brevity, since they are not applied in the test problems considered here.

%% file: section_pmc.tex
\subsubsection{Time discretization}\label{sec:pmc}
An explicit predictor/multi-corrector (PMC) time integration is used to march the numerical solutions in time~\cite{TSong:2015a,TSong:2018a}.
To this end, let us denote the time ordinate by $t$ and the spatial coordinate by $\bs{x}$.
Furthermore, the superscript $^n$ designates a variable associated with the $n^{\powth}$ time step $t^n$ and a subscript $_A$ designates a variable associated with a mesh node $\bs{x}_A$, $1\le A\le n_{\node}$, where $n_{\node}$ is the total number of nodes.

The discrete solution vector at $t^n$ is given by:
\begin{equation}\label{eq:pmc_fesol}
  \bs{U}(t^n,\bs{x}) \approx \bs{U}^n(\bs{x}) \eqdef \sum_{A=1}^{n_{\node}}\bs{U}_A^nN_A(\bs{x})\; ,
\end{equation}
where $N_A(\bs{x})$ is the piecewise linear shape function associated with the node $\bs{x}_A$ and $\bs{U}_A^n$ are the degrees of freedom at the same node.
Choosing $\bs{W}^h=\bs{e}_iN_A(\bs{x})$ for all $1\le A\le n_{\node}$ and $1\le i\le 1+d$ in the weak formulation~(\ref{eq:fom_weak_fem}), where $\bs{e}_i$ is the $i^\powth$ unit vector in $\mathbb{R}^{1+d}$, we obtain a system of ordinary differential equations:
\begin{equation}\label{eq:pmc_ode}
  \bb{M}\dot{\bb{U}}(t) + \mathcal{R}(\bb{U}) = 0\;,
\end{equation}
where $\bb{U}(t): \mathbb{R}^+ \mapsto \mathbb{R}^{(1+d)n_{\node}}$ contains all degrees of freedom at all nodes, i.e., $\bs{U}_A\in\mathbb{R}^{1+d}$ for all $1\le A\le n_{\node}$ and $\bb{M}$ is the diagonal lumped mass matrix.
The operator $\mathcal{R}: \mathbb{R}^{(1+d)n_{\node}}\mapsto\mathbb{R}^{(1+d)n_{\node}}$ contains the (spatial) residual at each node for each component of the solution variable and includes the standard finite element terms in continuous Galerkin formulation, the stabilization term, and the Shifted Boundary Method terms that arise in the transmission boundary condition at the surrogate boundary.

To update the solution vector from $\bb{U}^n$ to $\bb{U}^{n+1}$, the explicit PMC method can be considered as a fixed-point iteration approximation for the midpoint rule,
\begin{equation}\label{eq:pmc_midpoint}
  \bb{M}(\bb{U}^{n+1}-\bb{U}^n) + \delta t^n\mathcal{R}((\bb{U}^n+\bb{U}^{n+1})/2) = 0\;,
\end{equation}
where $\delta t^n = t^{n+1}-t^n>0$ is the time step size.
In PMC, one seeks successive approximations to $\bb{U}^{n+1}$, denoted by $\bb{U}^{(k)}$ where $k=0,1,\cdots$.
In particular $\bb{U}^{(0)}\eqdef\bb{U}^n$ is the ``predictor'' of $\bb{U}^{n+1}$ and once $\bb{U}^{(k)}, k\ge0$ is computed, one computes the next iterate $\bb{U}^{(k+1)}$ by:
\begin{equation}\label{eq:pmc_pmc}
  \bb{M}(\bb{U}^{(k+1)}-\bb{U}^n) + \delta t^n\mathcal{R}((\bb{U}^n+\bb{U}^{(k)})/2) = 0\;.
\end{equation}
In practice, one terminates the iteration after a preset number of correctors (each $\bb{U}^{(k)}$ with $k\ge1$ is known as a ``corrector''), i.e.,~\cref{eq:pmc_pmc} is performed for $k=0,\cdots,n_{\pmc}-1$ where typical values for $n_{\pmc}$ is between $2$ and $4$.
Finally, the update $\bb{U}^{n+1}=\bb{U}^{(n_{\pmc})}$ is applied.

\begin{remark}
If one sets $n_{\pmc}=1$, the method is equivalent to the explicit forward-Euler method, whereas $n_{\pmc}=2$ gives an implementation of the second-order Runge-Kutta scheme.
\end{remark}

%% file: section_rom.tex
\section{Reduced order model with a POD-Galerkin method}\label{sec:ROM}
%
The reduced order model proposed here is based on a POD-Galerkin approach. It means that the underlying system of equations is projected onto a linear subspace of smaller dimension spanned by a reduced number of global basis functions (POD modes). There are different techniques to generate this linear subspace and here we rely on the POD \cite{sirovich1987turbulence}. 
The overall methodology is based on the classic offline-online splitting approach \cite{handbook}, which is briefly recalled in what follows. 
\subsection*{Offline Stage}
During the offline stage we start with a parametric partial differential equation which is parame\-trized by means of a $p$-dimensional parameter vector $\bm{\mu} \in \mathcal{P}$. The full order model is then solved for a finite dimensional set of training points in $\{\bm{\mu}_i\}_{i=1}^{N_{{\mu}_\text{train}}} \subset \mathcal{P}$. 

In the current framework we construct one linear subspace which includes both parameter and time variations. This means that the snapshots matrix on which the POD is based is assembled as:
\begin{equation}
S_U = [U(t_1,\bm{\mu}_1), U(t_2,\bm{\mu}_1), \dots, U(t_{N_t},\bm{\mu}_{N_{\mu_\text{train}}}) )] \in \mathbb{R}^{N_h \times N_s},
\end{equation}
where $N_h = (1+d)n_{\text{node}}$ is the dimension of the FOM solution vector, $N_s = N_t\cdot N_{{\mu}_\text{train}}$ is the total number of stored snapshots where $N_t$ is the number of time steps in which we store the solution and $N_{{\mu}_\text{train}}$ is the number of training samples in the parameter space. The POD modes are generated using the method of snapshots as originally proposed by Sirovich \cite{sirovich1987turbulence} which relies on the solution of an eigenproblem on the correlation matrix:
\begin{equation}\label{eq:pod_cmat}
C = S_U^T S_U,
\end{equation}
and on the computation of the POD bases exploiting the resulting eigenvalues $\{\lambda_i\}_{i=1}^{N_s}$ ad eigenvectors $\{\psi_i\}_{i=1}^{N_s}$:
\begin{equation}
\phi_i = \frac{1}{\sqrt{\lambda_i}} S_U \psi_i.
\end{equation}
This operation results into to the POD space:
\begin{equation}\label{eq:pod_space}
\bm{\Phi}_{\text{full}} = \text{span} (\phi_1,\dots, \phi_{N_s}).
\end{equation}
\begin{remark}
Based on the eigenvalue decomposition of the correlation matrix $C$ it is possible to discard some of the modes and to create a POD space $\bm{\Phi}_{\text{r}}$ that includes a limited number of the computed POD modes. This space will be employed for the subsequent Galerkin projection. The correlation matrix has been computed relying onto the Frobenius inner product. Other options are possible (such as $L^2$ or $H^1$ norms) but in the current setting, due to the fact that the underlying background mesh is made of finite elements with a similar size, for the sake of simplicity, we have decided to rely on the Frobenius norm to calculate the POD modes. 
\end{remark}
\subsection{Snapshots Interpolation, the SBM-iROM}
\label{sec:rom_interp}
\begin{figure}[ht]\centering
\includegraphics[width=.8\textwidth]{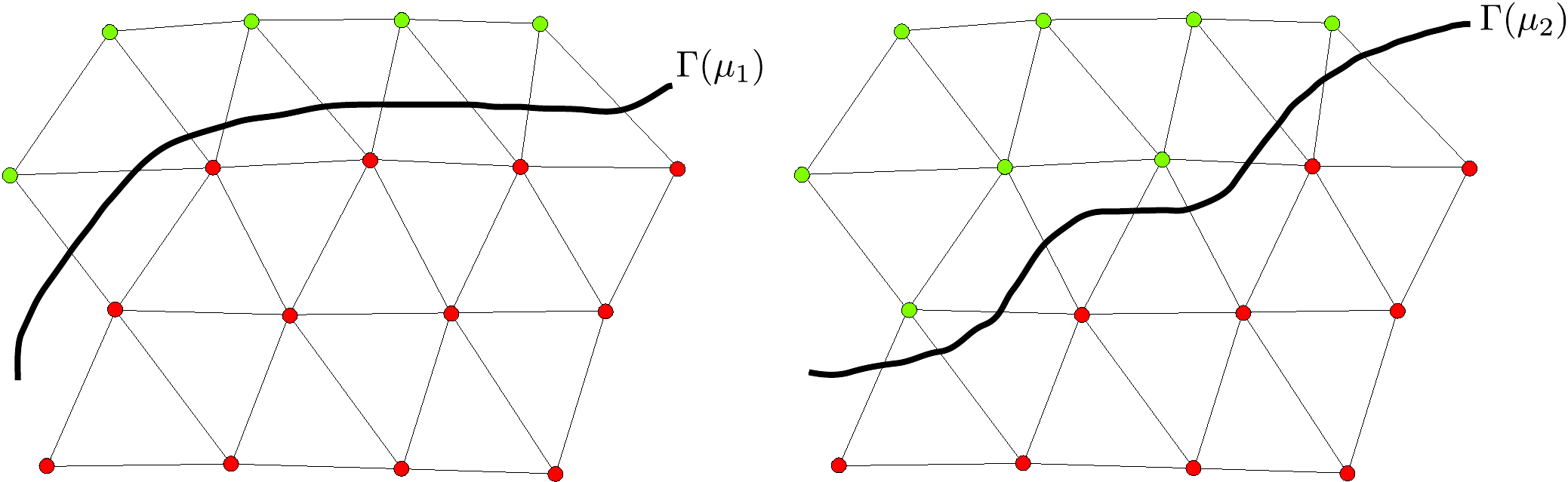}
\caption{Visualization of active (green nodes) and inactive nodes (red nodes) for two different parametrized geometries}
\label{fg:active_inactive}
\end{figure}

The procedure described above is rather straightforward for body-fitted meshes.
For unfitted meshes there are additional complexities that need to be
addressed. One of them is related to the inactive nodes that belong to the so-called \emph{ghost area}. These nodes are embedded into the body and do not play a role in the computation.
As depicted in Figure~\ref{fg:active_inactive}, the number and location of the
inactive nodes is changing depending on the shape of the parametrized geometry.
In order to create a global basis function that can be used for any new
parameter configuration it is necessary to handle also the inactive nodes. A
possible approach would be to set them to a constant value. However, this approach would introduce a discontinuity in the solution field on the jump between
active and inactive nodes. Therefore, we have decided to preprocess each
snapshot with an interpolation strategy, in order to avoid discontinuity caused
by deactivated elements. Another possible option to avoid this issue would be to
compute an harmonic extension from the boundary to the deactivated nodes
\cite{KaBaRO18}. The method proposed here has the advantage of ensuring a sufficient level of smoothness without the solve of the additional partial differential equation problem required by the harmonic extension.   The
interpolant function has been evaluated for each snapshot considering only the
active nodes and the values in the inactive nodes have been evaluated using the
interpolant function. Each solution field has been replaced by:
\begin{equation}\label{eq:interpolation}
U(\bm{x}) \to U_I(\bm{x}) = \sum_{k = 1}^{N_{\text{active}}}\omega_k \varphi (\norm{\bm{x} - \bm{x}_k}) + \sum_{j=1}^m b_j p_j (\bm{x}),
\end{equation}
where $\omega_k$ are weights that needs to be determined imposing the
interpolation condition, $\varphi (\norm{\bm{x} - \bm{x}_k})$ are radial basis
functions where $\bm{x}_k$ are the coordinates of the active nodes,
$\{p_i(\bm{x})\}_{i=1}^m$ are monomials that span the space of polynomials with
a specific degree. Adding polynomials to the RBF interpolant functions helps to
properly capture constant and linear features in the given data and ensures positive-definiteness of the RBF function, which in turn implies solvability of
the interpolation problem~\cite{buhmann2003radial}. The
coefficient vectors $\bm{\omega} = \left[\omega_1, \dots,
\omega_{N_\text{active}} \right]^T$ and $\bm{b} = \left[ b_1, \dots, b_m
\right]^T$ are obtained solving the following linear system of equations:
\begin{equation}\label{eq:rbf_system}
\begin{cases}
\left( \bm K + \sigma^2 \bm I \right) \bm \omega + \bm P \bm b = \bm d, \\
\bm P^T \bm \omega = \bm 0,
\end{cases}
\end{equation}
where $(\bm{K})_{ij} = \varphi(\norm{\bm{x}_i - \bm{x}_j})$, $(\bm{P})_{ij} = p_j(\bm{x}_i)$ and $\sigma$ is a smoothing parameter that eventually relaxes the interpolation condition in order to smoothen the ``interpolant'' function. For the specific case has been set to $\sigma = 0$. Selecting a RBF which is positive definite of order $q$ with $\bm P$ that has full column rank, the solution is unique provided that the degree of the monomial terms is at least $m=q-1$ \cite{Fasshauer2007,buhmann2003radial}. In the specific case, for each snapshot and solution field, the interpolant function has been constructed using observation only on the active nodes. The solution field on the inactive nodes has been reconstructed using the previously defined interpolant function:
\begin{equation}\label{eq:}
U(\bm{x}_{i_{\text{inactive}}}) = U_I(\bm{x}_{i_{\text{inactive}}}).
\end{equation}
Here we decided to rely on second-order polyharmonic splines that have the expression:
\begin{equation}\label{eq:rbf_basis}
\varphi = (\epsilon r)^2 \log{ (\epsilon r)},
\end{equation}
where $r = \norm{\bm{x} - \bm{c}}$ with $\bm{x}$ and $\bm{c}$ the evaluation points and the radial basis center respectively. $\epsilon$ is a scalar shape parameter that in this case is set $\epsilon = 1.0$. This function is conditionally positive definite provided that the order of the monomial is at least equal to $2$ \cite{Fasshauer2007,buhmann2003radial}. Therefore, in order to ensure the uniqueness of the solution we have chosen to use $m=2$ in expression \eqref{eq:interpolation}. 
\begin{figure}[ht]\centering
\begin{minipage}{\textwidth}
\centering
\begin{minipage}{0.48\textwidth}
\includegraphics[width=\textwidth]{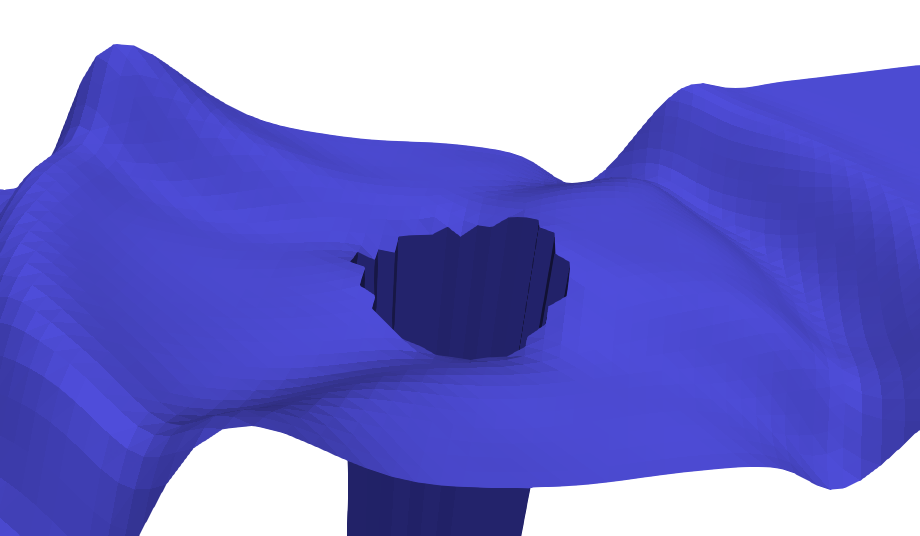}
\end{minipage}
\begin{minipage}{0.48\textwidth}
\includegraphics[width=\textwidth]{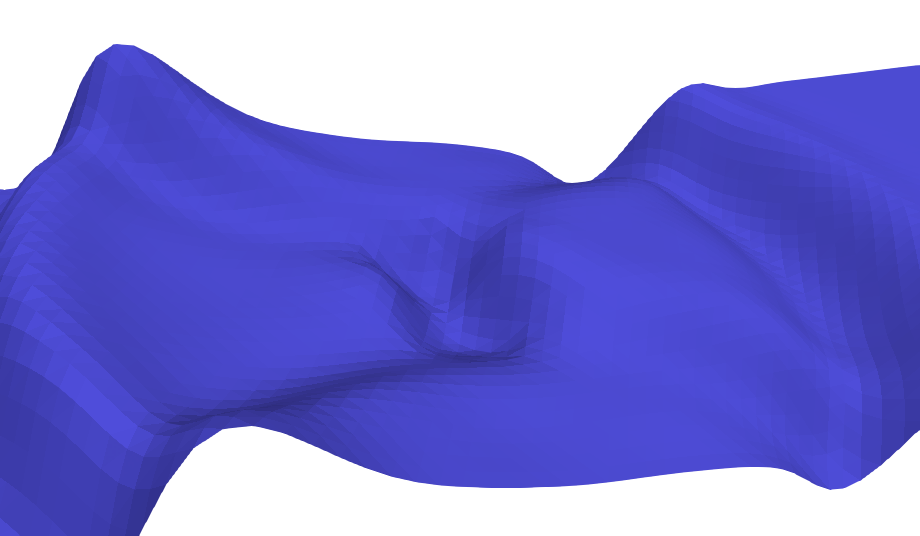}
\end{minipage}
\end{minipage}
\caption{Visualization of a solution snapshot for the height of the water field, with and without interpolation procedure}
\label{fg:with_without_interpolation}
\end{figure}
Using this setting, prior to the POD basis computation, the solution snapshots are replaced by the corresponding interpolant function. This operation permits to have a continuous solution field. We denote this model as SBM-iROM (Shifted Boundary Method with interpolation Reduced Order Model). In Figure~\ref{fg:with_without_interpolation} it is shown, for the height of the water field, a snapshot before and after the interpolation procedure. The hole on the left part of the picture corresponds to the inactive nodes which belongs to the embedded geometry and are set to zero by the full order model solver.

\subsection{POD-Galerkin Projection}
Once the POD bases are computed one can use them to approximate the solution field  $\bm U \in \mathbb{R}^{N_h}$ with its low dimensional representation $\bm U_r$ which is defined as a linear combination of $N_r$ global basis functions $\bm \Phi_r = \mbox{span} (\phi_1, \dots, \phi_{N_r})$:
\begin{equation}
\bm U \approx \bm U_r = \bm \Phi_r \bm a,
\end{equation}
where the basis functions $\phi_i \in \mathbb{R}^{N_h}$ are computed by POD applied on the interpolated snapshots. The coefficients of the POD expansion are then retrieved by means of Galerkin projection of the original system of equations onto the space spanned by the POD modes. Also at the reduced order level we have decided to use exactly the same time marching as the one employed at the full order level. This results in:
\begin{equation}
\bm \Phi_r^T \bb M \bm \Phi_r (\bm a^{(k+1)} - \bm a^n) + \delta t^n \bm \Phi_r^T \mathcal{R}(\bm \Phi_r (\bm a^{n}+\bm a^{(k)})/2) = 0.
\end{equation}
That can be reformulated as:
\begin{equation}\label{eq:update_a}
\bb M_r (\bm a^{(k+1)} - \bm a^n) + \delta t^n \bm \Phi_r^T \mathcal{R}(\bm \Phi_r (\bm a^{n}+\bm a^{(k)})/2) = 0.
\end{equation}
In the expression above, $\bb M_r = \bm \Phi_r^T \bb M \bm \Phi_r \in \mathbb{R}^{N_r \times N_r}$  can be precomputed and does not depend on the input parameters $\bm \mu$. The problem can be therefore expressed in terms of reduced coefficients $\bm a$ with the only difficulty that we will have to assemble the residual $\mathcal{R}$ also at every iteration of the ROM problem. Moreover, the residual term $\mathcal R$ has a nonlinear dependency with respect to the input parameter $\bm{\mu}$ which parametrizes the embedded geometry. Such possibly expensive residual computation could be substituted by point-wise evaluation of the residual function in some selected points of the domain using an hyper-reduction technique. However, in this work, since the main concern is to test the applicability of the methodology and the application of the proposed interpolation strategy, we have decided to assemble the full residual and to project it onto the reduced basis spaces at each iteration. 
%
%
%

%% file: section_num.tex
\section{Numerical experiments}\label{sec:num_exp}
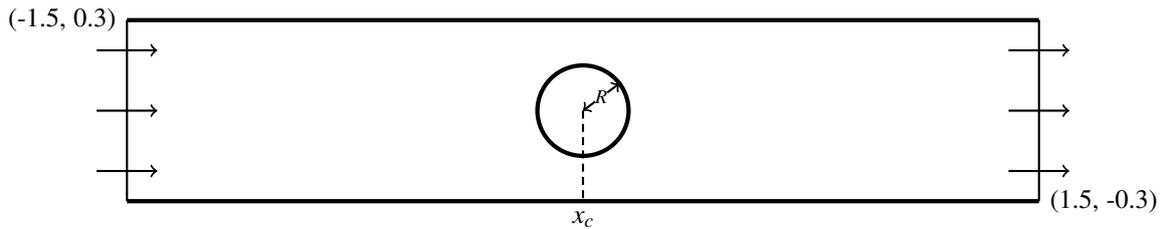
\begin{figure}[ht]\centering
  \begin{tikzpicture}[line width=1.6, scale=4.0]
    \draw (-1.5,-0.3) -- (1.5,-0.3);
    \draw (-1.5, 0.3) -- (1.5, 0.3);
    \draw [line width=0.8] (-1.5,-0.3) -- (-1.5,0.3);
    \draw [line width=0.8] ( 1.5,-0.3) -- ( 1.5,0.3);
    \draw ( 0.0, 0.0) circle (0.15);
    \foreach \y in {-0.2,0.0,0.2}
    {
      \draw [line width=0.8, ->] (-1.6,\y) -- (-1.4,\y);
      \draw [line width=0.8, ->] ( 1.4,\y) -- ( 1.6,\y);
    }
    \draw (-1.5,0.3) node [left] {(-1.5,\;0.3)};
    \draw (1.5,-0.3) node [right] {(1.5,\;-0.3)};
    \draw [densely dashed, line width=0.8] (0.0,0.0) -- (0.0,-0.3) node [below] {$x_c$};
    \draw [<-, line width=0.8] (0.0,0.0) -- (0.04,0.03);
    \draw [->, line width=0.8] (0.08,0.06) -- (0.12,0.09);
    \draw (0.06,0.045) node {$\scriptstyle R$};
  \end{tikzpicture}
  \caption{Symmetric SWE flow past a stationary cylinder in a straight channel, together with the two parameters characterizing the configuration in the numerical examples, namely the cylinder radius $R$ and the $x$-coordinate of the cylinder center $x_c$.}
  \label{fg:num_cyl_setup}
\end{figure}

We consider three test cases featuring symmetric SWE flow past a stationary cylinder, based on the configuration in Figure~\ref{fg:num_cyl_setup}.
The straight channel is represented by the computational domain $\Omega=[-1.5,\;1.5]\times[0.3,\;0.3]$.
Denoting the unit outer normal vector to $\partial\Omega$ by $\bs{n}$, the upper and lower boundary conditions are given by slippery walls, for which $\bs{v}\cdot\bs{n}=0$, the left boundary is set to a constant inflow flux $h\bs{v}\cdot\bs{n}=-0.02$, and the right boundary condition is a constant outflow with flux $h\bs{v}\cdot\bs{n}=0.02$.
The flow direction is indicated by the arrows in the same figure.
The cylinder location determined by two parameters, namely the radius $R$ and the $x$-coordinate of its center $x_c$.
In all computations, the initial condition is given by the uniform flow condition with $h=0.2$ and $\bs{v}=(0.1,\;0.0)$ at all active nodes.
The first two test cases consist of geometrical parameterization using a one-dimensional parameter space with either $x_c$ or $R$ fixed, respectively.
The third test case consists of a geometrical parameterization with a two-dimensional parameter space where both $x_c$ and $R$ are left free.
All results for the test problems are obtained using a background triangular mesh with 5,419 vertices and 10,476 elements; the average element edge size is $0.02$. 


In all test cases, we generate the snapshots using FOM solutions that are obtained by a fixed CFL number $\alpha_{\cfl}=0.5$.
Denoting the numerical solution at time $t^n$ by $\bs{U}^n$ as before, then the snapshots are picked as $\bs{U}^{kn_{\freq}}$, $k=1,2,\cdots$, where $n_{\freq}$ is the sampling frequency; the solution at the terminal time is always picked as a snapshot, if not sampled already.
To assess the performance of the reduced-order models, we instead compute both the ROM and FOM solutions using a fixed time step size $\Delta t=0.002$, which corresponds to a CFL number slightly smaller than $0.5$, to avoid interpolation error at different time steps.

For each test, we consider two sets of ROM computations: the first set uses the unprocessed POD basis vectors\footnote{That is, a constant value zero is filled at all inactive nodes in the snapshots, see Section~\ref{sec:rom_interp} for more details.} and the second set uses POD basis vectors with interpolated values at inactive nodes.
In all tables and plots, the two sets are designated by the standard SBM-ROM (without interpolation) and with interpolation (SBM-iROM), respectively.
Given a prescribed energy threshold $\mu_{\textrm{pod}}$, the number of POD modes used in the ROM computations is determined by the smallest number $n$ such that:
\begin{displaymath}
  \frac{\sum_{i=1}^n\lambda_i}{\sum_{i=1}^{N_{s}}\lambda_i} \ge \mu_{\textrm{pod}}\;,
\end{displaymath}
where $\lambda_1\ge\lambda_2\ge\cdots\ge\lambda_{N_{s}}$ are all the non-zero eigenvalues of the correlation matrix given in (\ref{eq:pod_cmat}).
As the test is convection dominated, we choose thresholds that are very close to unity, with the particular number of POD modes summarized in Table~\ref{tb:num_cyl_numpod}.
In the first row of the table, we also indicate the sampling frequency $n_{\freq}$ for each test.
\begin{table}[H]\centering
  \caption{Number of POD modes for each ROM computation in all three test cases. ``SBM-iROM'' and ``SBM-ROM'' stand for ROM with and without interpolation, respectively. The row with $\mu_{\textrm{pod}}=1$ corresponds to the sizes of full sets of POD modes.}
  \label{tb:num_cyl_numpod}
  \begin{tabular}{@{}lcccccccc@{}}
    \toprule[.5mm]
    & \multicolumn{2}{l}{Test 1 ($n_{\freq}=2$)} & & \multicolumn{2}{l}{Test 2 ($n_{\freq}=10$)} & & \multicolumn{2}{l}{Test 3 ($n_{\freq}=10$)} \\ \cmidrule[.2mm](lr){2-3} \cmidrule[.2mm](lr){5-6} \cmidrule[.2mm](l){8-9}
    $\mu_{\textrm{pod}}$ & SBM-iROM & SBM-ROM & & SBM-iROM & SBM-ROM & & SBM-iROM & SBM-ROM \\ \cmidrule[.3mm](l){2-9} 
    $1-10^{-5}$ & 15 & 16 & & 63 & 67 & & 83 & 107 \\
    $1-10^{-6}$ & 29 & 30 & & 120 & 127 & & 160 & 200 \\
    $1-10^{-7}$ & 52 & 53 & & 178 & 181 & & 275 & 350 \\
    $1-10^{-8}$ & 71 & 72 & & 223 & 223 & & 465 & 535 \\
    $1-10^{-9}$ & 85 & 86 & & 256 & 256 & & 646 & 679 \\
    $1$ & 561 & 561 & & 344 & 344 & & 1032 & 1032 \\ 
    \bottomrule[.5mm]
  \end{tabular}
\end{table}
The performance of each set of ROM computations is demonstrated both qualitatively and quantitatively:
\begin{itemize}
  \item Qualitatively, the final water height solution is plotted and compared among the FOM computation and two ROM computations, for the latter the number of POD modes is determined by $\mu_{\textrm{pod}}=1-10^{-6}$.
  \item Quantitatively, we compute and tabulate the relative error in Frobenius norm of ROM computations using the POD modes determined using the thresholds $\mu_{\textrm{pod}}=1-10^{-5}$, $1-10^{-6}$, $1-10^{-7}$, $1-10^{-8}$, and $1-10^{-9}$, and compare them to the projected FOM solutions 
  (let the interpolated basis functions be $\bs{\Phi}_r$ as before, the projected solution of the full-order solution $\bs{U}_{\fom}$ is given by $(\bs{\Phi}_r^T\bs{\Phi}_r)^{-1}\bs{\Phi}_r^T\bs{U}_{\fom}$).
  In particular, the Frobenius norm of a generic solution vector $\bs{U}\in\mathbb{R}^{n_{\node}}$ is defined as:
  \begin{equation}\label{eq:num_cyl_frob_active}
    \nrm{\bs{U}}_{\frob(\tilde{\mathcal{D}})} \eqdef \left(\sum_{A=1,\, \bs{x}_A\in\tilde{\mathcal{D}}}^{n_{\node}}U_A^2\right)^{\frac{1}{2}}\;.
  \end{equation}
  To this end, the relative space-time Frobenius error of the ROM solution is computed as:
  \begin{equation}\label{eq:num_cyl_froberr_rom}
    \left(\int_0^T\nrm{\bs{U}_{\fom}-\bs{U}_{\rom}}_{\frob(\tilde{\mathcal{D}})}^2dt\right)^{\frac{1}{2}}\bigg/\left(\int_0^T\nrm{\bs{U}_{\fom}}_{\frob(\tilde{\mathcal{D}})}^2dt\right)^{\frac{1}{2}}\;;
  \end{equation}
  whereas the relative space-time Frobenius error of the projected FOM solution is given by:
  \begin{equation}\label{eq:num_cyl_froberr_fomproj}
    \left(\int_0^T\nrm{\bs{U}_{\fom}-(\bs{\Phi}_r^T\bs{\Phi}_r)^{-1}\bs{\Phi}_r^T\bs{U}_{\fom}}_{\frob(\tilde{\mathcal{D}})}^2dt\right)^{\frac{1}{2}}\bigg/\left(\int_0^T\nrm{{\bs{U}}_{\fom}}_{\frob(\tilde{\mathcal{D}})}^2dt\right)^{\frac{1}{2}}\;.
  \end{equation}
  In both~(\ref{eq:num_cyl_froberr_rom}) and~(\ref{eq:num_cyl_froberr_fomproj}), the time integral is approximated by weighted sum of discrete solutions in a straight forward manner. 
\end{itemize}

\subsection{Test 1: Geometrical parameterization with varying cylinder radius}
\label{sec:num_cyl_r}
In this test case we fix $x_c=0.0$ and generate the FOM snapshots using three radii $R=0.1$, $R=0.15$, and $R=0.2$; in all tests the computation is performed until $T=0.8$.
A total number of 561 snapshots are created from the three FOM computations with a sampling frequency $n_{\freq}=2$, and we assess the performance of ROM by computing the flow past cylinders with radii $R=0.08$, $R=0.13$, $R=0.17$, and $R=0.22$.

\begin{figure}[h]\centering
  \includegraphics[width=.48\textwidth]{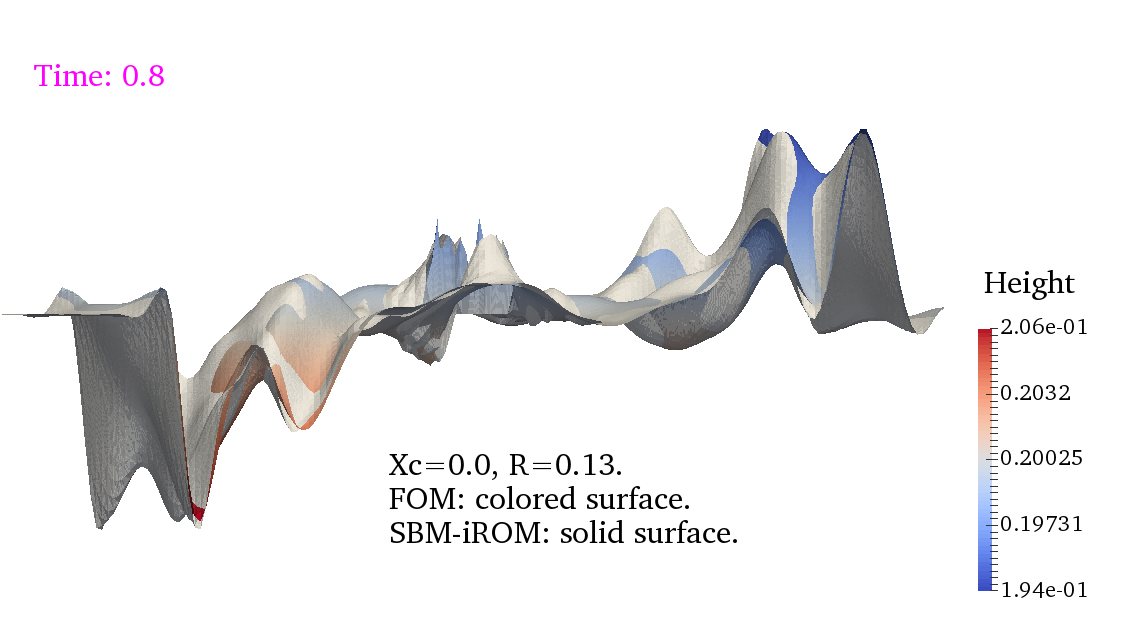}
  \includegraphics[width=.48\textwidth]{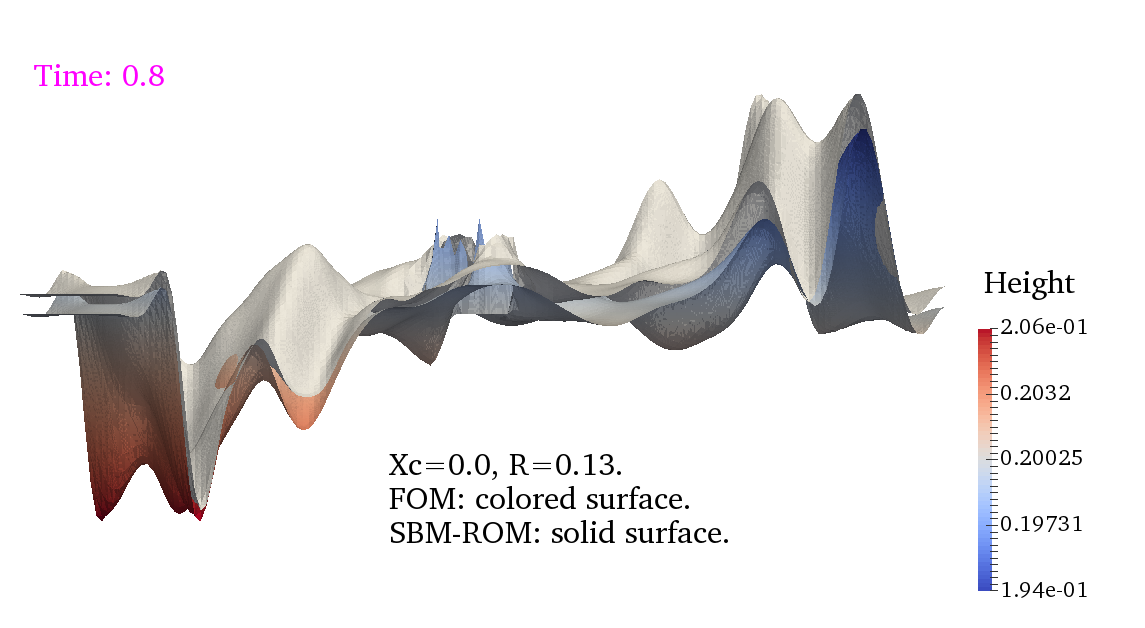}
  \caption{The solution surfaces for water height at $T=0.8$ in the case $R=0.13$ of Test 1 for SBM-iROM (left panel) and SBM-ROM (right panel). In both plots, the ROM solution (solid surface) is plotted on top of the FOM one (color surface).}
  \label{fg:num_cyl_r_warp_r0d13}
\end{figure}

Among the four, $R=0.08$ is the most challenging one in the sense that several active nodes are inactive in all snapshots, hence it is not surprising to see that all ROM computations without interpolation fail to deliver any reasonable solutions.
For the other three radii tested, all active nodes are also active in some snapshots and thus a solution is obtained whether the interpolation is employed or not; in these cases, we still observe that no interpolation leads to significantly worse ROM solution, as demonstrated by the height solution surfaces for the case $R=0.13$ in Figure~\ref{fg:num_cyl_r_warp_r0d13}.
In particular, SBM-ROM leads to an overall ``shift'' in the entire computational domain whereas the SBM-iROM solution agrees much better with the FOM one.
It worthies noting here that in other two tests (that is when $x_c$ is varying), the SBM-ROM solutions demonstrate oscillation with very large magnitude; hence one expects no significant improvement for ROM without interpolation, even if some non-zero value is used to fill all the inactive nodes in the snapshots.

\begin{figure}[h]\centering
  \includegraphics[trim=0.0in 3.0in 0.0in 0.0in, clip, width=.48\textwidth]{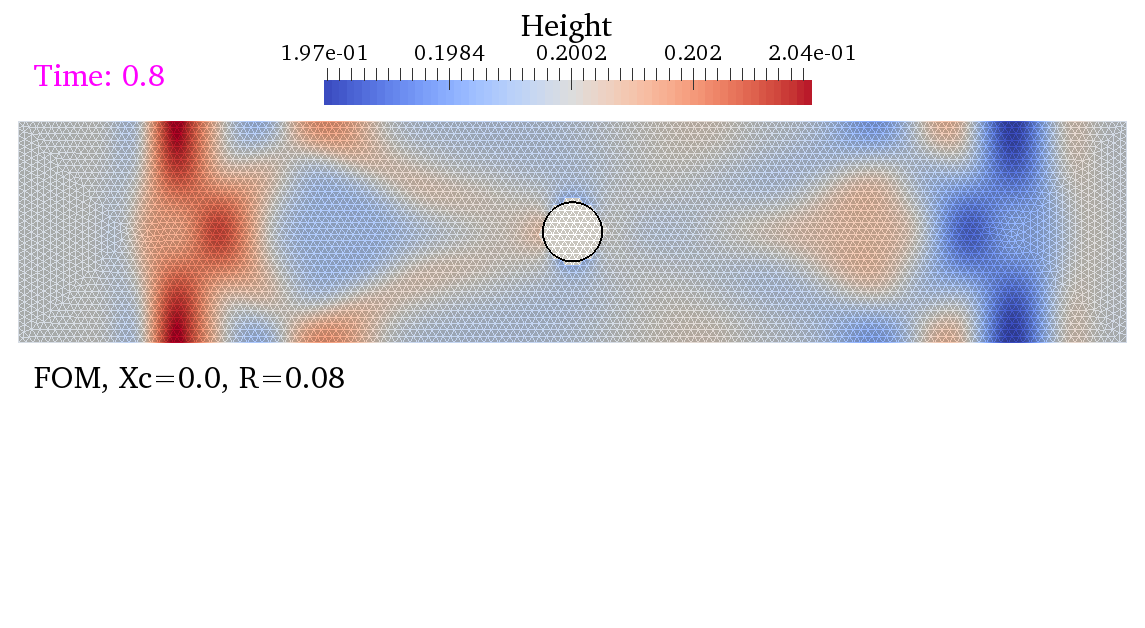} 
  \includegraphics[trim=0.0in 3.0in 0.0in 0.0in, clip, width=.48\textwidth]{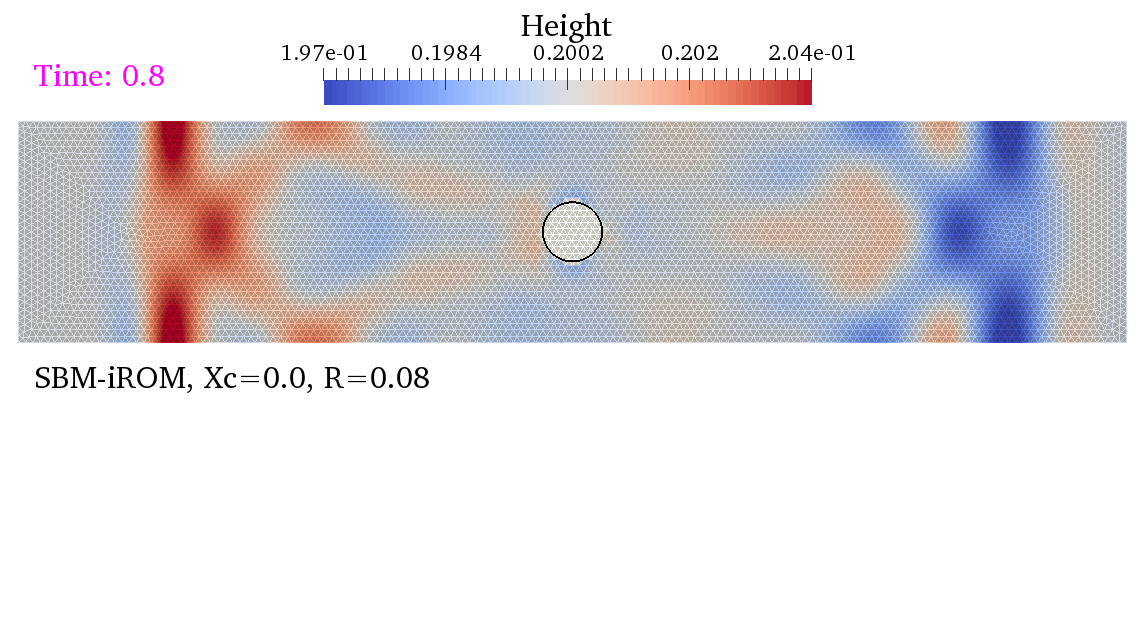} \\
  \includegraphics[trim=0.0in 3.0in 0.0in 0.0in, clip, width=.48\textwidth]{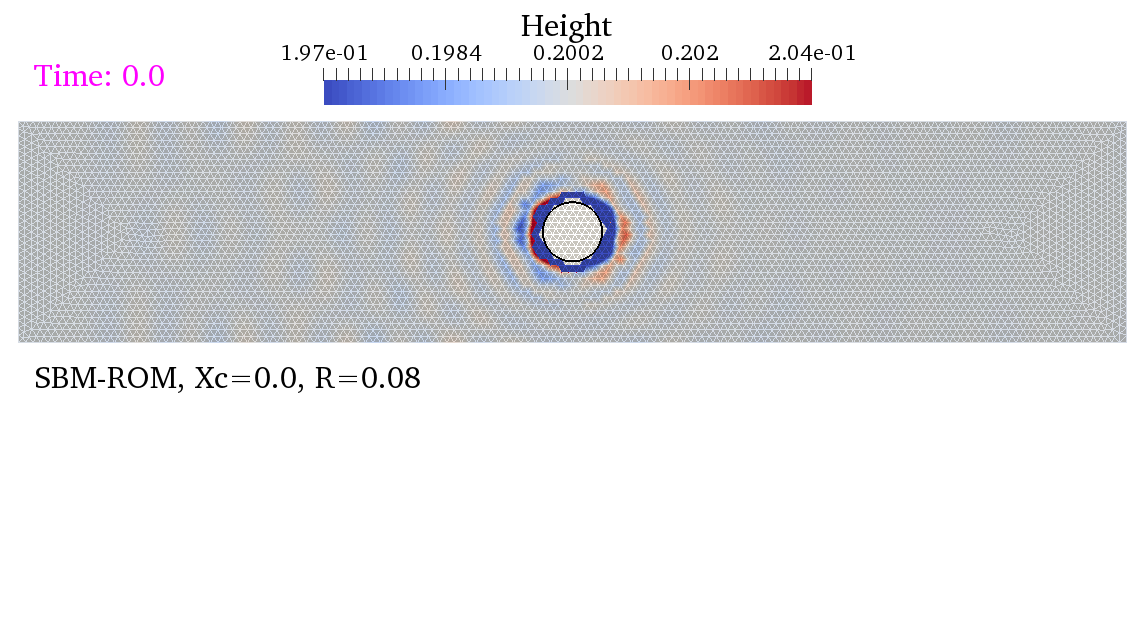} 
  \includegraphics[trim=0.0in 3.0in 0.0in 0.0in, clip, width=.48\textwidth]{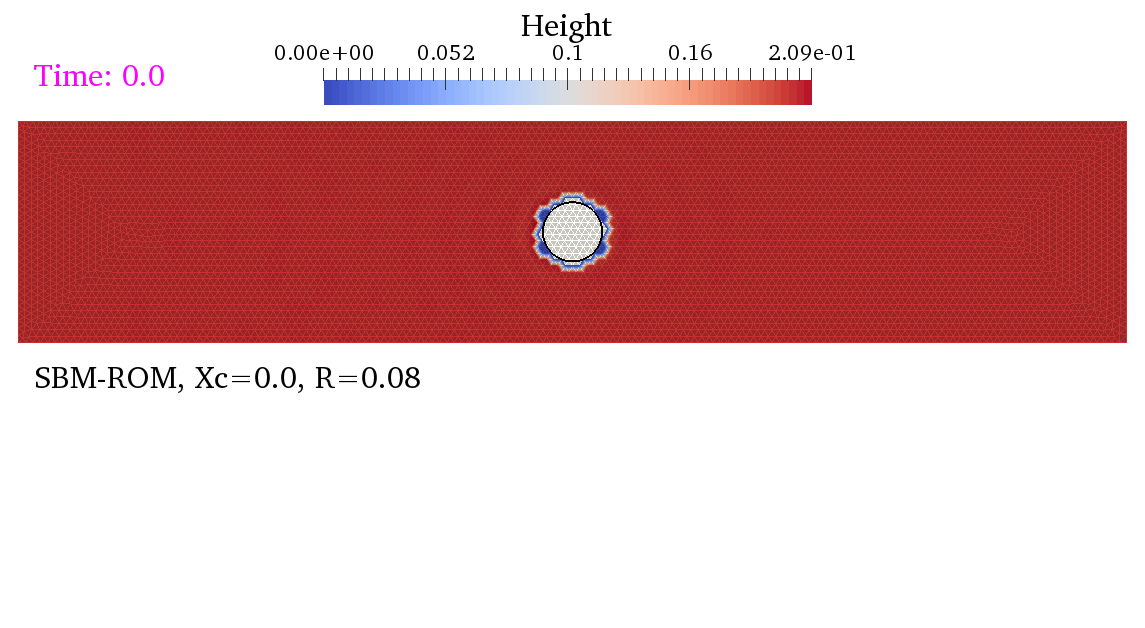} 
  \caption{The height ($h$) in the case of $R=0.08$ computed by FOM (top left), SBM-iROM (top right), and initial height of unsuccessful ROM computations without interpolation (bottom row). 
  The legend range is set according to the FOM computation except the bottom right panel, where the SBM-ROM solution is plotted using its own range.
  All ROM computations are performed using POD modes corresponding to $\mu_{\textrm{pod}}=1-10^{-6}$.}
  \label{fg:num_cyl_r0d08_sol}
\end{figure}

\begin{figure}[h]\centering
  \includegraphics[trim=0.0in 3.0in 0.0in 0.0in, clip, width=.48\textwidth]{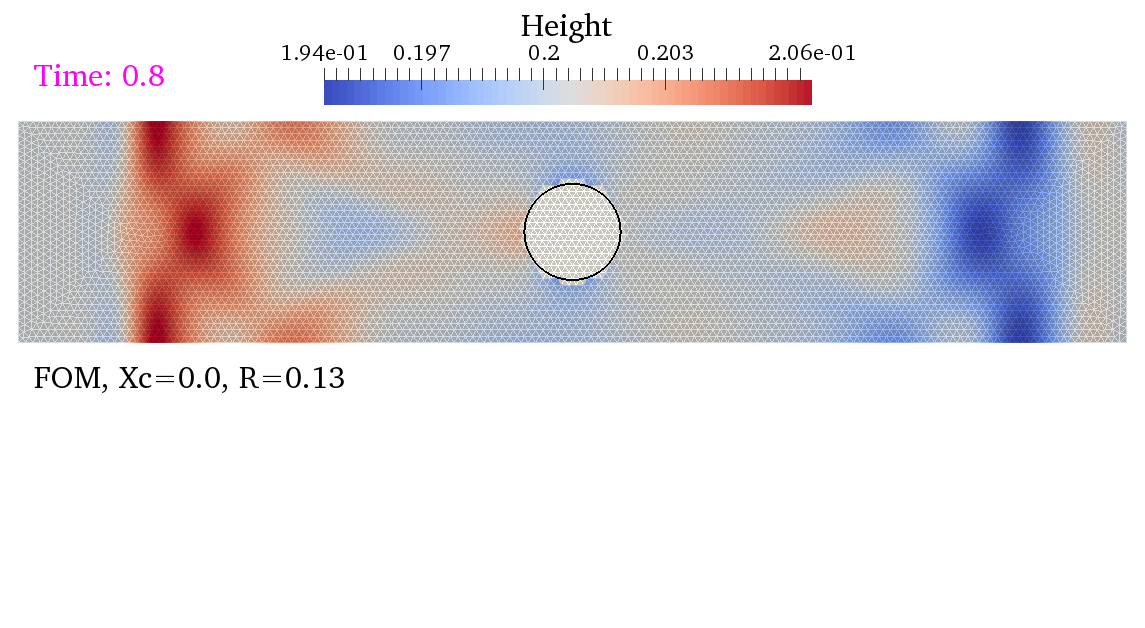} 
  \includegraphics[trim=0.0in 3.0in 0.0in 0.0in, clip, width=.48\textwidth]{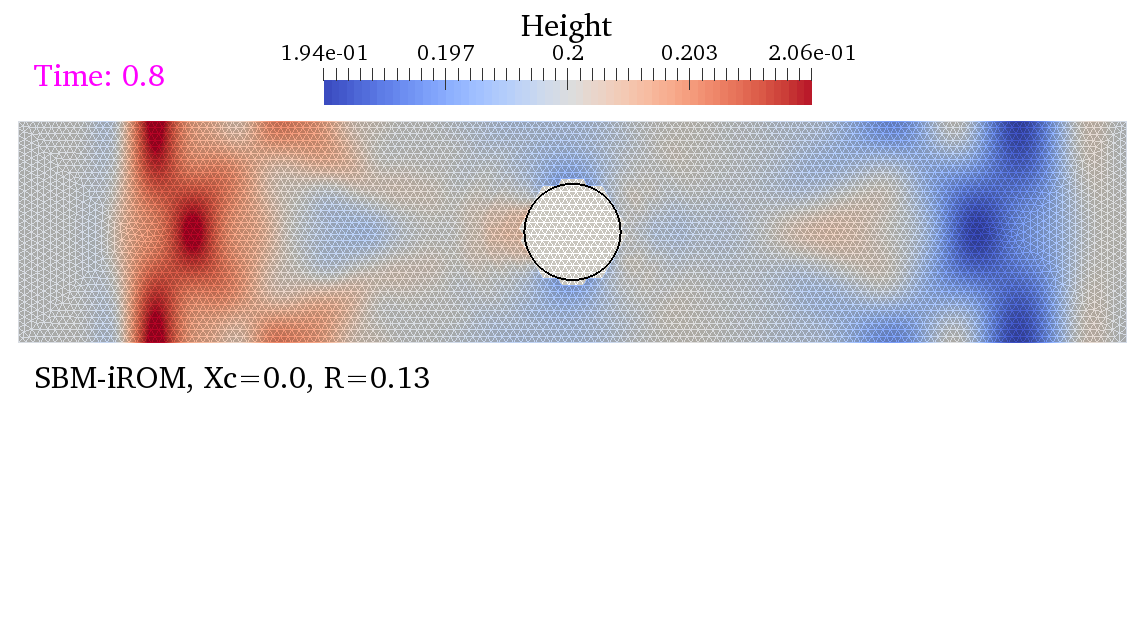} \\
  \includegraphics[trim=0.0in 3.0in 0.0in 0.0in, clip, width=.48\textwidth]{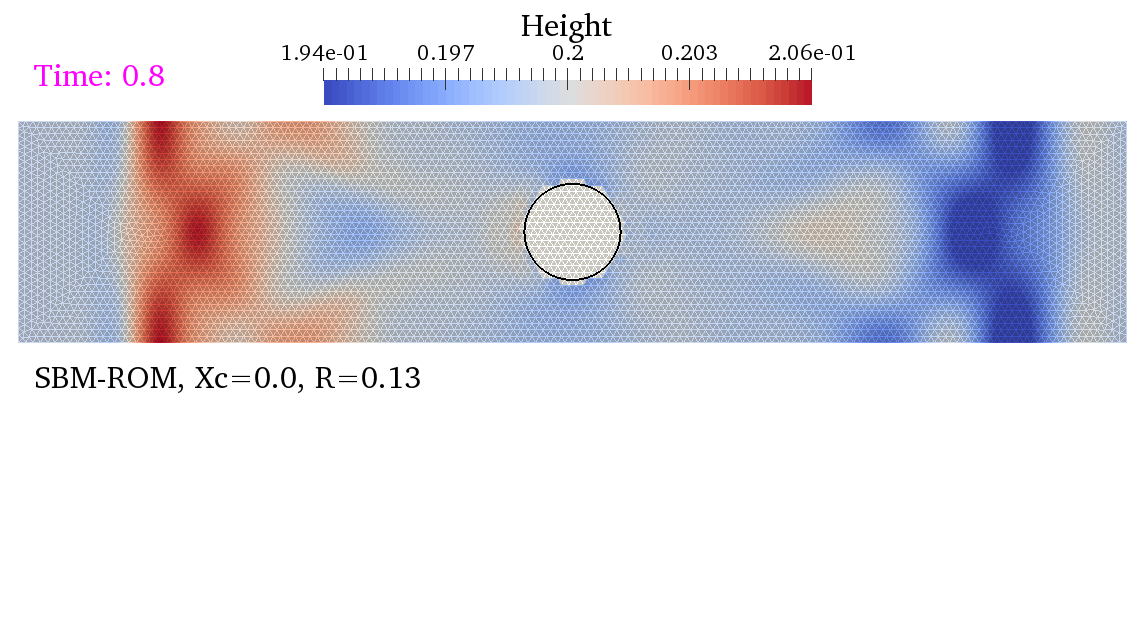} 
  \includegraphics[trim=0.0in 3.0in 0.0in 0.0in, clip, width=.48\textwidth]{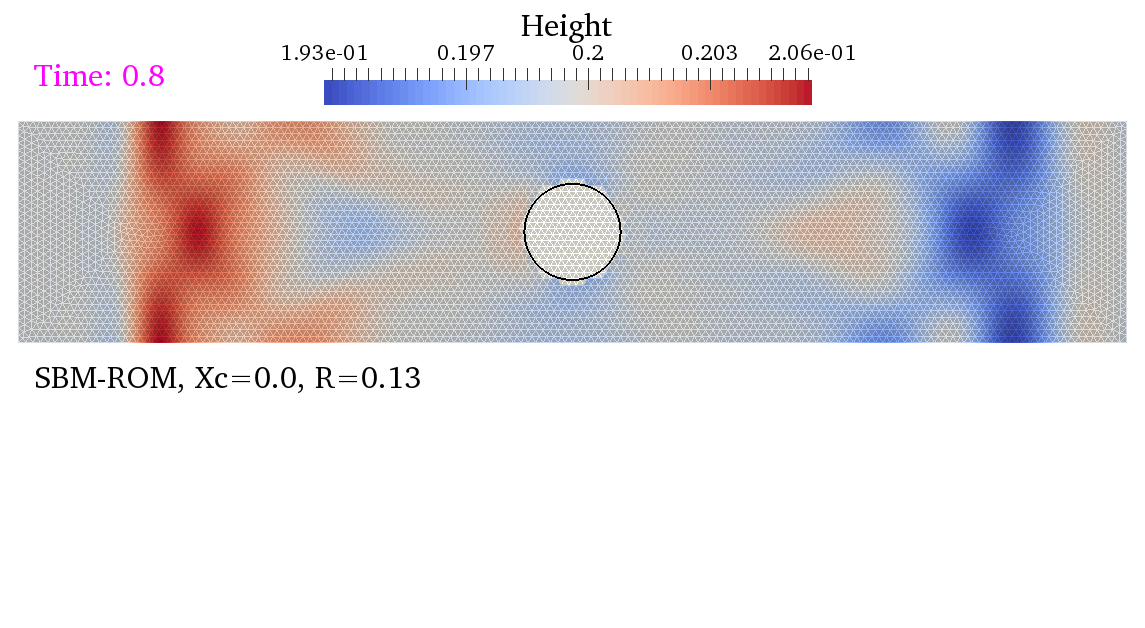} 
  \caption{The height ($h$) in the case of $R=0.13$ computed by FOM (top left), SBM-iROM (top right), and SBM-ROM (bottom row). 
  The legend range is set according to the FOM computation except the bottom right panel, where the SBM-ROM solution is plotted using its own range.
  All ROM computations are performed using POD modes corresponding to $\mu_{\textrm{pod}}=1-10^{-6}$.}
  \label{fg:num_cyl_r0d13_sol}
\end{figure}

\begin{figure}[h]\centering
  \includegraphics[trim=0.0in 3.0in 0.0in 0.0in, clip, width=.48\textwidth]{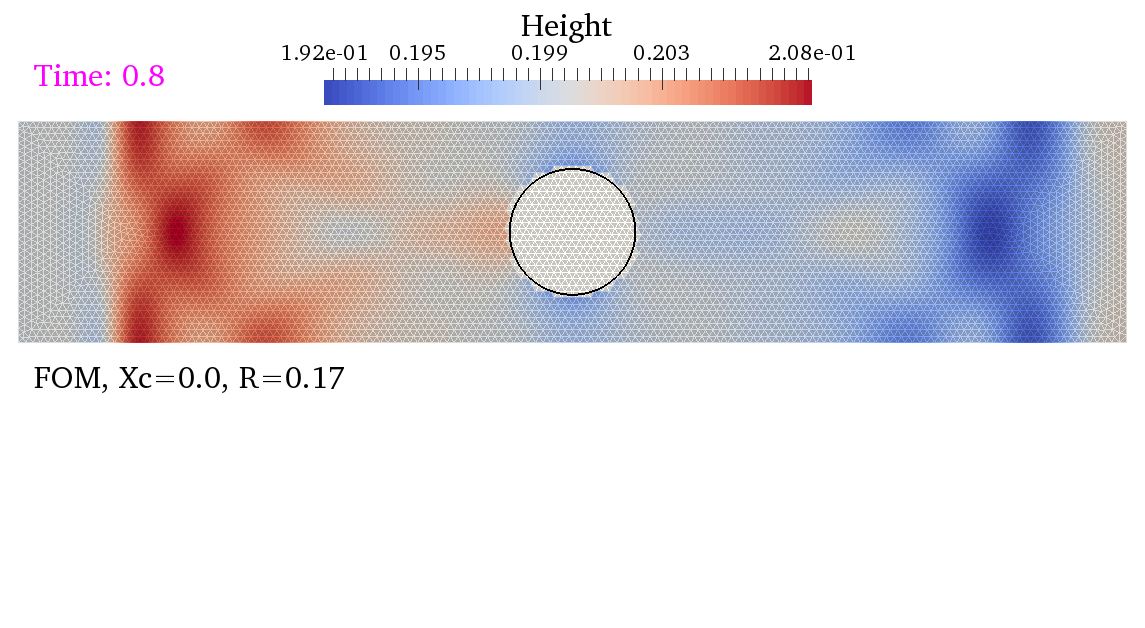} 
  \includegraphics[trim=0.0in 3.0in 0.0in 0.0in, clip, width=.48\textwidth]{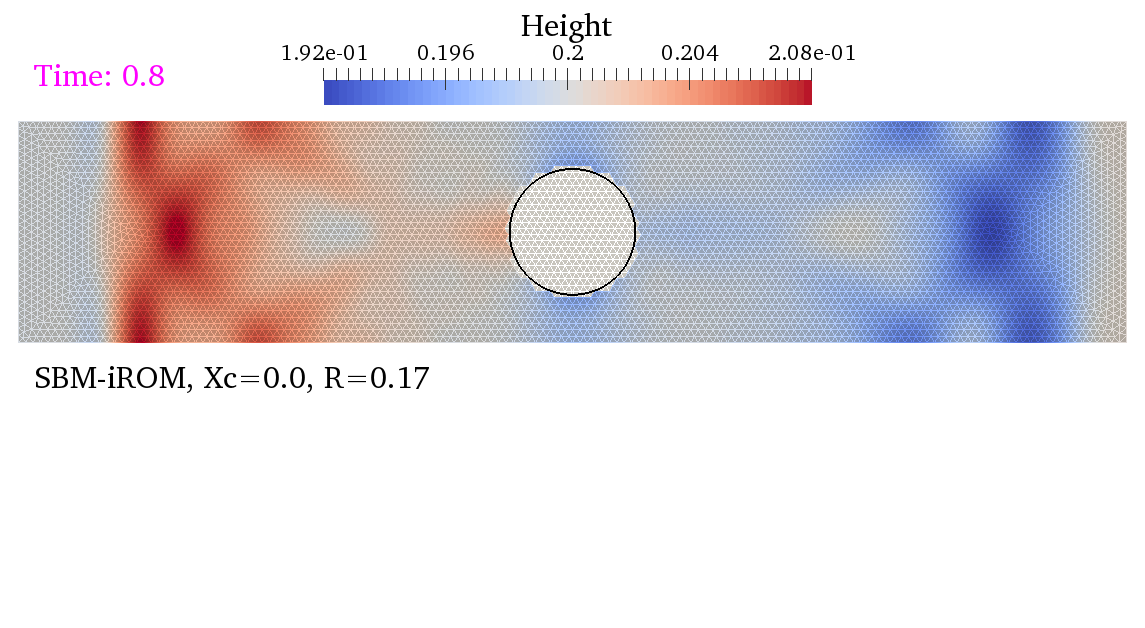} \\
  \includegraphics[trim=0.0in 3.0in 0.0in 0.0in, clip, width=.48\textwidth]{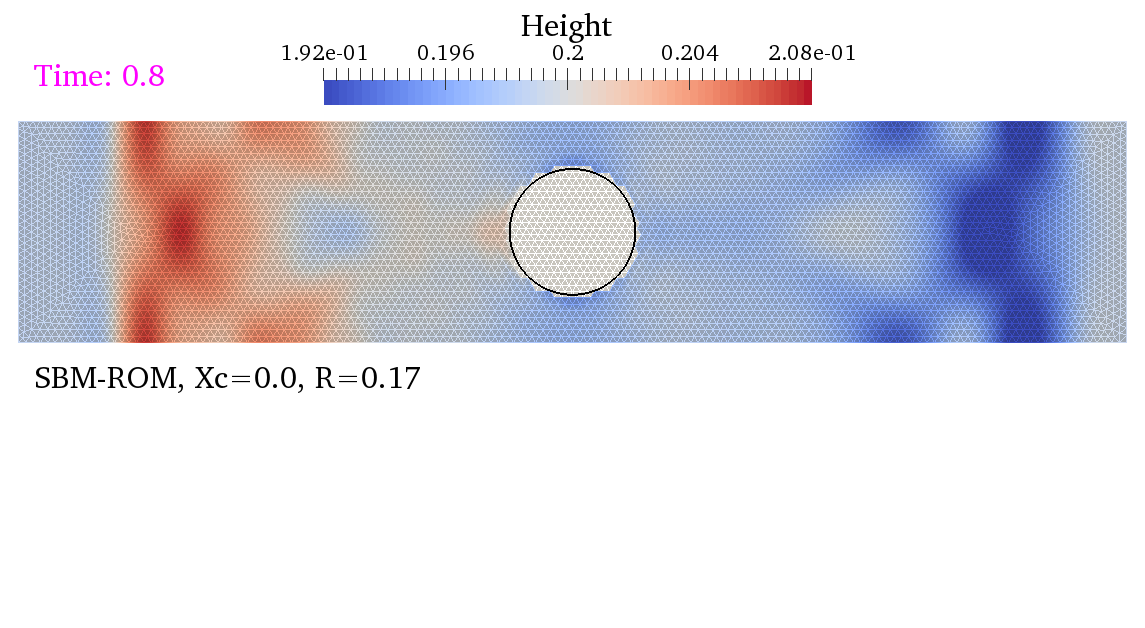} 
  \includegraphics[trim=0.0in 3.0in 0.0in 0.0in, clip, width=.48\textwidth]{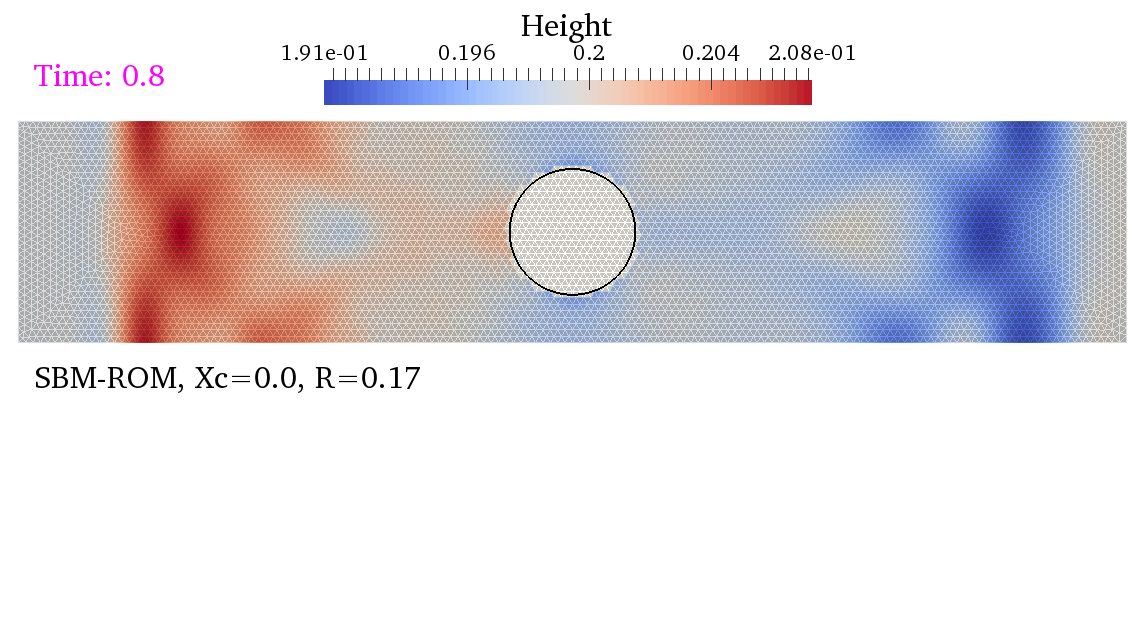} 
  \caption{The height ($h$) in the case of $R=0.17$ computed by FOM (top left), SBM-iROM (top right), and SBM-ROM (bottom row). 
  The legend range is set according to the FOM computation except the bottom right panel, where the SBM-ROM solution is plotted using its own range.
  All ROM computations are performed using POD modes corresponding to $\mu_{\textrm{pod}}=1-10^{-6}$.}
  \label{fg:num_cyl_r0d17_sol}
\end{figure}

\begin{figure}[h]\centering
  \includegraphics[trim=0.0in 3.0in 0.0in 0.0in, clip, width=.48\textwidth]{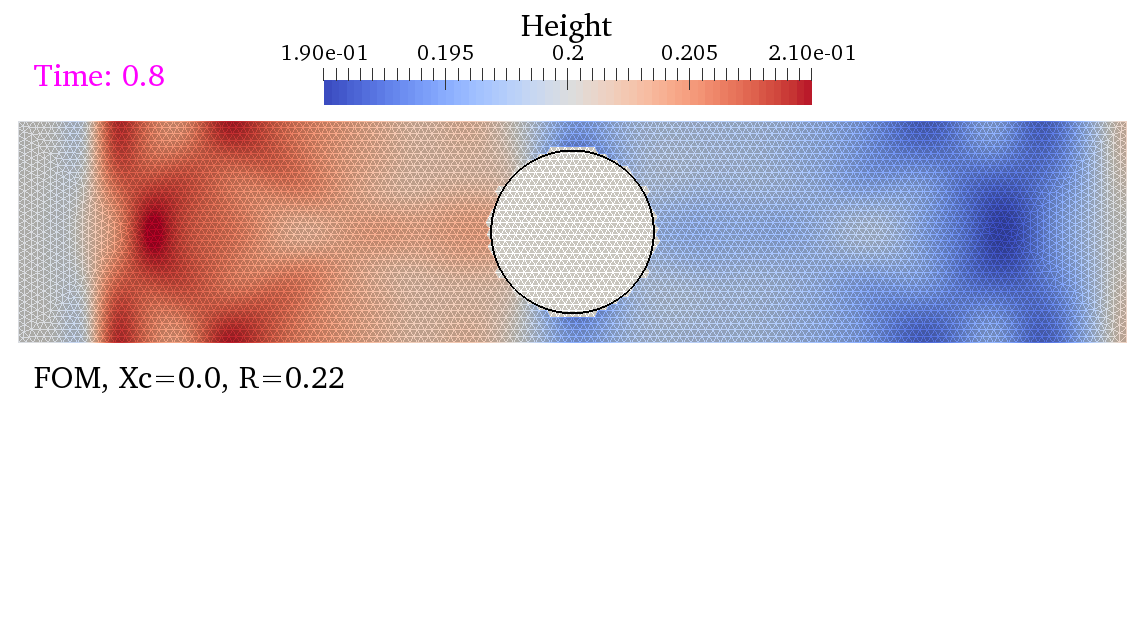} 
  \includegraphics[trim=0.0in 3.0in 0.0in 0.0in, clip, width=.48\textwidth]{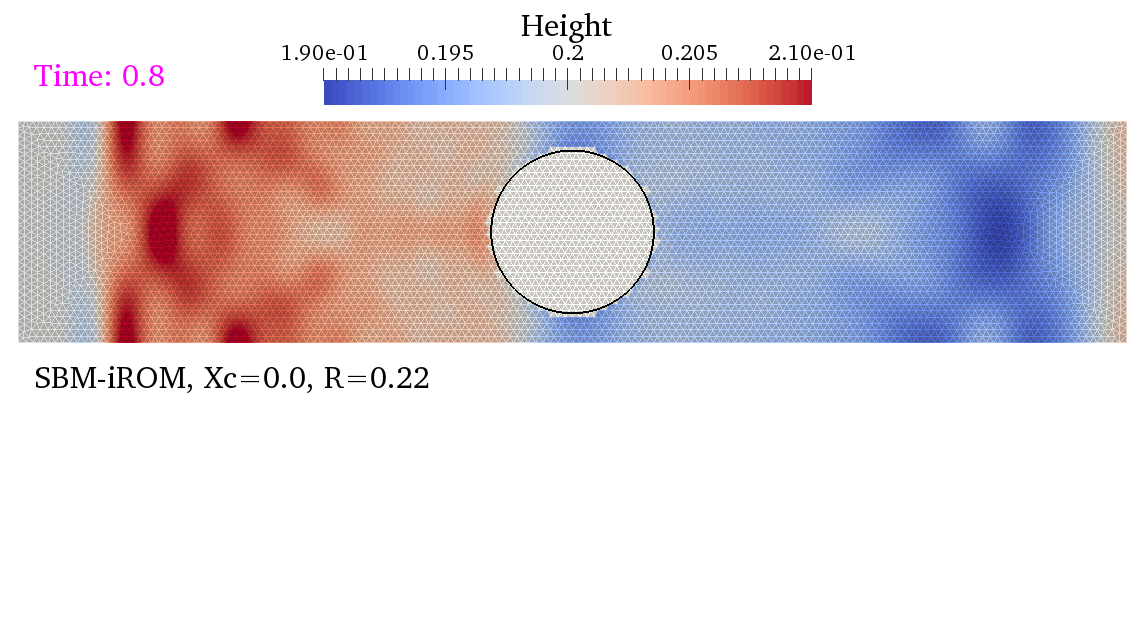} \\
  \includegraphics[trim=0.0in 3.0in 0.0in 0.0in, clip, width=.48\textwidth]{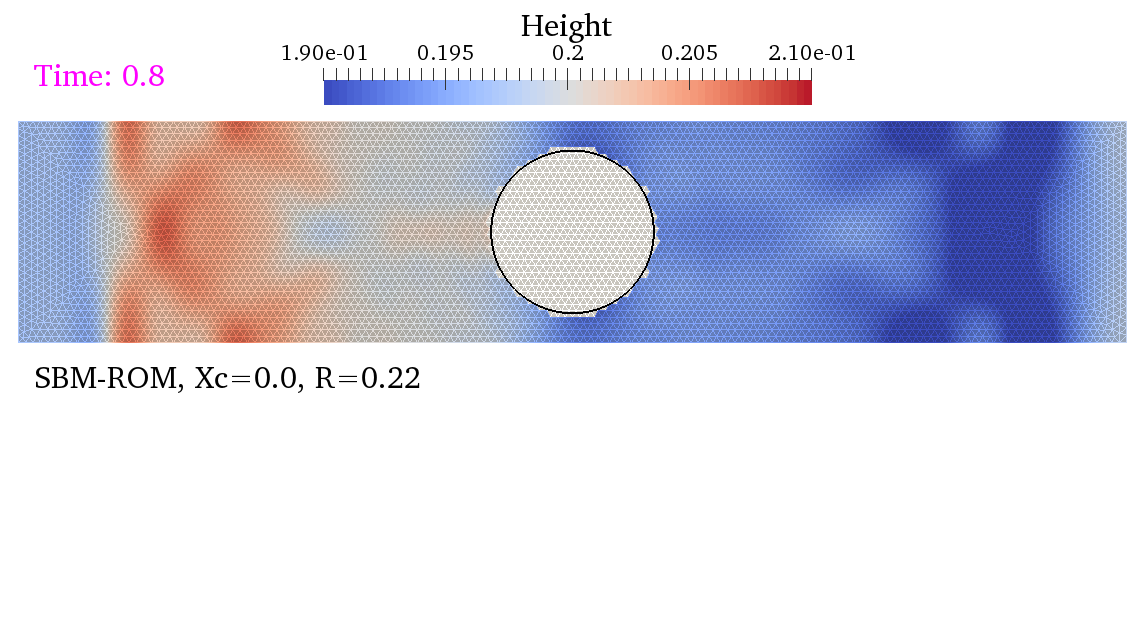} 
  \includegraphics[trim=0.0in 3.0in 0.0in 0.0in, clip, width=.48\textwidth]{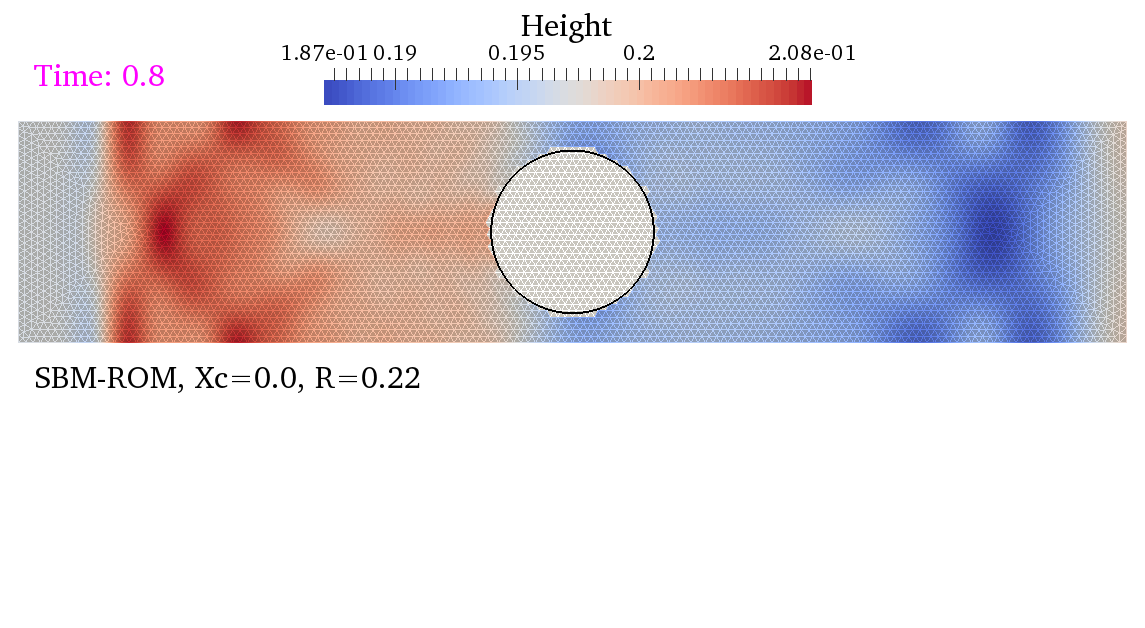} 
  \caption{The height ($h$) in the case of $R=0.22$ computed by FOM (top left), SBM-iROM (top right), and SBM-ROM (bottom row). 
  The legend range is set according to the FOM computation except the bottom right panel, where the SBM-ROM solution is plotted using its own range.
  All ROM computations are performed using POD modes corresponding to $\mu_{\textrm{pod}}=1-10^{-6}$.}
  \label{fg:num_cyl_r0d22_sol}
\end{figure}

Next in Figures~\ref{fg:num_cyl_r0d08_sol}--\ref{fg:num_cyl_r0d22_sol}, we plot the water height solutions at terminal time $T=0.8$ in all four cases.
As mentioned at the beginning of this section, the ROM computations are performed using POD modes corresponding to the energy threshold $\mu_{\textrm{pod}}=1-10^{-6}$, that is, 29 modes for SBM-iROM and 30 modes for SBM-ROM, see also Table~\ref{tb:num_cyl_numpod}.
Note that when $R=0.08$, SBM-ROM is unsuccessful and thus only the initial data is plotted (see the bottom row of Figure~\ref{fg:num_cyl_r0d08_sol}).
Furthermore, because the SBM-ROM solutions are very different from the FOM one, we plot them both in the scale of the FOM solution and in the scale of its own in the bottom left panel and the bottom right panel, respectively, in these figures. 

\begin{figure}[h]\centering
  \begin{subfigure}[t]{.4\textwidth}\centering
    \includegraphics[width=\textwidth]{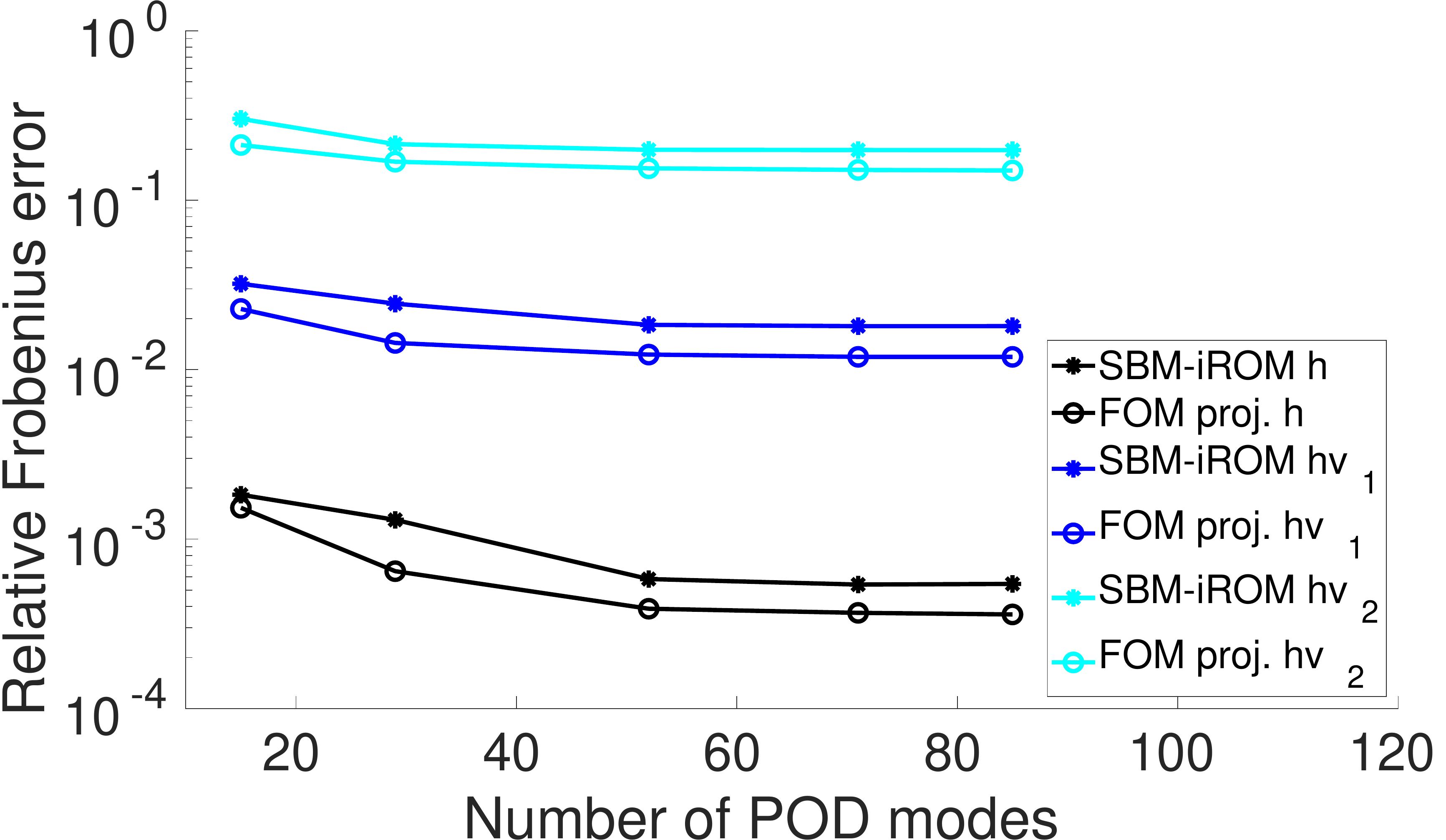}
    \caption{$R=0.08$. No solution is computed by uROM.}
    \label{fg:num_cyl_r_froberr_r0d08}
  \end{subfigure}
  \begin{subfigure}[t]{.4\textwidth}\centering
    \includegraphics[width=\textwidth]{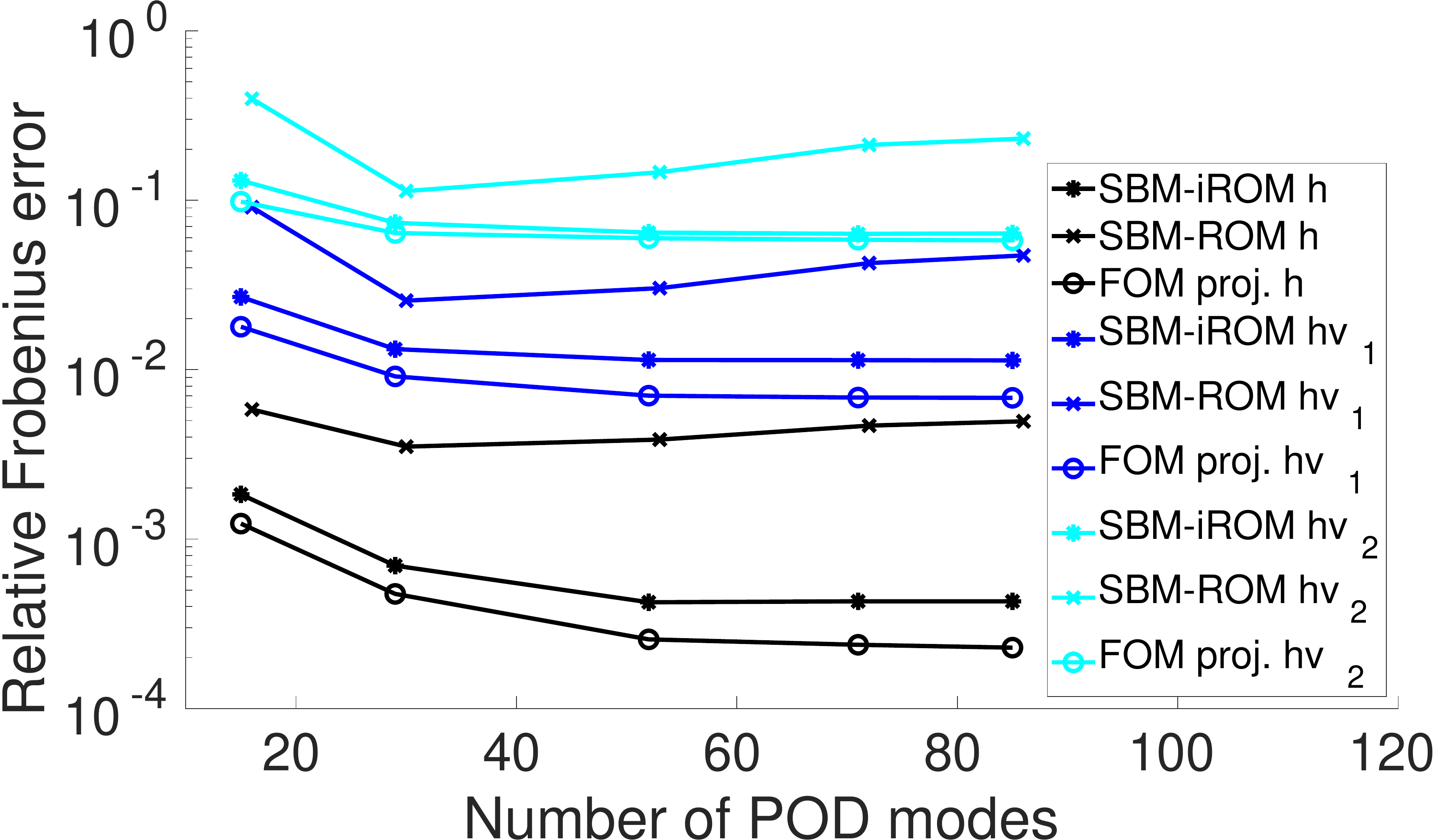}
    \caption{$R=0.13$.}
    \label{fg:num_cyl_r_froberr_r0d13}
  \end{subfigure} \\
  \begin{subfigure}[t]{.4\textwidth}\centering
    \includegraphics[width=\textwidth]{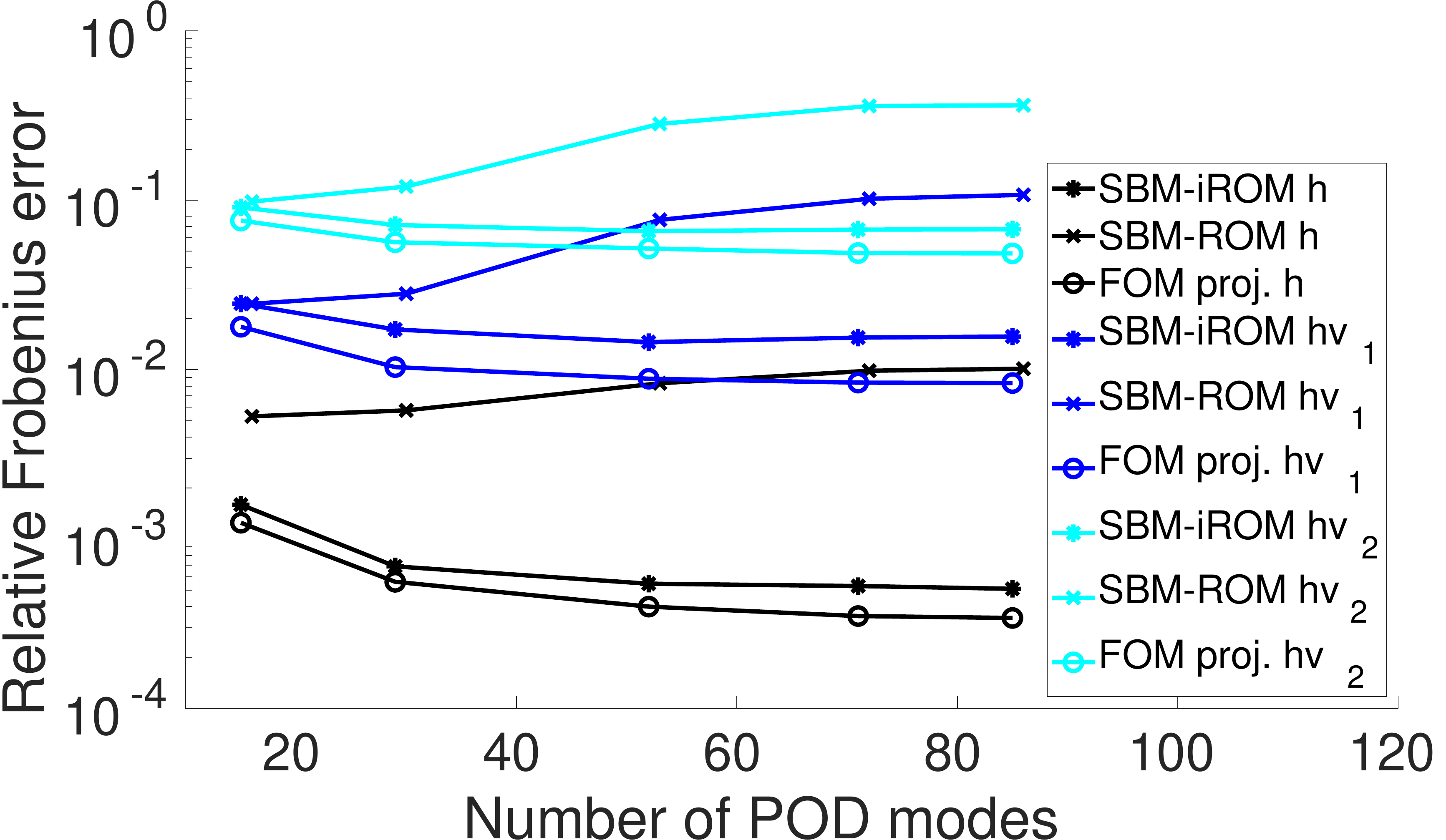}
    \caption{$R=0.17$.}
    \label{fg:num_cyl_r_froberr_r0d17}
  \end{subfigure}
  \begin{subfigure}[t]{.4\textwidth}\centering
    \includegraphics[width=\textwidth]{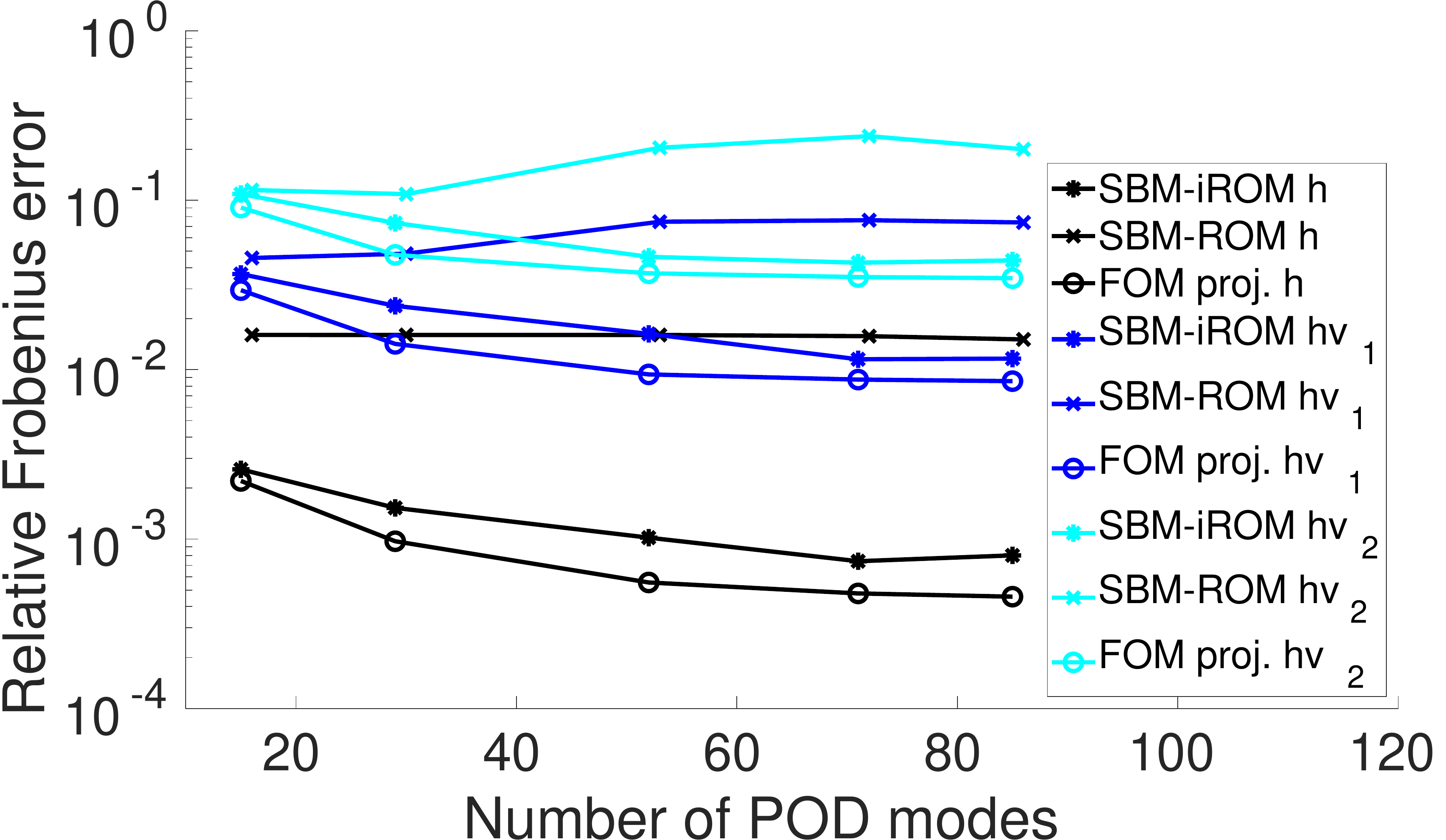}
    \caption{$R=0.22$.}
    \label{fg:num_cyl_r_froberr_r0d22}
  \end{subfigure}
  \caption{The relative Frobenius errors of Test 1 by SBM-iROM (marked by $\ast$), SBM-ROM (marked by x), and FOM projection error (FOM proj., marked by o) in logarithmic scale plotted against the number of POD modes.
  The errors in $h$, $hv_1$, $hv_2$ are plotted in black, blue, and cyan colors, respectively.}
  \label{fg:num_cyl_r_froberr}
\end{figure}

The relative space-time Frobenius errors of ROM computations using different number of POD modes and the projected FOM solutions are summarized in Table~\ref{tb:num_cyl_r_froberr_r0d08}--Table~\ref{tb:num_cyl_r_froberr_r0d22} for the four radii, respectively.
\begin{table}[h]\centering
  \caption{The relative space-time Frobenius errors of the projected FOM solutions and ROM computations in the case $R=0.08$ of Test 1. Note that SBM-ROM fails to deliver any solution.}
  \label{tb:num_cyl_r_froberr_r0d08}
  \begin{tabular}{@{}lccccccccccc@{}}
  \toprule[.5mm]
  & \multicolumn{3}{l}{FOM projection} & & \multicolumn{3}{l}{SBM-iROM} & & \multicolumn{3}{l}{SBM-ROM} \\ \cmidrule[.2mm](lr){2-4} \cmidrule[.2mm](lr){6-8} \cmidrule[.2mm](l){10-12}
  $\mu_{\textrm{pod}}$ & $h$ & $hv_1$ & $hv_2$ & & $h$ & $hv_1$ & $hv_2$ & & $h$ & $hv_1$ & $hv_2$ \\ \cmidrule[.3mm](l){2-12}
  $1-10^{-5}$ & 1.53e-3 & 2.28e-2 & 2.12e-1 & & 1.83e-3 & 3.21e-2 & 3.02e-1 & & -- & -- & -- \\
  $1-10^{-6}$ & 6.46e-4 & 1.44e-2 & 1.69e-1 & & 1.30e-3 & 2.45e-2 & 2.14e-1 & & -- & -- & -- \\
  $1-10^{-7}$ & 3.88e-4 & 1.23e-2 & 1.54e-1 & & 5.81e-4 & 1.84e-2 & 1.98e-1 & & -- & -- & -- \\
  $1-10^{-8}$ & 3.67e-4 & 1.19e-2 & 1.51e-1 & & 5.39e-4 & 1.81e-2 & 1.98e-1 & & -- & -- & -- \\
  $1-10^{-9}$ & 3.59e-4 & 1.19e-2 & 1.50e-1 & & 5.44e-4 & 1.81e-2 & 1.97e-1 & & -- & -- & -- \\
  \bottomrule[.4mm]
  \end{tabular}
\end{table}
\begin{table}[h]\centering
  \caption{The relative space-time Frobenius errors of the projected FOM solutions and ROM computations in the case $R=0.13$ of Test 1.}
  \label{tb:num_cyl_r_froberr_r0d13}
  \begin{tabular}{@{}lccccccccccc@{}}
  \toprule[.5mm]
  & \multicolumn{3}{l}{FOM projection} & & \multicolumn{3}{l}{SBM-iROM} & & \multicolumn{3}{l}{SBM-ROM} \\ \cmidrule[.2mm](lr){2-4} \cmidrule[.2mm](lr){6-8} \cmidrule[.2mm](l){10-12}
  $\mu_{\textrm{pod}}$ & $h$ & $hv_1$ & $hv_2$ & & $h$ & $hv_1$ & $hv_2$ & & $h$ & $hv_1$ & $hv_2$ \\ \cmidrule[.3mm](l){2-12}
  $1-10^{-5}$ & 1.23e-3 & 1.79e-2 & 9.81e-2 & & 1.84e-3 & 2.69e-2 & 1.31e-1 & & 5.82e-3 & 9.09e-2 & 3.97e-1 \\
  $1-10^{-6}$ & 4.75e-4 & 9.12e-3 & 6.40e-2 & & 6.97e-4 & 1.32e-2 & 7.35e-2 & & 3.51e-3 & 2.55e-2 & 1.13e-1 \\
  $1-10^{-7}$ & 2.56e-4 & 7.02e-3 & 5.97e-2 & & 4.24e-4 & 1.14e-2 & 6.44e-2 & & 3.87e-3 & 3.03e-2 & 1.46e-1 \\
  $1-10^{-8}$ & 2.38e-4 & 6.84e-3 & 5.85e-2 & & 4.29e-4 & 1.14e-2 & 6.33e-2 & & 4.67e-3 & 4.25e-2 & 2.12e-1 \\
  $1-10^{-9}$ & 2.29e-4 & 6.82e-3 & 5.80e-2 & & 4.29e-4 & 1.13e-2 & 6.37e-2 & & 4.96e-3 & 4.72e-2 & 2.31e-1 \\
  \bottomrule[.4mm]
  \end{tabular}
\end{table}
\begin{table}[h]\centering
  \caption{The relative space-time Frobenius errors of the projected FOM solutions and ROM computations in the case $R=0.17$ of Test 1.}
  \label{tb:num_cyl_r_froberr_r0d17}
  \begin{tabular}{@{}lccccccccccc@{}}
  \toprule[.5mm]
  & \multicolumn{3}{l}{FOM projection} & & \multicolumn{3}{l}{SBM-iROM} & & \multicolumn{3}{l}{SBM-ROM} \\ \cmidrule[.2mm](lr){2-4} \cmidrule[.2mm](lr){6-8} \cmidrule[.2mm](l){10-12}
  $\mu_{\textrm{pod}}$ & $h$ & $hv_1$ & $hv_2$ & & $h$ & $hv_1$ & $hv_2$ & & $h$ & $hv_1$ & $hv_2$ \\ \cmidrule[.3mm](l){2-12}
  $1-10^{-5}$ & 1.25e-3 & 1.79e-2 & 7.58e-2 & & 1.60e-3 & 2.45e-2 & 9.07e-2 & & 5.30e-3 & 2.45e-2 & 9.81e-2 \\
  $1-10^{-6}$ & 5.60e-4 & 1.04e-2 & 5.65e-2 & & 6.90e-4 & 1.73e-2 & 7.14e-2 & & 5.75e-3 & 2.81e-2 & 1.20e-1 \\
  $1-10^{-7}$ & 3.99e-4 & 8.84e-3 & 5.19e-2 & & 5.45e-4 & 1.45e-2 & 6.57e-2 & & 8.31e-3 & 7.65e-2 & 2.82e-1 \\
  $1-10^{-8}$ & 3.51e-4 & 8.38e-3 & 4.87e-2 & & 5.28e-4 & 1.55e-2 & 6.70e-2 & & 9.85e-3 & 1.02e-1 & 3.60e-1 \\
  $1-10^{-9}$ & 3.42e-4 & 8.34e-3 & 4.85e-2 & & 5.09e-4 & 1.57e-2 & 6.73e-2 & & 1.01e-2 & 1.08e-1 & 3.63e-1 \\
  \bottomrule[.4mm]
  \end{tabular}
\end{table}
\begin{table}[h]\centering
  \caption{The relative space-time Frobenius errors of the projected FOM solutions and ROM computations in the case $R=0.22$ of Test 1.}
  \label{tb:num_cyl_r_froberr_r0d22}
  \begin{tabular}{@{}lccccccccccc@{}}
  \toprule[.5mm]
  & \multicolumn{3}{l}{FOM projection} & & \multicolumn{3}{l}{SBM-iROM} & & \multicolumn{3}{l}{SBM-ROM} \\ \cmidrule[.2mm](lr){2-4} \cmidrule[.2mm](lr){6-8} \cmidrule[.2mm](l){10-12}
  $\mu_{\textrm{pod}}$ & $h$ & $hv_1$ & $hv_2$ & & $h$ & $hv_1$ & $hv_2$ & & $h$ & $hv_1$ & $hv_2$ \\ \cmidrule[.3mm](l){2-12}
  $1-10^{-5}$ & 2.20e-3 & 2.95e-2 & 9.05e-2 & & 2.58e-3 & 3.67e-2 & 1.09e-1 & & 1.60e-2 & 4.56e-2 & 1.15e-1 \\
  $1-10^{-6}$ & 9.71e-4 & 1.42e-2 & 4.76e-2 & & 1.53e-3 & 2.38e-2 & 7.34e-2 & & 1.60e-2 & 4.84e-2 & 1.09e-1 \\
  $1-10^{-7}$ & 5.55e-4 & 9.37e-3 & 3.70e-2 & & 1.02e-3 & 1.63e-2 & 4.63e-2 & & 1.60e-2 & 7.46e-2 & 2.04e-1 \\
  $1-10^{-8}$ & 4.78e-4 & 8.74e-3 & 3.52e-2 & & 7.40e-4 & 1.15e-2 & 4.28e-2 & & 1.58e-2 & 7.62e-2 & 2.39e-1 \\
  $1-10^{-9}$ & 4.57e-4 & 8.55e-3 & 3.47e-2 & & 8.02e-4 & 1.16e-2 & 4.41e-2 & & 1.51e-2 & 7.39e-2 & 2.00e-1 \\
  \bottomrule[.4mm]
  \end{tabular}
\end{table}
These errors are also plotted against the number of POD modes in Figure~\ref{fg:num_cyl_r_froberr}.
The correspondence between the number of POD modes and $\mu_{\textrm{pod}}$ is documented in Table~\ref{tb:num_cyl_numpod}.
From these figures we see that the errors by SBM-iROM are very close to the projection error by the same set of POD modes, and these errors typically decrease as the number of modes increases.
In contrast, without interpolation on the one hand the errors by the ROM solutions are far from the projection error, and they seem to increase as a larger number of POD modes is used.

\subsection{Test 2: Geometrical parameterization with varying cylinder center location}
\label{sec:num_cyl_x}
In the second test case we fix $R=0.15$ and generate the FOM snapshots using three cylinder locations with $x_c=-0.5$, $x_c=0.0$, and $x_c=0.5$, with termination times being $T=0.5$, $T=0.8$, and $T=0.5$, respectively.
The shorter simulation period is selected to avoid interaction between the reflected waves and the left or right boundaries.
A total number of 344 snapshots are created from the three FOM computations with a sampling frequency $n_{\freq}=10$, and we assess the performance of ROM by computing the flow past cylinders with the same radius and locations $x_c=-0.65$, $x_c=-0.15$, $x_c=0.3$, and $x_c=0.8$, with termination times given by $T=0.4$, $T=0.7$, $T=0.6$, and $T=0.3$, respectively.

\begin{figure}[h]\centering
  \includegraphics[width=.48\textwidth]{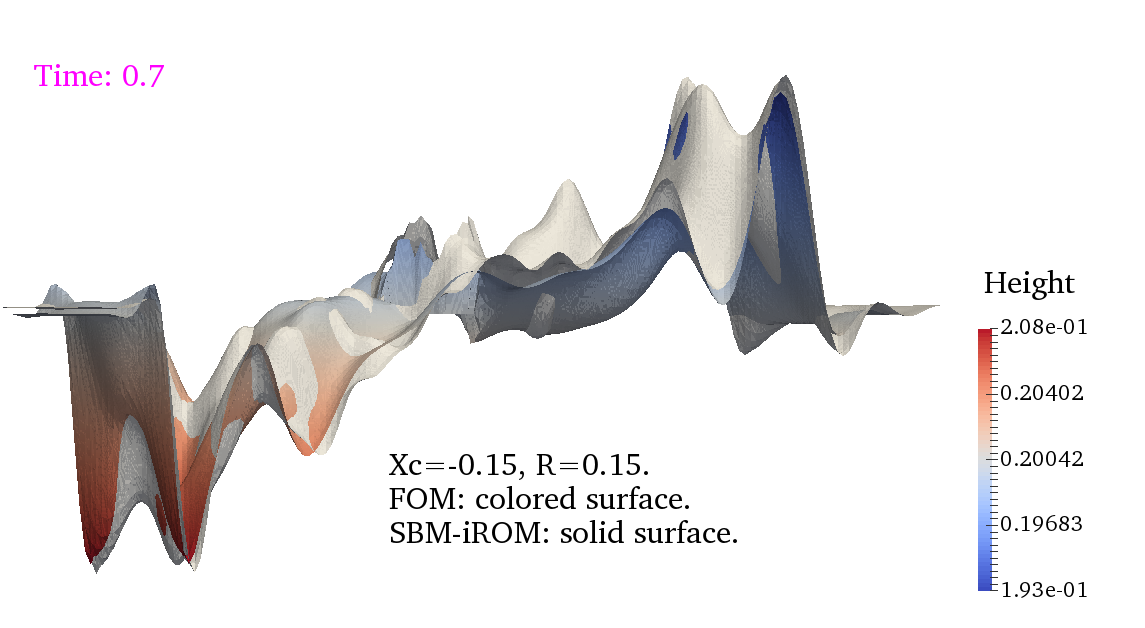}
  \includegraphics[width=.48\textwidth]{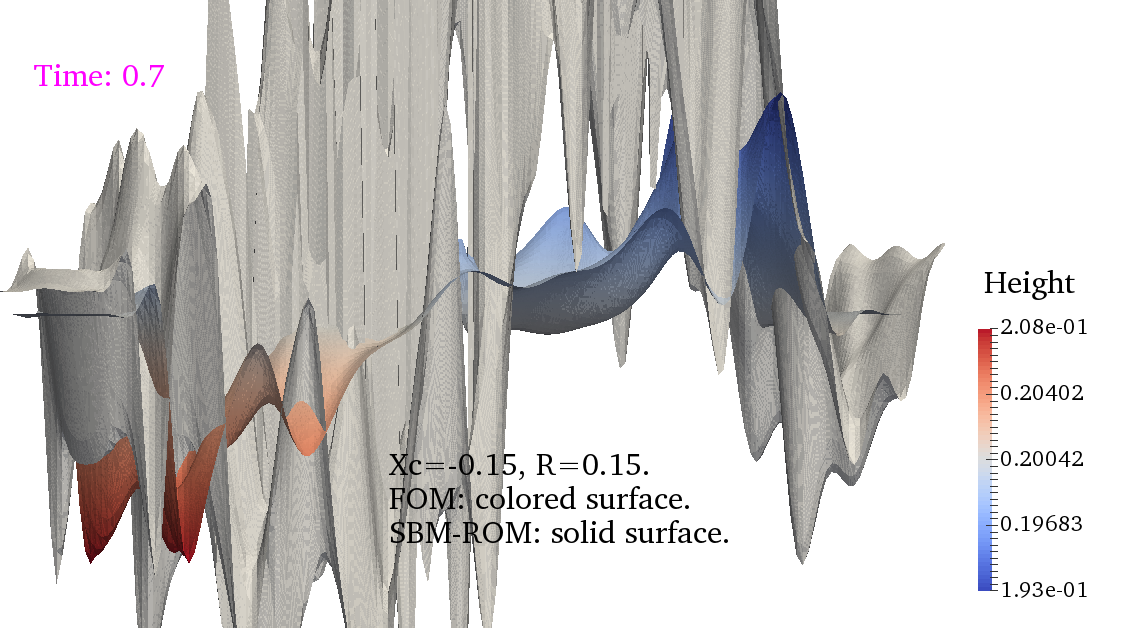}
  \caption{The solution surfaces for water height at $T=0.7$ in the case $x_c=-0.15$ of Test 2 for SBM-iROM (left panel) and SBM-ROM (right panel). In both plots, the ROM solution (solid surface) is plotted on top of the FOM one (color surface).}
  \label{fg:num_cyl_x_warp_x-0d15}
\end{figure}

With varying $x_c$, SBM-ROM is highly unstable, as demonstrated by the height solution surfaces for the case $x_c=-0.15$ in Figure~\ref{fg:num_cyl_x_warp_x-0d15}; in contrast, the SBM-iROM solution agrees very well with the FOM one.

\begin{figure}[h]\centering
  \includegraphics[trim=0.0in 3.0in 0.0in 0.0in, clip, width=.48\textwidth]{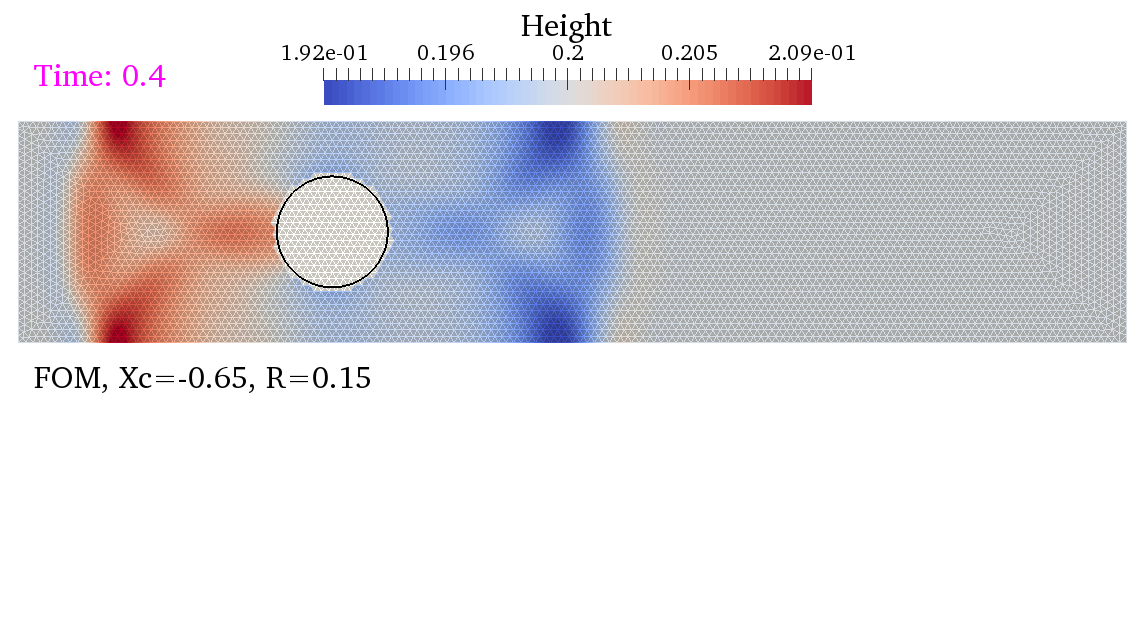} 
  \includegraphics[trim=0.0in 3.0in 0.0in 0.0in, clip, width=.48\textwidth]{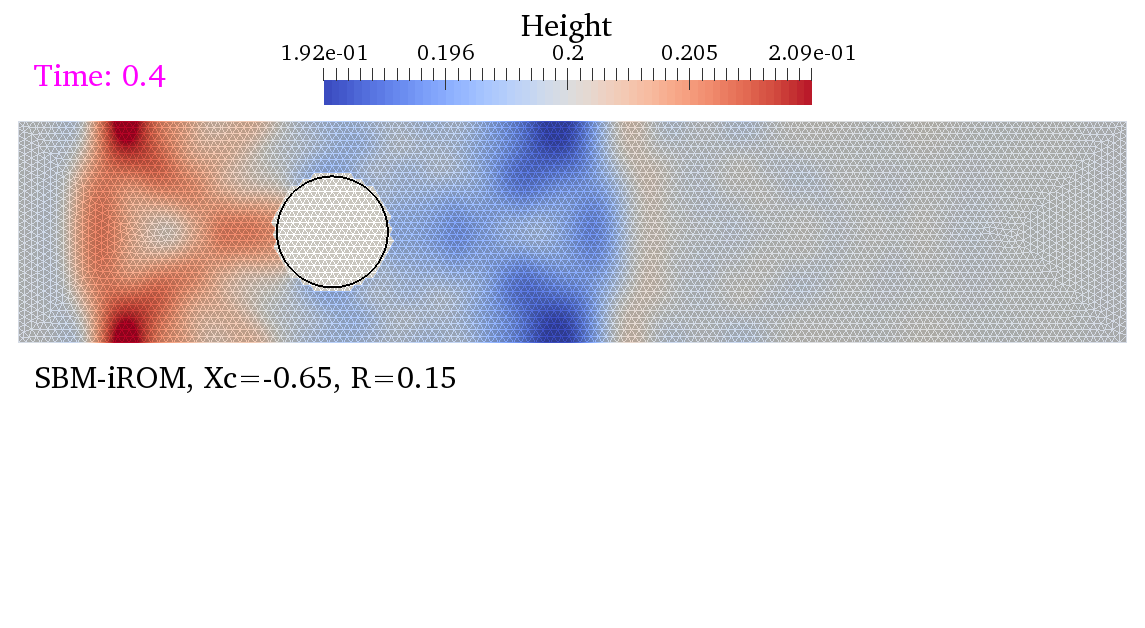} 
  \caption{The height ($h$) in the case of $x_c=-0.65$ computed by FOM (left panel) and SBM-iROM (right panel), with the legend range set according to the FOM computation.
  The ROM computation is performed using POD modes corresponding to $\mu_{\textrm{pod}}=1-10^{-6}$.}
  \label{fg:num_cyl_x-0d65_sol}
\end{figure}

\begin{figure}[h]\centering
  \includegraphics[trim=0.0in 3.0in 0.0in 0.0in, clip, width=.48\textwidth]{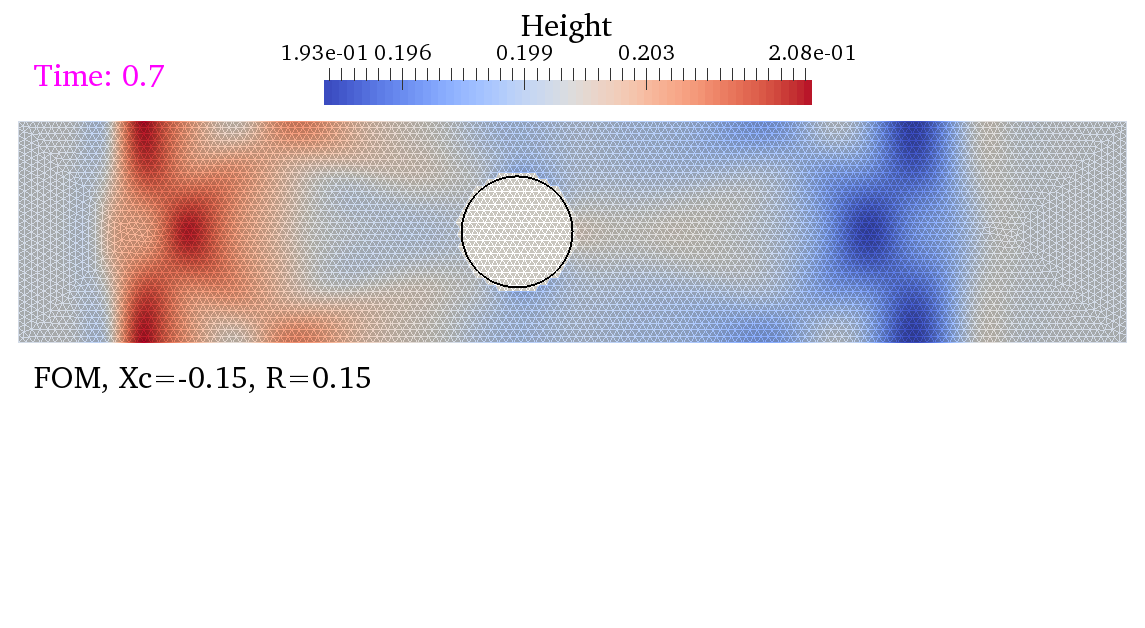} 
  \includegraphics[trim=0.0in 3.0in 0.0in 0.0in, clip, width=.48\textwidth]{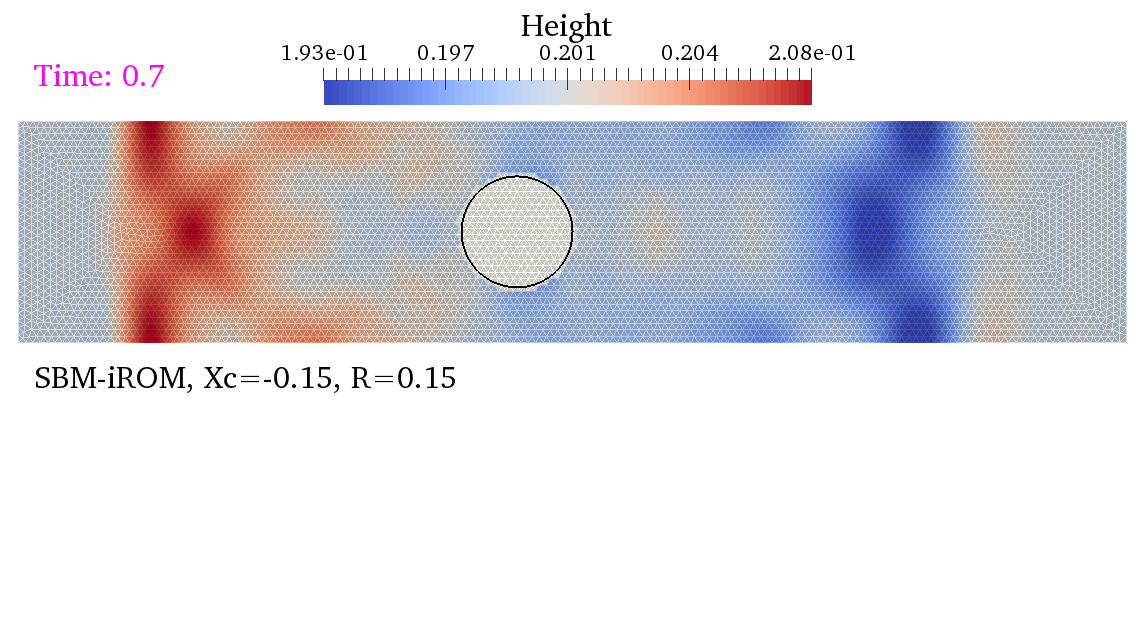} 
  \caption{The height ($h$) in the case of $x_c=-0.15$ computed by FOM (left panel) and SBM-iROM (right panel), with the legend range set according to the FOM computation.
  The ROM computation is performed using POD modes corresponding to $\mu_{\textrm{pod}}=1-10^{-6}$.}
  \label{fg:num_cyl_x-0d15_sol}
\end{figure}

\begin{figure}[h]\centering
  \includegraphics[trim=0.0in 3.0in 0.0in 0.0in, clip, width=.48\textwidth]{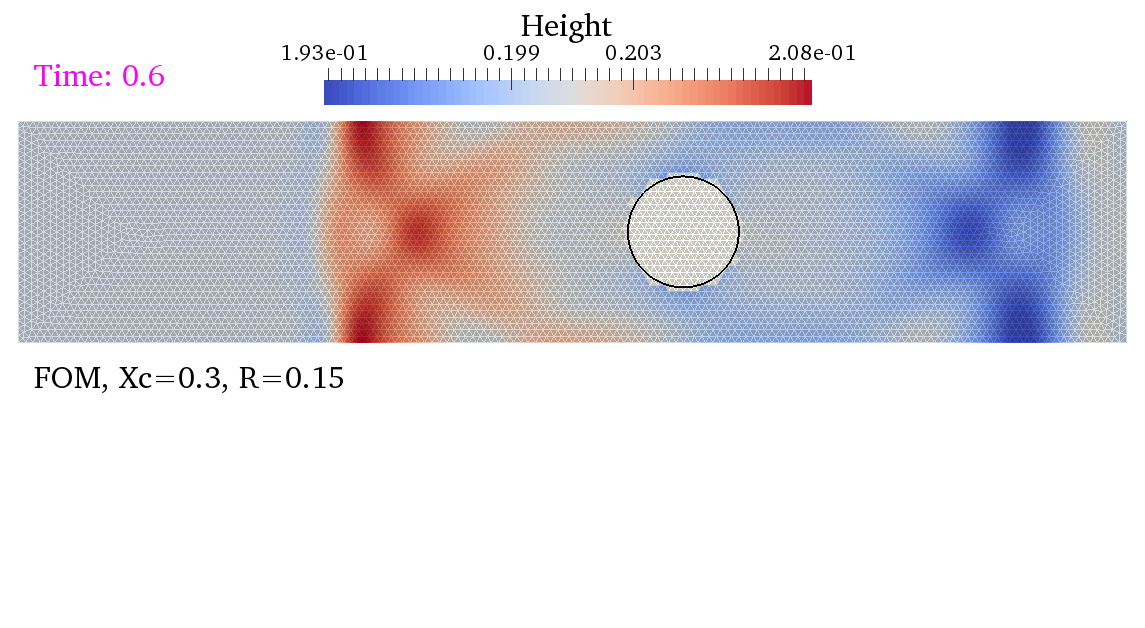} 
  \includegraphics[trim=0.0in 3.0in 0.0in 0.0in, clip, width=.48\textwidth]{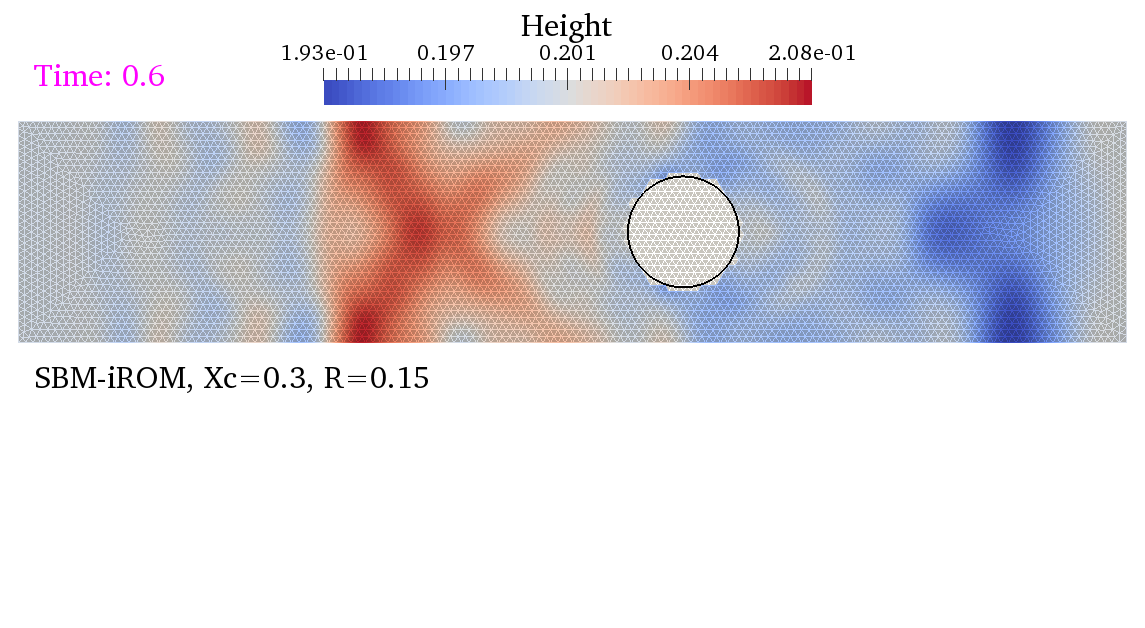} 
  \caption{The height ($h$) in the case of $x_c=0.3$ computed by FOM (left panel) and SBM-iROM (right panel), with the legend range set according to the FOM computation.
  The ROM computation is performed using POD modes corresponding to $\mu_{\textrm{pod}}=1-10^{-6}$.}
  \label{fg:num_cyl_x0d3_sol}
\end{figure}

\begin{figure}[h]\centering
  \includegraphics[trim=0.0in 3.0in 0.0in 0.0in, clip, width=.48\textwidth]{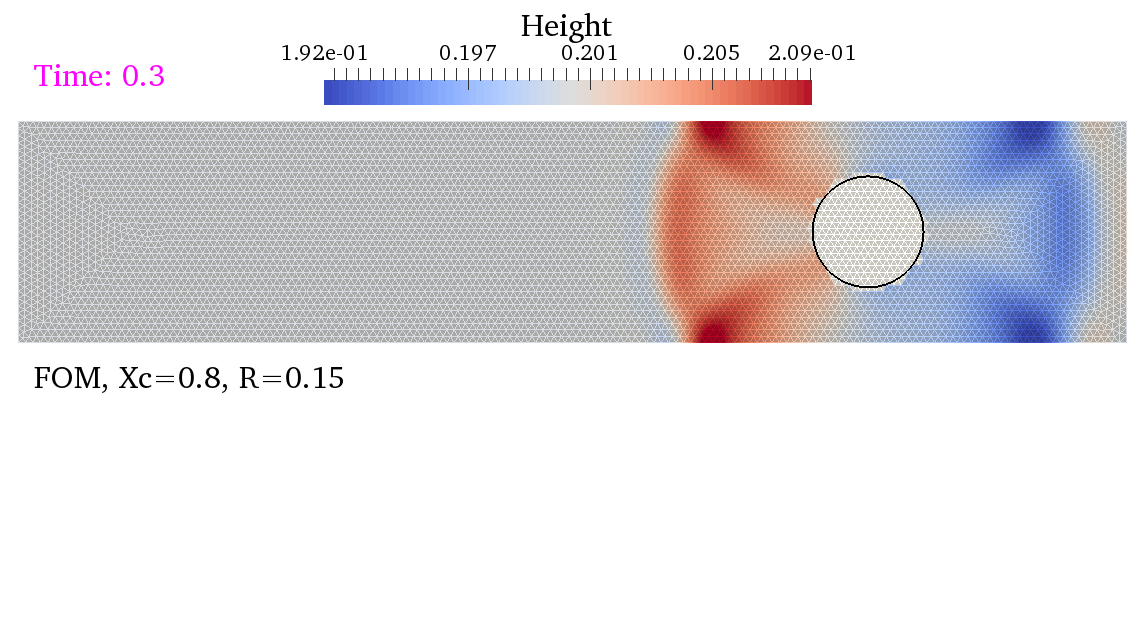} 
  \includegraphics[trim=0.0in 3.0in 0.0in 0.0in, clip, width=.48\textwidth]{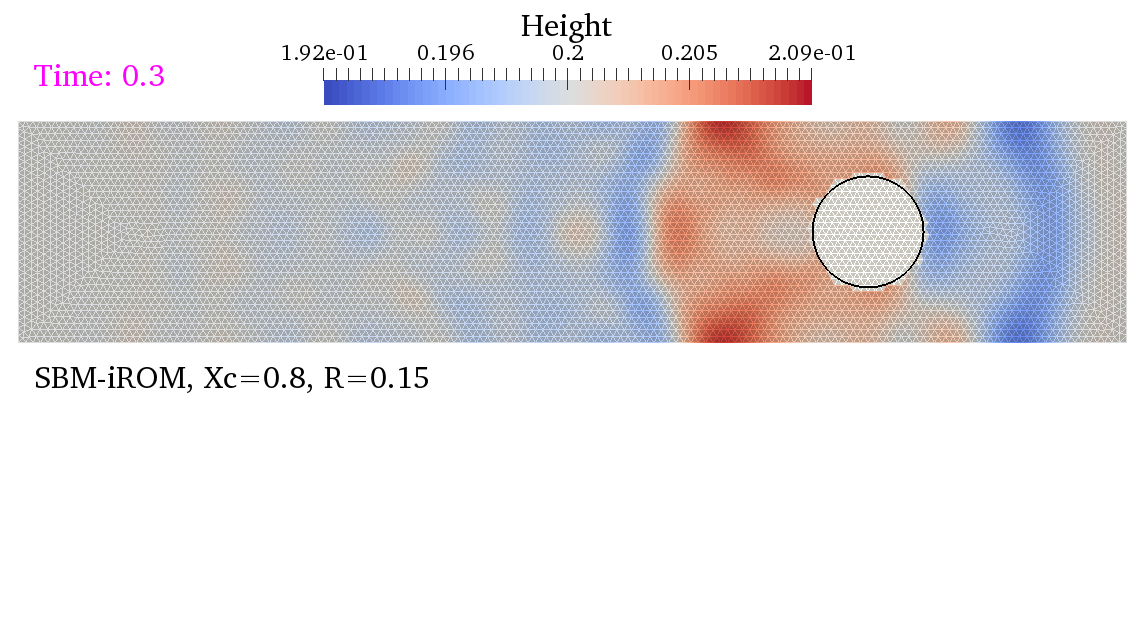} 
  \caption{The height ($h$) in the case of $x_c=0.8$ computed by FOM (left panel) and SBM-iROM (right panel), with the legend range set according to the FOM computation.
  The ROM computation is performed using POD modes corresponding to $\mu_{\textrm{pod}}=1-10^{-6}$.}
  \label{fg:num_cyl_x0d8_sol}
\end{figure}

\begin{figure}[h]\centering
  \begin{subfigure}[t]{.4\textwidth}\centering
    \includegraphics[width=\textwidth]{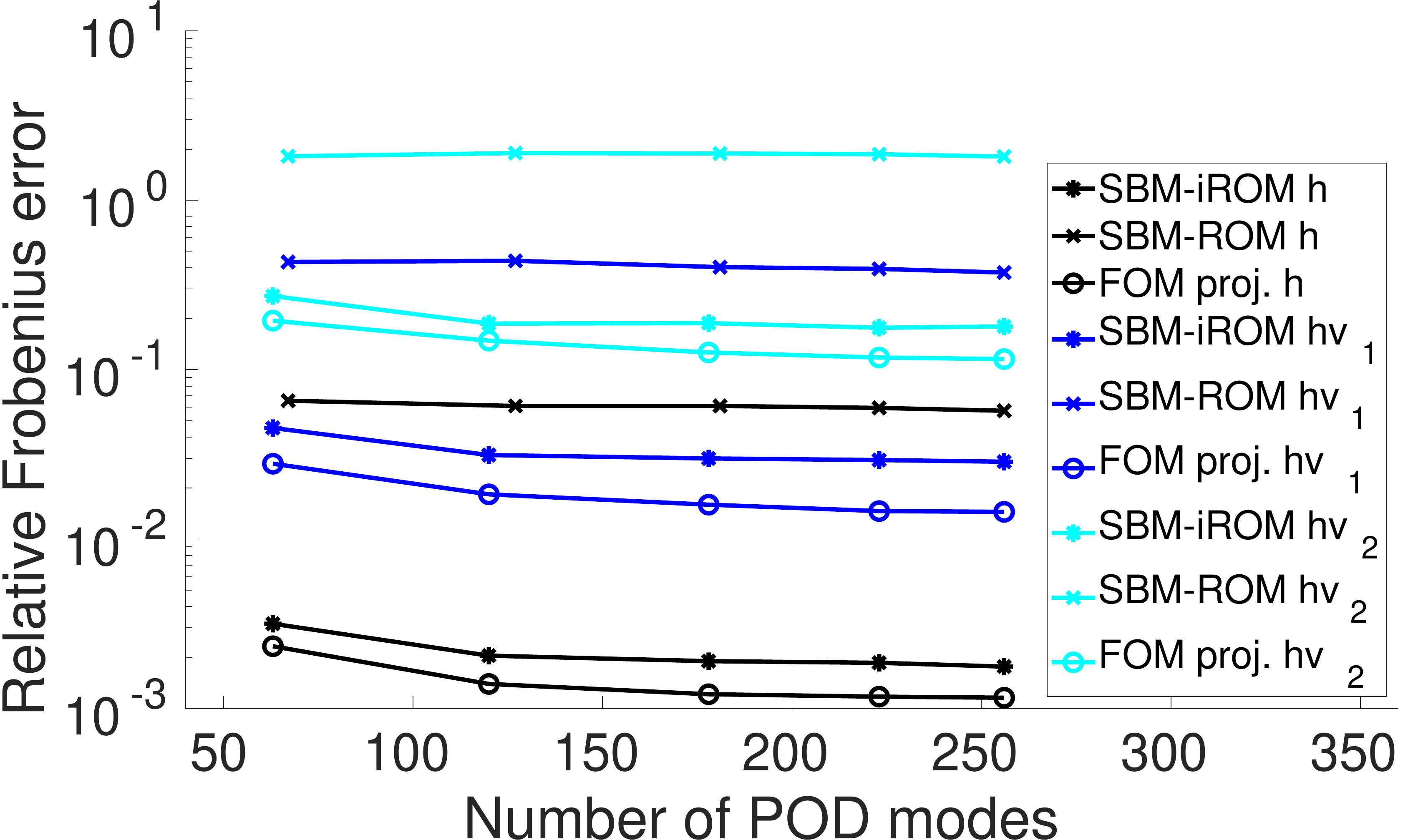}
    \caption{$x_c=-0.65$.}
    \label{fg:num_cyl_x_froberr_x-0d65}
  \end{subfigure}
  \begin{subfigure}[t]{.4\textwidth}\centering
    \includegraphics[width=\textwidth]{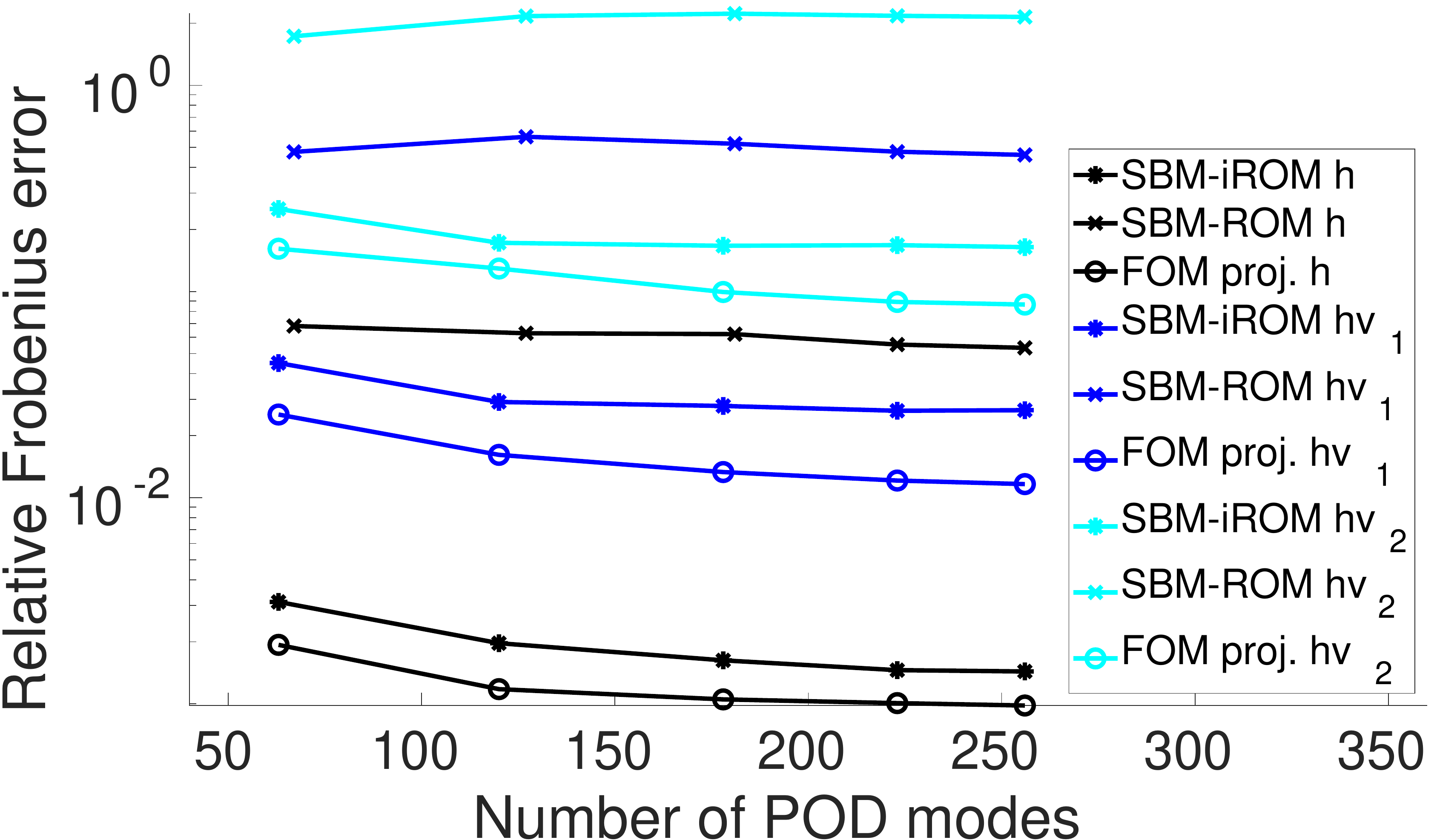}
    \caption{$x_c=-0.15$.}
    \label{fg:num_cyl_x_froberr_x-0d15}
  \end{subfigure} \\
  \begin{subfigure}[t]{.4\textwidth}\centering
    \includegraphics[width=\textwidth]{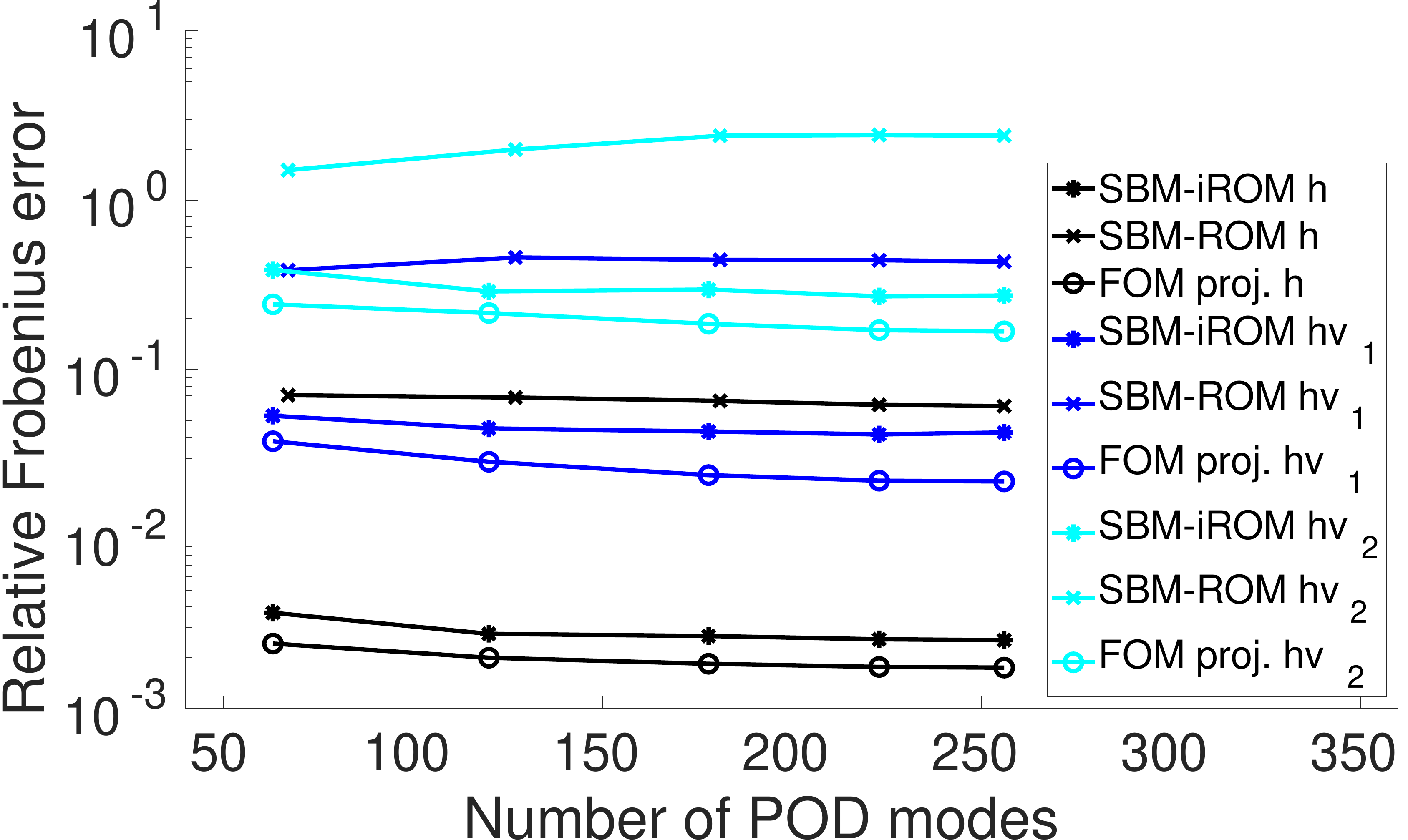}
    \caption{$x_c=0.3$.}
    \label{fg:num_cyl_x_froberr_x0d3}
  \end{subfigure}
  \begin{subfigure}[t]{.4\textwidth}\centering
    \includegraphics[width=\textwidth]{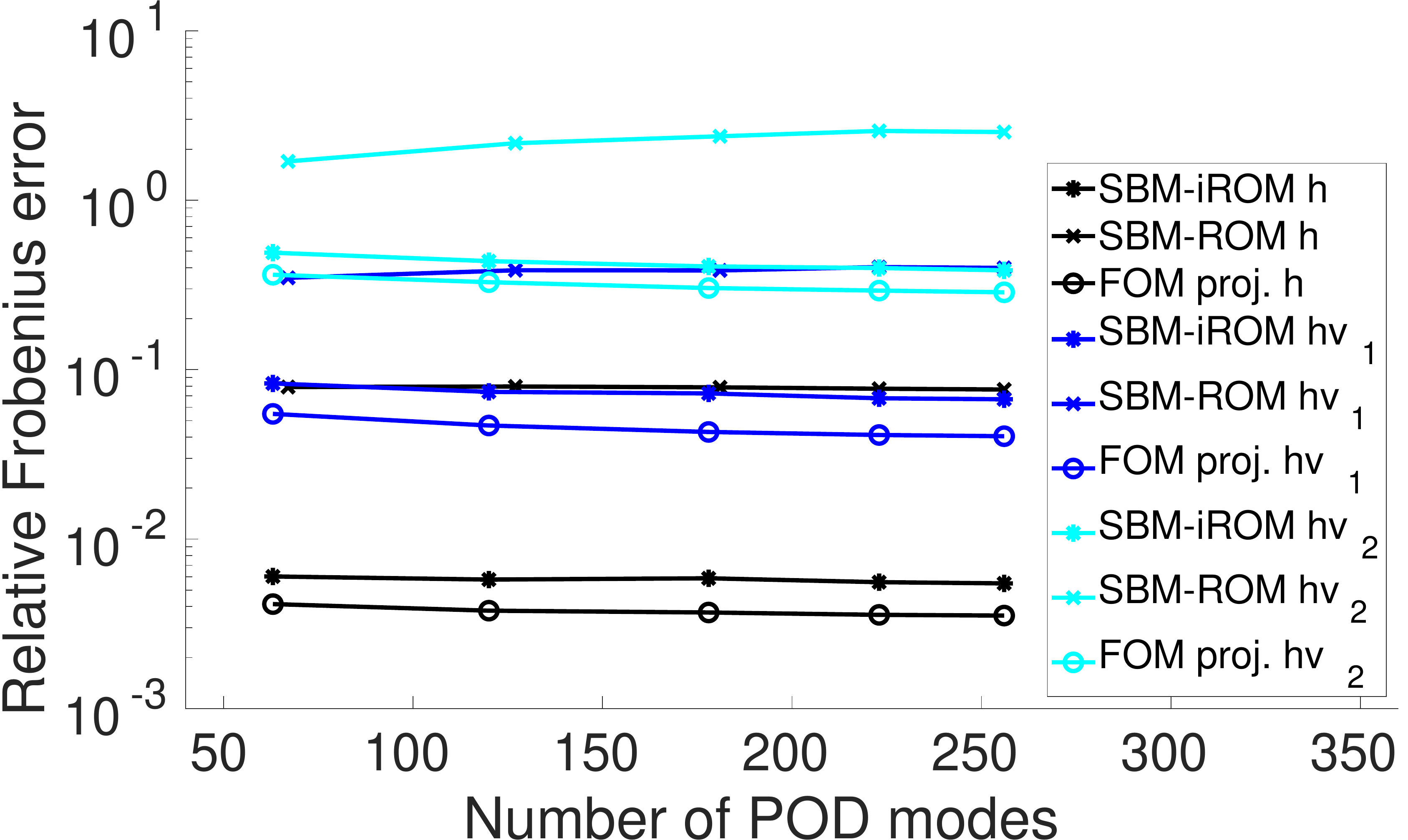}
    \caption{$x_c=0.8$.}
    \label{fg:num_cyl_x_froberr_x0d8}
  \end{subfigure}
  \caption{The relative Frobenius errors of Test 2 by SBM-iROM (marked by $\ast$), SBM-ROM (marked by x), and FOM projection error (FOM proj., marked by o) in logarithmic scale plotted against the number of POD modes.
  The errors in $h$, $hv_1$, $hv_2$ are plotted in black, blue, and cyan colors, respectively.}
  \label{fg:num_cyl_x_froberr}
\end{figure}

Next in Figures~\ref{fg:num_cyl_x-0d65_sol}--\ref{fg:num_cyl_x0d8_sol}, we plot the water height solutions at terminal times using FOM and SBM-iROM; the SBM-ROM solutions are omitted as they are highly unstable, as also seen by the relative space-time Frobenius errors reported in Tables~\ref{tb:num_cyl_x_froberr_xc-0d65}--\ref{tb:num_cyl_x_froberr_xc0d8} and Figure~\ref{fg:num_cyl_x_froberr}.
\begin{table}[h]\centering
  \caption{The relative space-time Frobenius errors of the projected FOM solutions and ROM computations in the case $x_c=-0.65$ of Test 2.}
  \label{tb:num_cyl_x_froberr_xc-0d65}
  \begin{tabular}{@{}lccccccccccc@{}}
  \toprule[.5mm]
  & \multicolumn{3}{l}{FOM projection} & & \multicolumn{3}{l}{SBM-iROM} & & \multicolumn{3}{l}{SBM-ROM} \\ \cmidrule[.2mm](lr){2-4} \cmidrule[.2mm](lr){6-8} \cmidrule[.2mm](l){10-12}
  $\mu_{\textrm{pod}}$ & $h$ & $hv_1$ & $hv_2$ & & $h$ & $hv_1$ & $hv_2$ & & $h$ & $hv_1$ & $hv_2$ \\ \cmidrule[.3mm](l){2-12}
  $1-10^{-5}$ & 2.33e-3 & 2.78e-2 & 1.95e-1 & & 3.17e-3 & 4.53e-2 & 2.72e-1 & & 6.56e-2 & 4.32e-1 & 1.82e+0 \\
  $1-10^{-6}$ & 1.40e-3 & 1.84e-2 & 1.48e-1 & & 2.05e-3 & 3.13e-2 & 1.87e-1 & & 6.10e-2 & 4.39e-1 & 1.90e+0 \\
  $1-10^{-7}$ & 1.21e-3 & 1.60e-2 & 1.26e-1 & & 1.90e-3 & 2.99e-2 & 1.88e-1 & & 6.11e-2 & 4.03e-1 & 1.89e+0 \\
  $1-10^{-8}$ & 1.17e-3 & 1.46e-2 & 1.18e-1 & & 1.86e-3 & 2.93e-2 & 1.77e-1 & & 5.94e-2 & 3.93e-1 & 1.87e+0 \\
  $1-10^{-9}$ & 1.16e-3 & 1.45e-2 & 1.15e-1 & & 1.77e-3 & 2.86e-2 & 1.80e-1 & & 5.72e-2 & 3.74e-1 & 1.81e+0 \\
  \bottomrule[.4mm]
  \end{tabular}
\end{table}
\begin{table}[h]\centering
  \caption{The relative space-time Frobenius errors of the projected FOM solutions and ROM computations in the case $x_c=-0.15$ of Test 2.}
  \label{tb:num_cyl_x_froberr_xc-0d15}
  \begin{tabular}{@{}lccccccccccc@{}}
  \toprule[.5mm]
  & \multicolumn{3}{l}{FOM projection} & & \multicolumn{3}{l}{SBM-iROM} & & \multicolumn{3}{l}{SBM-ROM} \\ \cmidrule[.2mm](lr){2-4} \cmidrule[.2mm](lr){6-8} \cmidrule[.2mm](l){10-12}
  $\mu_{\textrm{pod}}$ & $h$ & $hv_1$ & $hv_2$ & & $h$ & $hv_1$ & $hv_2$ & & $h$ & $hv_1$ & $hv_2$ \\ \cmidrule[.3mm](l){2-12}
  $1-10^{-5}$ & 1.93e-3 & 2.53e-2 & 1.61e-1 & & 3.12e-3 & 4.50e-2 & 2.51e-1 & & 6.80e-2 & 4.76e-1 & 1.73e+0 \\
  $1-10^{-6}$ & 1.18e-3 & 1.61e-2 & 1.29e-1 & & 1.97e-3 & 2.91e-2 & 1.72e-1 & & 6.27e-2 & 5.63e-1 & 2.17e+0 \\
  $1-10^{-7}$ & 1.05e-3 & 1.33e-2 & 9.94e-2 & & 1.63e-3 & 2.78e-2 & 1.67e-1 & & 6.21e-2 & 5.21e-1 & 2.23e+0 \\
  $1-10^{-8}$ & 1.01e-3 & 1.21e-2 & 8.90e-2 & & 1.46e-3 & 2.64e-2 & 1.68e-1 & & 5.53e-2 & 4.76e-1 & 2.17e+0 \\
  $1-10^{-9}$ & 9.82e-4 & 1.16e-2 & 8.64e-2 & & 1.44e-3 & 2.66e-2 & 1.64e-1 & & 5.33e-2 & 4.60e-1 & 2.15e+0 \\
  \bottomrule[.4mm]
  \end{tabular}
\end{table}
\begin{table}[h]\centering
  \caption{The relative space-time Frobenius errors of the projected FOM solutions and ROM computations in the case $x_c=0.3$ of Test 2.}
  \label{tb:num_cyl_x_froberr_xc0d3}
  \begin{tabular}{@{}lccccccccccc@{}}
  \toprule[.5mm]
  & \multicolumn{3}{l}{FOM projection} & & \multicolumn{3}{l}{SBM-iROM} & & \multicolumn{3}{l}{SBM-ROM} \\ \cmidrule[.2mm](lr){2-4} \cmidrule[.2mm](lr){6-8} \cmidrule[.2mm](l){10-12}
  $\mu_{\textrm{pod}}$ & $h$ & $hv_1$ & $hv_2$ & & $h$ & $hv_1$ & $hv_2$ & & $h$ & $hv_1$ & $hv_2$ \\ \cmidrule[.3mm](l){2-12}
  $1-10^{-5}$ & 2.42e-3 & 3.78e-2 & 2.43e-1 & & 3.67e-3 & 5.34e-2 & 3.88e-1 & & 7.06e-2 & 3.86e-1 & 1.50e+0 \\
  $1-10^{-6}$ & 1.99e-3 & 2.86e-2 & 2.16e-1 & & 2.76e-3 & 4.50e-2 & 2.90e-1 & & 6.85e-2 & 4.59e-1 & 1.99e+0 \\
  $1-10^{-7}$ & 1.84e-3 & 2.38e-2 & 1.86e-1 & & 2.68e-3 & 4.31e-2 & 2.97e-1 & & 6.54e-2 & 4.45e-1 & 2.40e+0 \\
  $1-10^{-8}$ & 1.76e-3 & 2.21e-2 & 1.71e-1 & & 2.56e-3 & 4.15e-2 & 2.71e-1 & & 6.20e-2 & 4.43e-1 & 2.42e+0 \\
  $1-10^{-9}$ & 1.74e-3 & 2.19e-2 & 1.68e-1 & & 2.53e-3 & 4.26e-2 & 2.74e-1 & & 6.09e-2 & 4.34e-1 & 2.40e+0 \\
  \bottomrule[.4mm]
  \end{tabular}
\end{table}
\begin{table}[h]\centering
  \caption{The relative space-time Frobenius errors of the projected FOM solutions and ROM computations in the case $x_c=0.8$ of Test 2.}
  \label{tb:num_cyl_x_froberr_xc0d8}
  \begin{tabular}{@{}lccccccccccc@{}}
  \toprule[.5mm]
  & \multicolumn{3}{l}{FOM projection} & & \multicolumn{3}{l}{SBM-iROM} & & \multicolumn{3}{l}{SBM-ROM} \\ \cmidrule[.2mm](lr){2-4} \cmidrule[.2mm](lr){6-8} \cmidrule[.2mm](l){10-12}
  $\mu_{\textrm{pod}}$ & $h$ & $hv_1$ & $hv_2$ & & $h$ & $hv_1$ & $hv_2$ & & $h$ & $hv_1$ & $hv_2$ \\ \cmidrule[.3mm](l){2-12}
  $1-10^{-5}$ & 4.13e-3 & 5.48e-2 & 3.63e-1 & & 6.02e-3 & 8.29e-2 & 4.90e-1 & & 7.90e-2 & 3.48e-1 & 1.69e+0 \\
  $1-10^{-6}$ & 3.78e-3 & 4.68e-2 & 3.28e-1 & & 5.77e-3 & 7.39e-2 & 4.37e-1 & & 7.96e-2 & 3.86e-1 & 2.17e+0 \\
  $1-10^{-7}$ & 3.69e-3 & 4.29e-2 & 3.03e-1 & & 5.87e-3 & 7.24e-2 & 4.07e-1 & & 7.86e-2 & 3.85e-1 & 2.38e+0 \\
  $1-10^{-8}$ & 3.57e-3 & 4.11e-2 & 2.93e-1 & & 5.57e-3 & 6.78e-2 & 3.98e-1 & & 7.72e-2 & 4.03e-1 & 2.56e+0 \\
  $1-10^{-9}$ & 3.54e-3 & 4.05e-2 & 2.86e-1 & & 5.48e-3 & 6.69e-2 & 3.86e-1 & & 7.64e-2 & 3.99e-1 & 2.52e+0 \\
  \bottomrule[.4mm]
  \end{tabular}
\end{table}

\subsection{Test 3: Geometrical parameterization study with two-dimensional parameter space}
\label{sec:num_cyl_xr}
\begin{figure}[h]\centering
  \includegraphics[width=.48\textwidth]{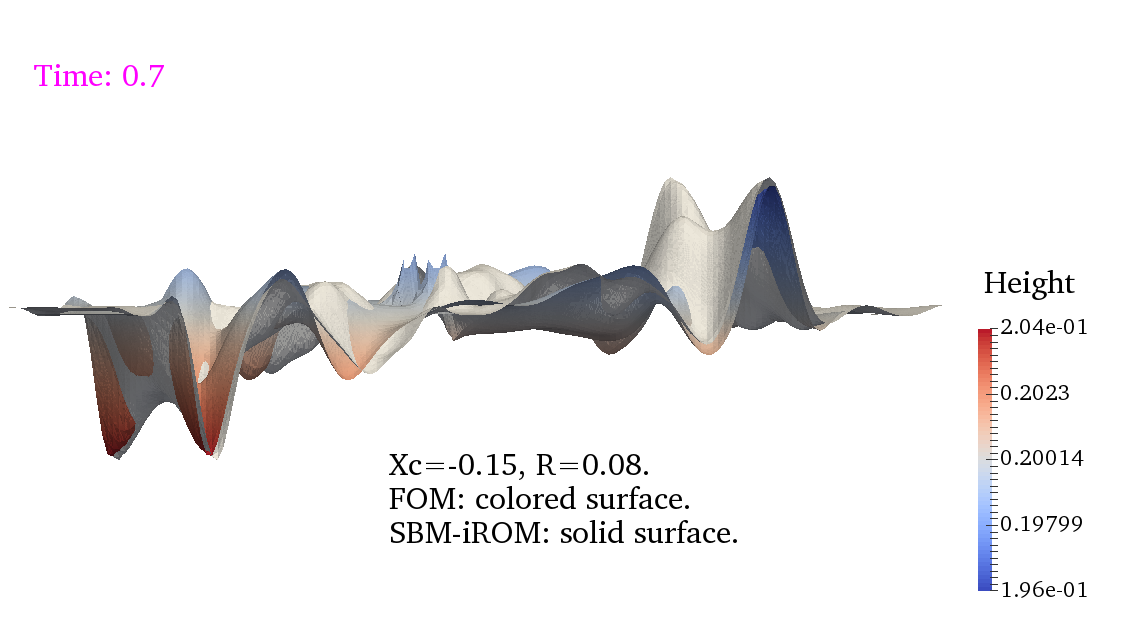}
  \includegraphics[width=.48\textwidth]{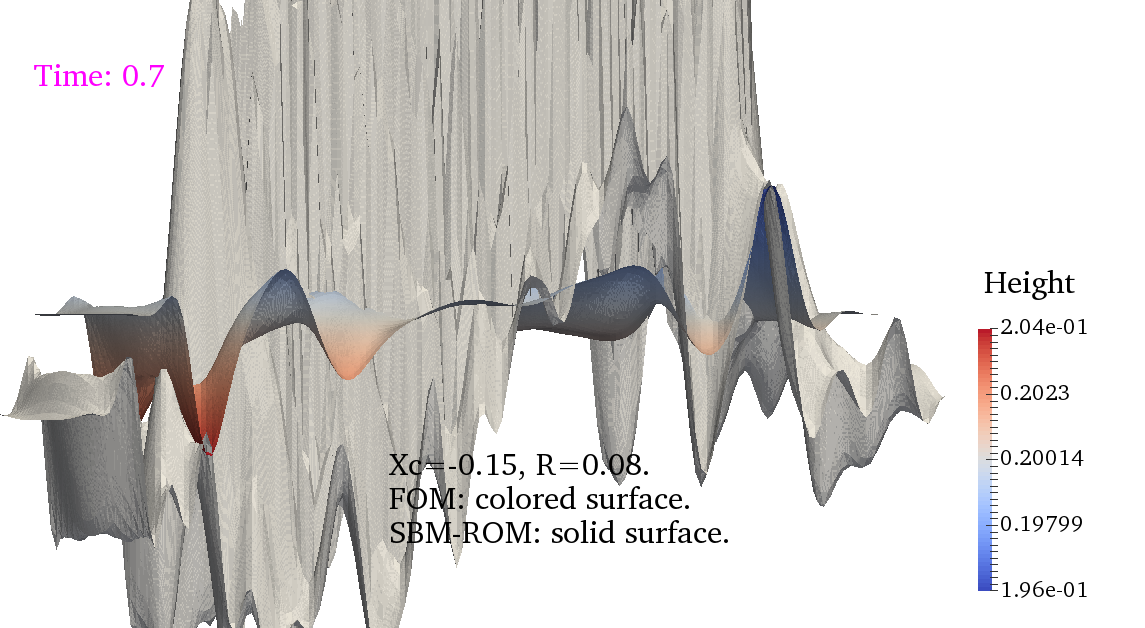}
  \caption{The solution surfaces for water height at $T=0.7$ in the case $R=0.08$ and $x_c=-0.15$ of Test 3 for SBM-iROM (left panel) and SBM-ROM (right panel). In both plots, the ROM solution (solid surface) is plotted on top of the FOM one (color surface).}
  \label{fg:num_cyl_xr_warp_r0d08_x-0d15}
\end{figure}

\begin{figure}[h]\centering
  \includegraphics[width=.48\textwidth]{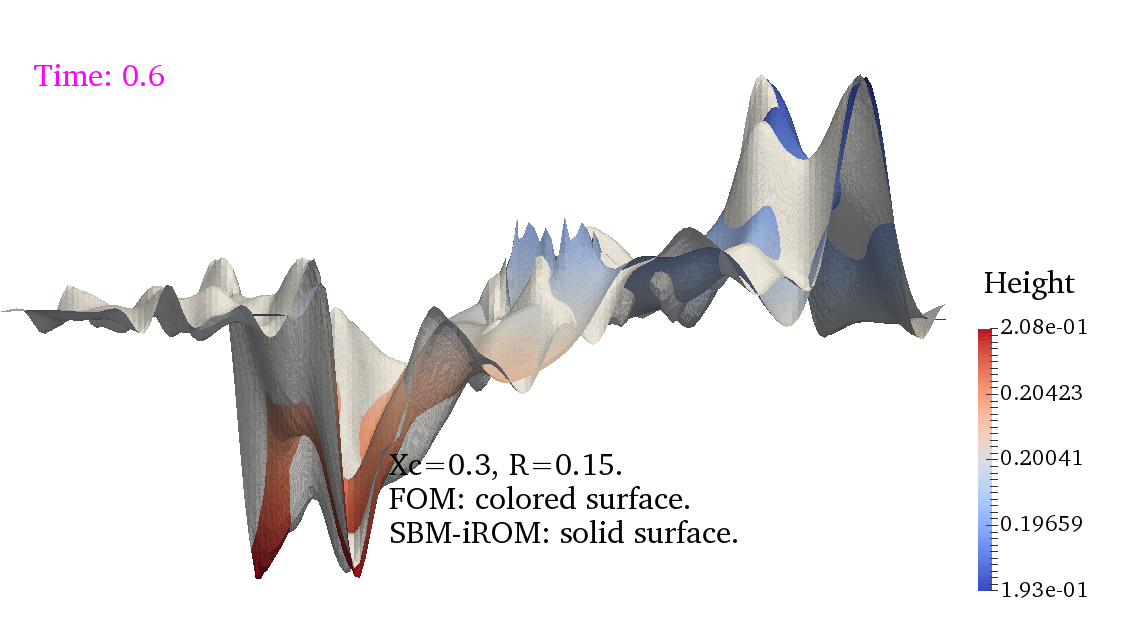}
  \includegraphics[width=.48\textwidth]{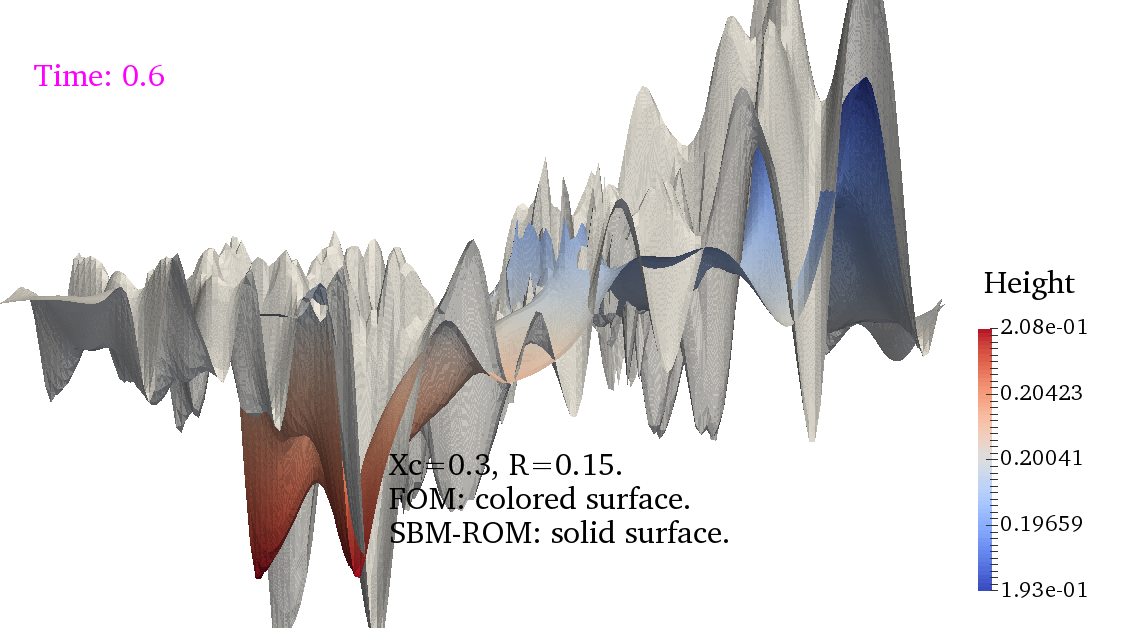}
  \caption{The solution surfaces for water height at $T=0.6$ in the case $R=0.15$ and $x_c=0.3$ of Test 3 for SBM-iROM (left panel) and SBM-ROM (right panel). In both plots, the ROM solution (solid surface) is plotted on top of the FOM one (color surface).}
  \label{fg:num_cyl_xr_warp_r0d15_x0d3}
\end{figure}

\begin{figure}[h]\centering
  \begin{subfigure}[t]{\textwidth}\centering
    \includegraphics[trim=0.0in 3.0in 0.0in 0.0in, clip, width=.48\textwidth]{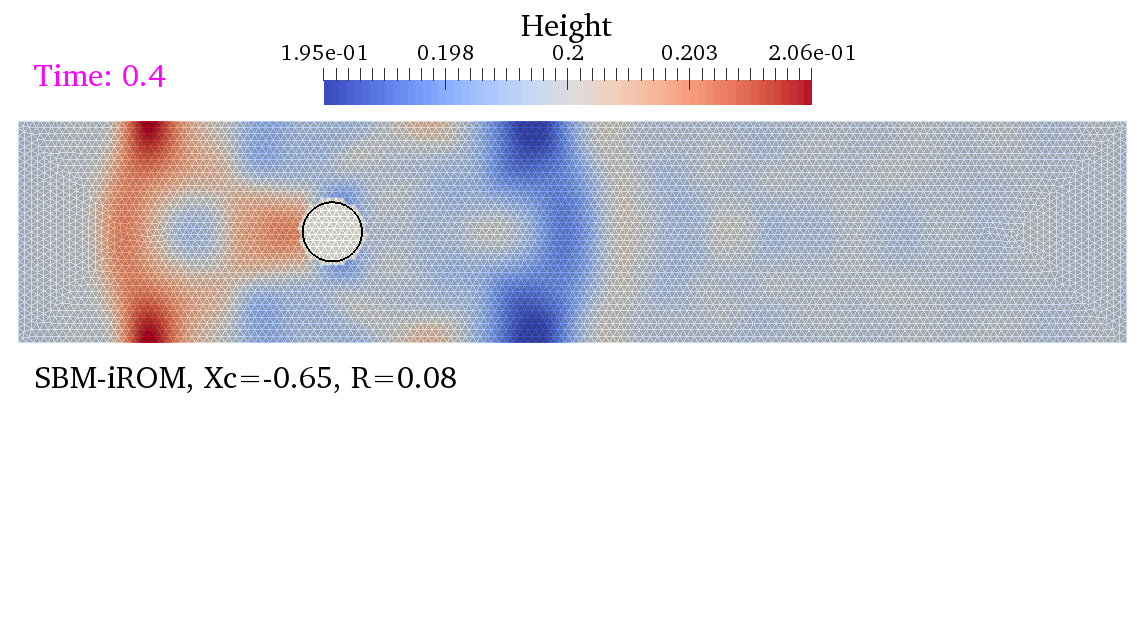}
    \includegraphics[trim=0.0in 3.0in 0.0in 0.0in, clip, width=.48\textwidth]{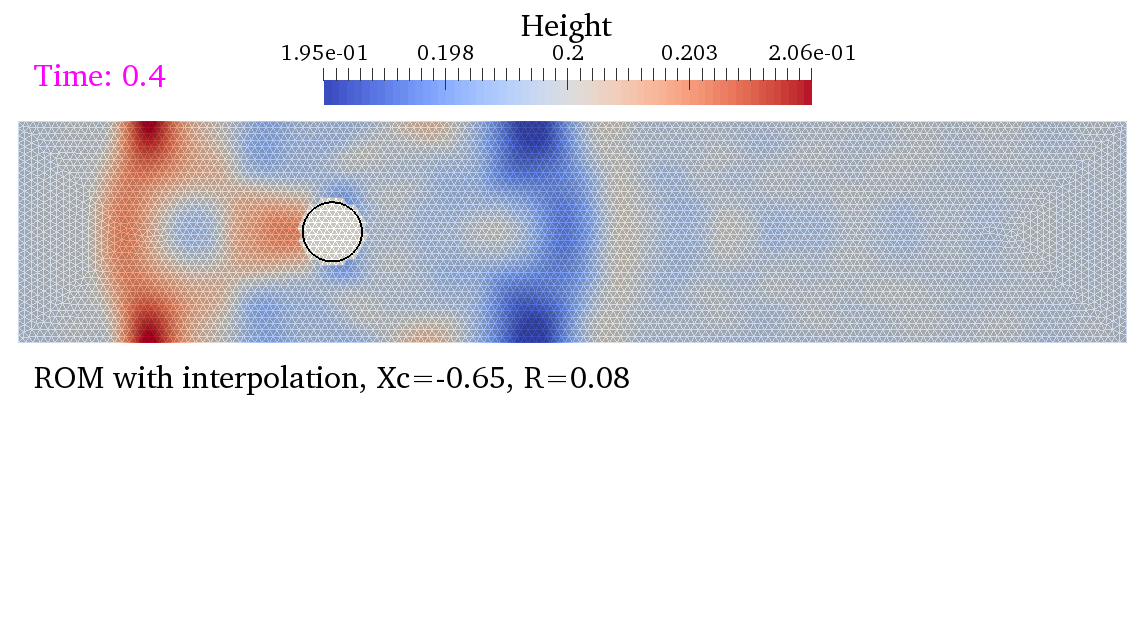}
    \caption{$R=0.08$.}
  \end{subfigure}
  \begin{subfigure}[t]{\textwidth}\centering
    \includegraphics[trim=0.0in 3.0in 0.0in 0.0in, clip, width=.48\textwidth]{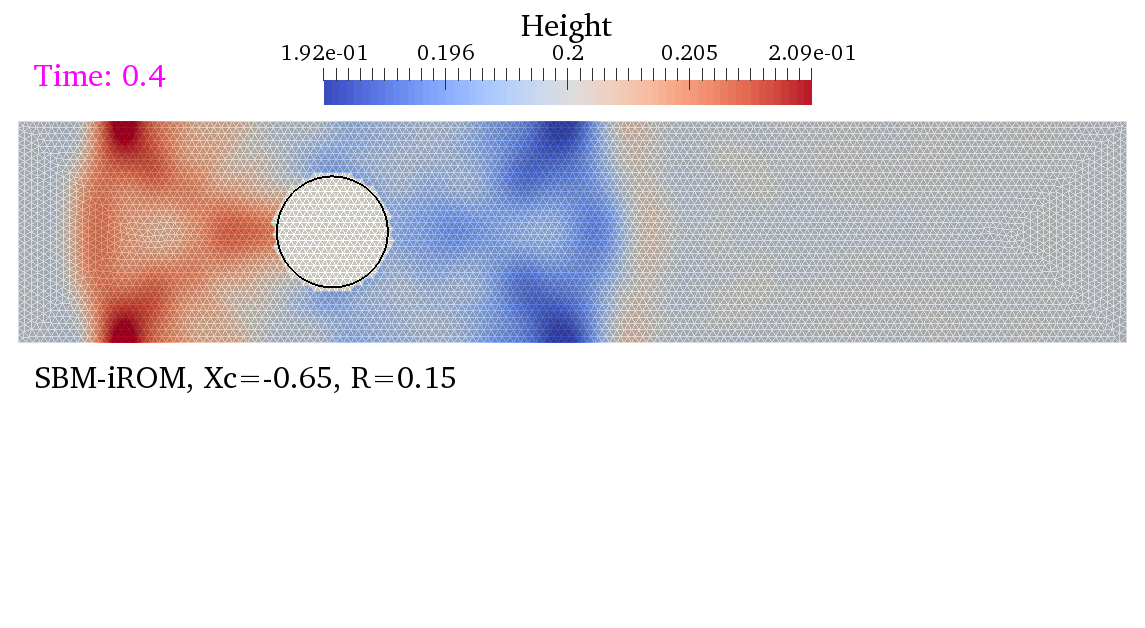}
    \includegraphics[trim=0.0in 3.0in 0.0in 0.0in, clip, width=.48\textwidth]{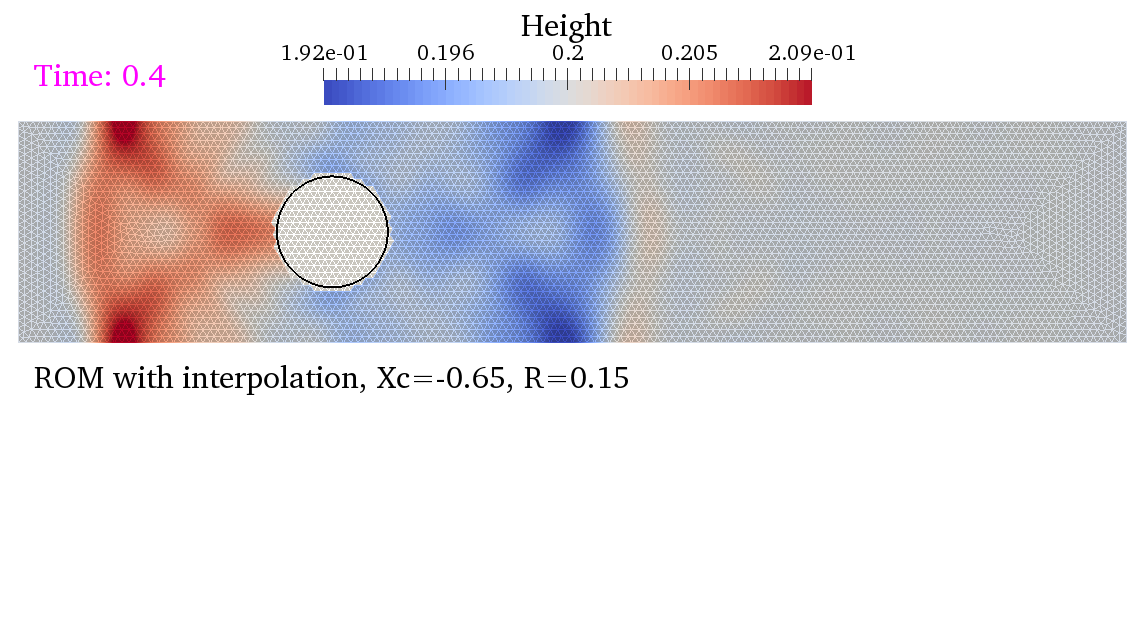}
    \caption{$R=0.15$.}
  \end{subfigure}
  \begin{subfigure}[t]{\textwidth}\centering
    \includegraphics[trim=0.0in 3.0in 0.0in 0.0in, clip, width=.48\textwidth]{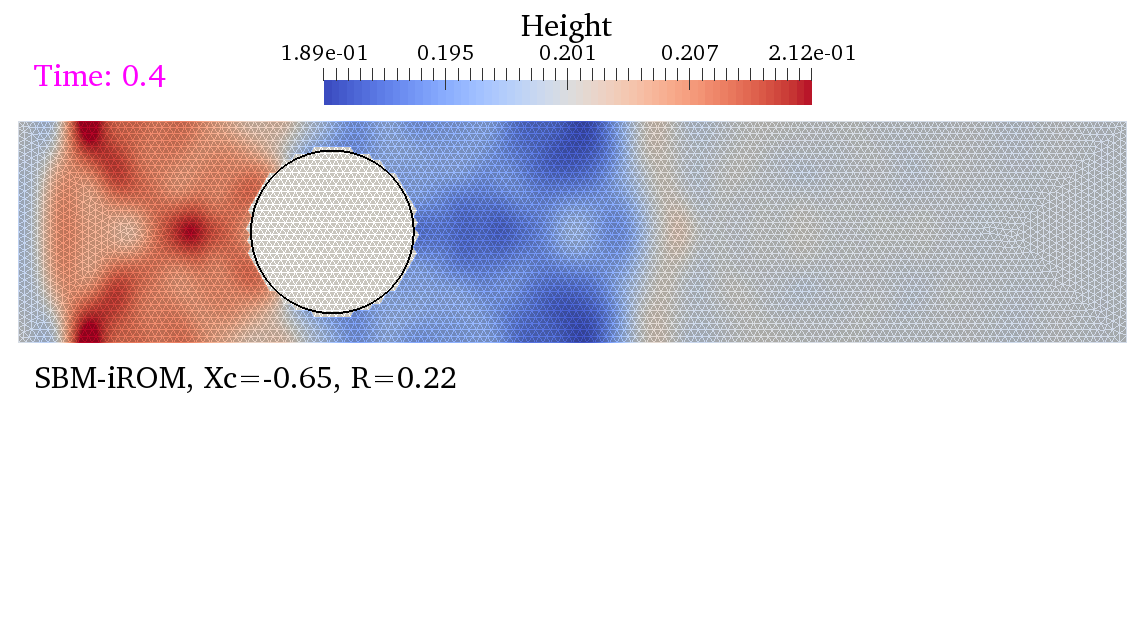}
    \includegraphics[trim=0.0in 3.0in 0.0in 0.0in, clip, width=.48\textwidth]{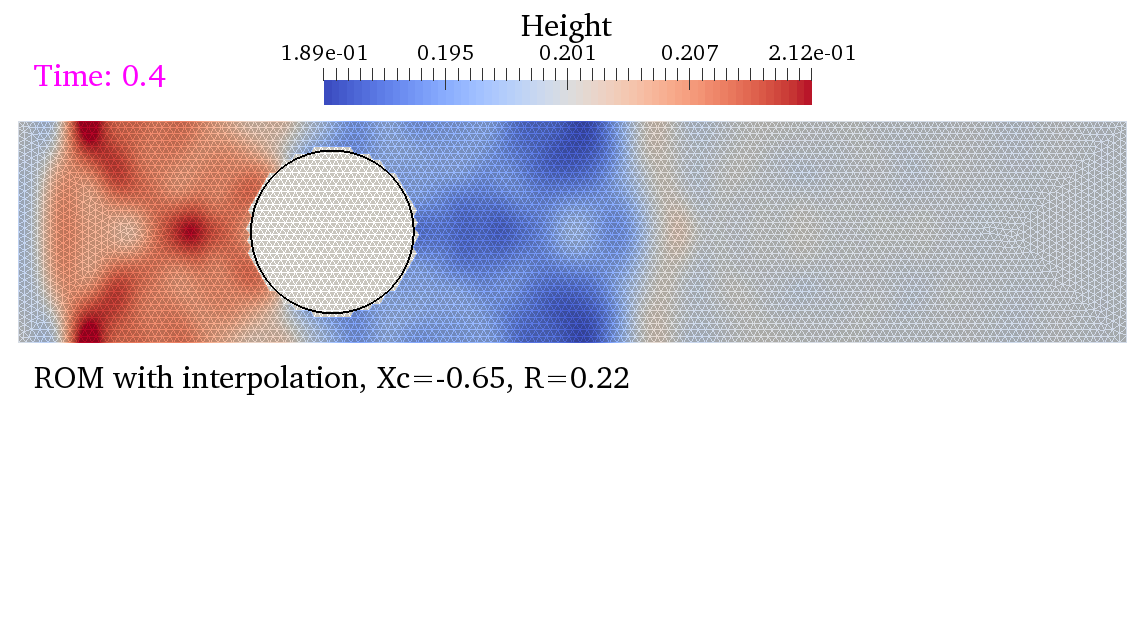}
    \caption{$R=0.22$.}
  \end{subfigure}
  \caption{The height ($h$) in the case of $R=0.08$ (top row), $R=0.15$ (middle row), or $R=0.22$ (bottom row), and $x_c=-0.65$ computed by FOM (left panels) and SBM-iROM (right panels), with the legend range set according to the FOM computations.
  The ROM computations are performed using POD modes corresponding to $\mu_{\textrm{pod}}=1-10^{-6}$.}
  \label{fg:num_cyl_xr_x-0d65_sol}
\end{figure}

\begin{figure}[h]\centering
  \begin{subfigure}[t]{\textwidth}\centering
    \includegraphics[trim=0.0in 3.0in 0.0in 0.0in, clip, width=.48\textwidth]{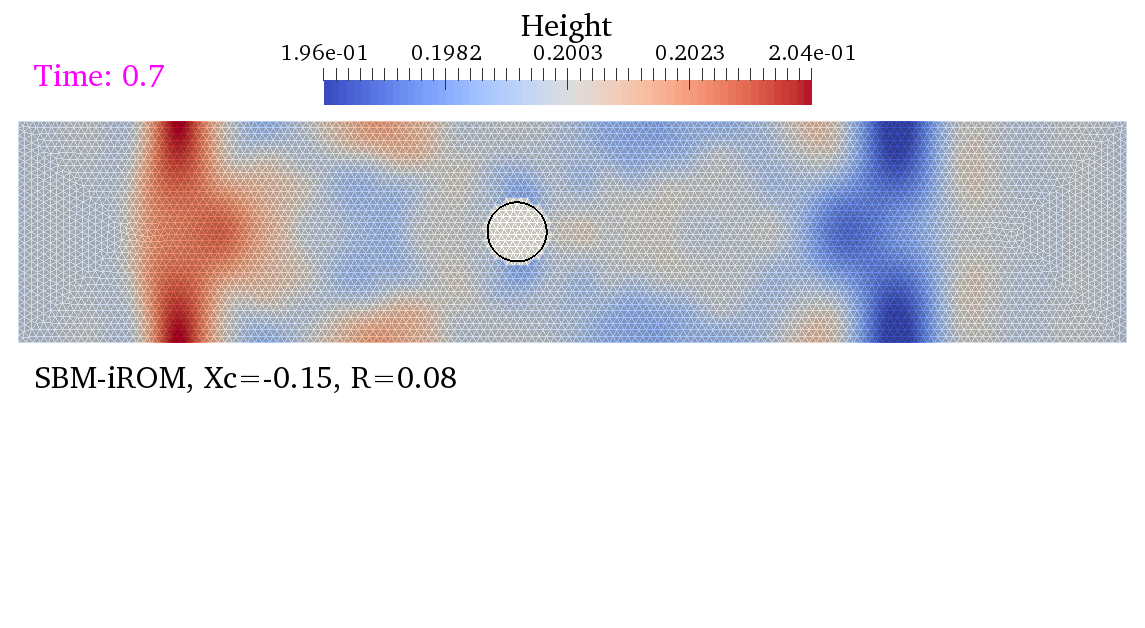}
    \includegraphics[trim=0.0in 3.0in 0.0in 0.0in, clip, width=.48\textwidth]{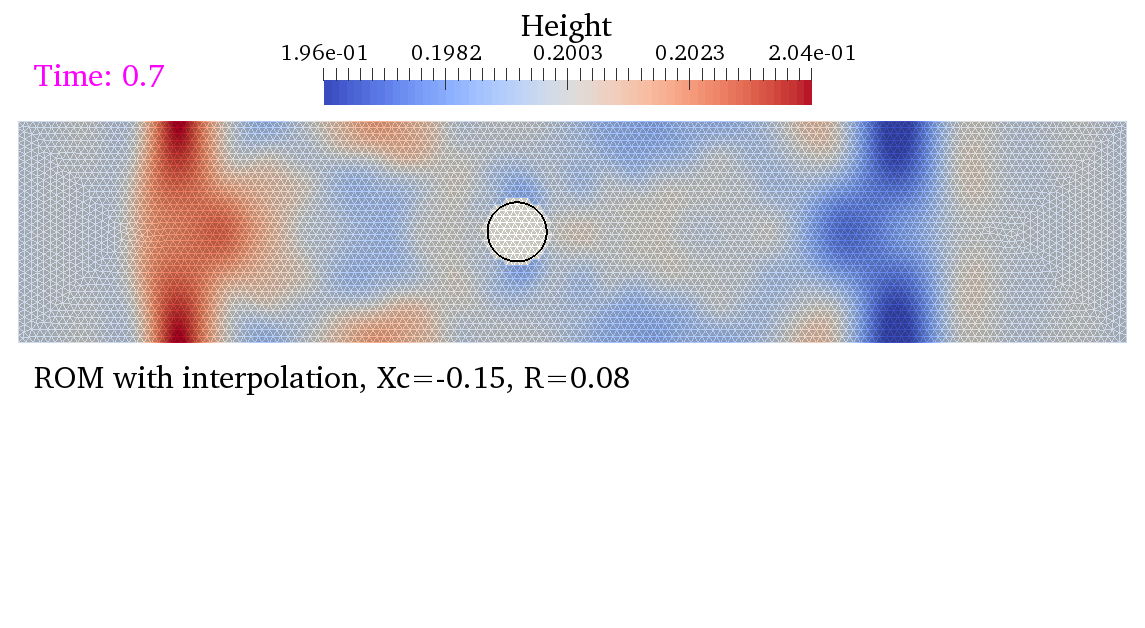}
    \caption{$R=0.08$.}
  \end{subfigure}
  \begin{subfigure}[t]{\textwidth}\centering
    \includegraphics[trim=0.0in 3.0in 0.0in 0.0in, clip, width=.48\textwidth]{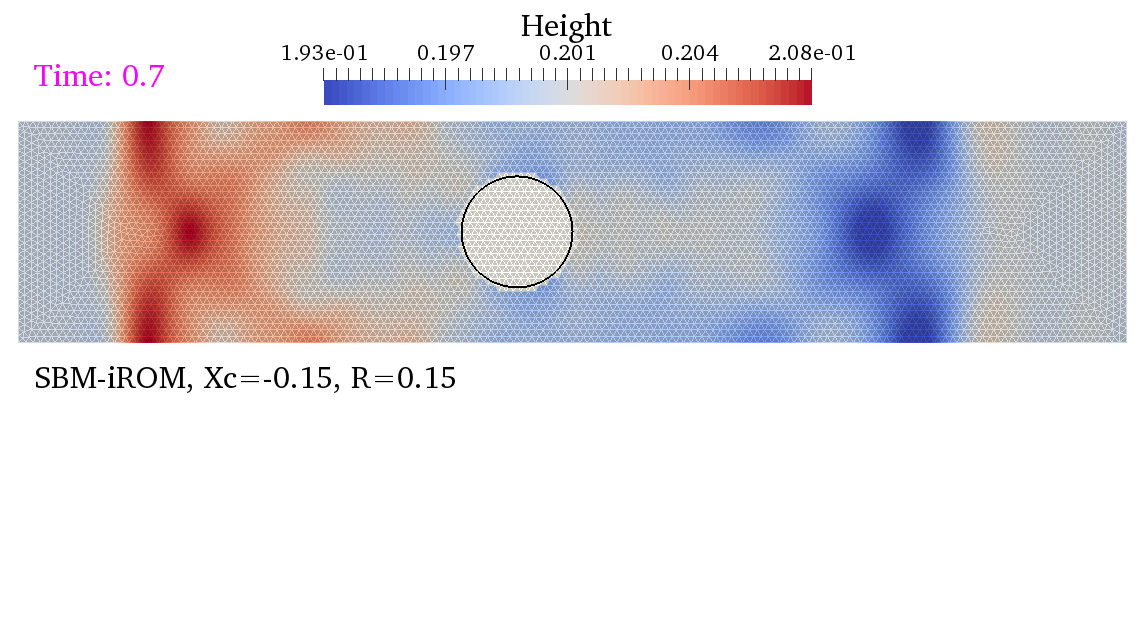}
    \includegraphics[trim=0.0in 3.0in 0.0in 0.0in, clip, width=.48\textwidth]{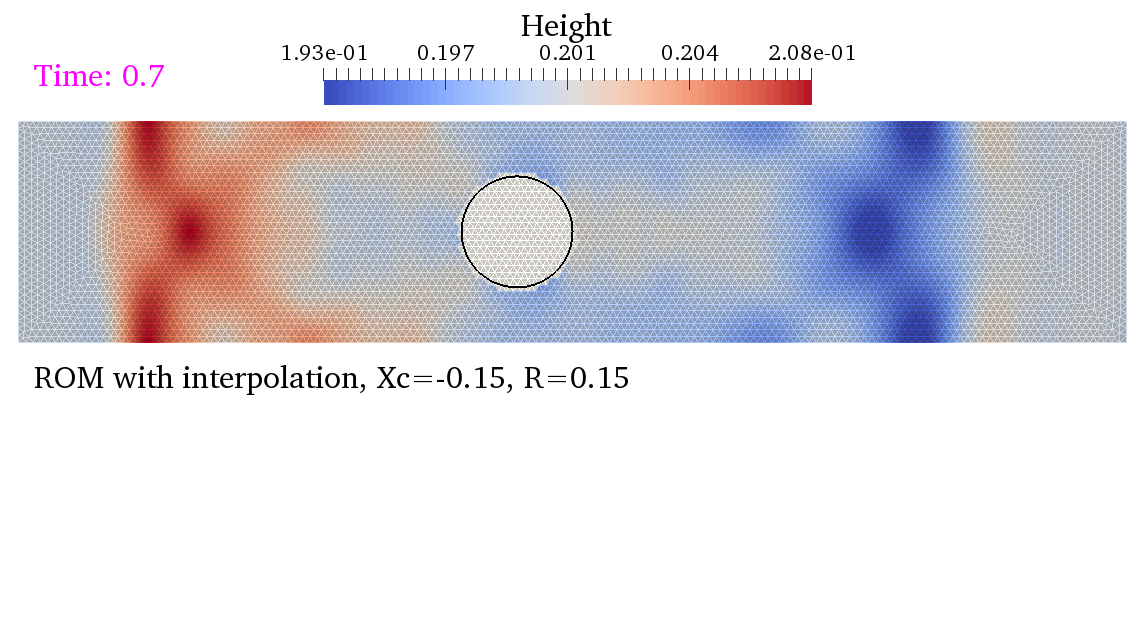}
    \caption{$R=0.15$.}
  \end{subfigure}
  \begin{subfigure}[t]{\textwidth}\centering
    \includegraphics[trim=0.0in 3.0in 0.0in 0.0in, clip, width=.48\textwidth]{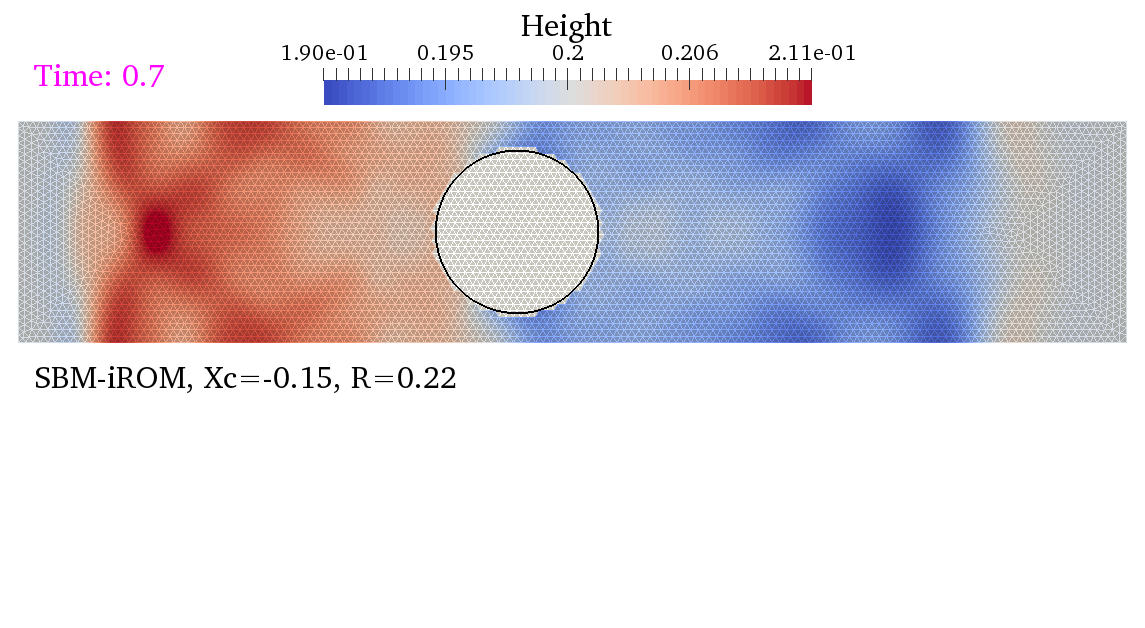}
    \includegraphics[trim=0.0in 3.0in 0.0in 0.0in, clip, width=.48\textwidth]{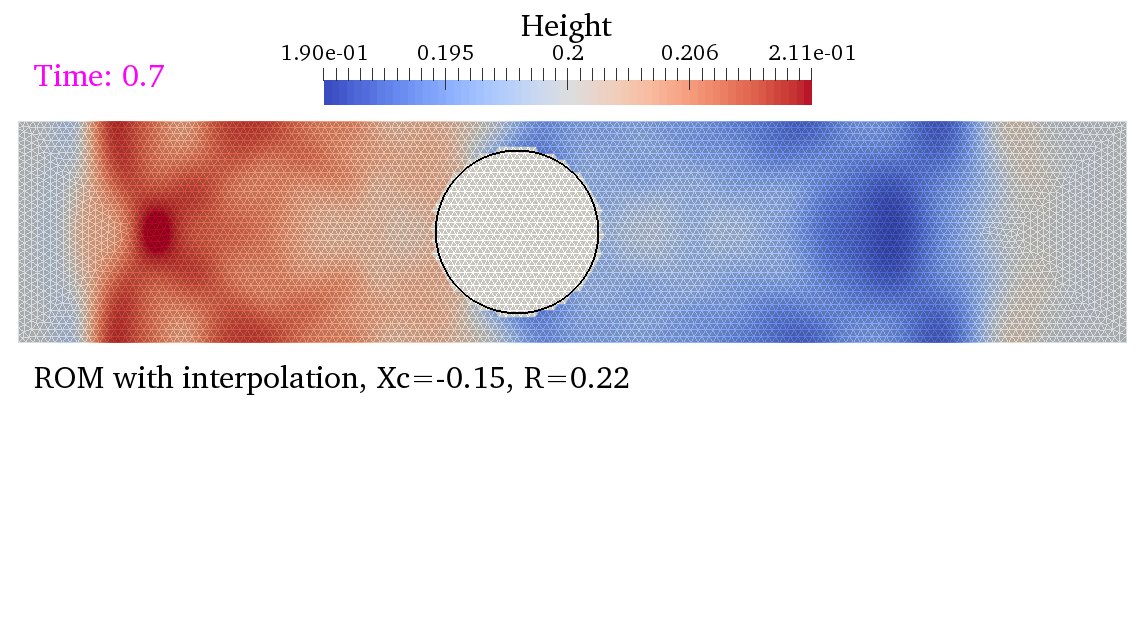}
    \caption{$R=0.22$.}
  \end{subfigure}
  \caption{The height ($h$) in the case of $R=0.08$ (top row), $R=0.15$ (middle row), or $R=0.22$ (bottom row), and $x_c=-0.15$ computed by FOM (left panels) and SBM-iROM (right panels), with the legend range set according to the FOM computations.
  The ROM computations are performed using POD modes corresponding to $\mu_{\textrm{pod}}=1-10^{-6}$.}
  \label{fg:num_cyl_xr_x-0d15_sol}
\end{figure}

\begin{figure}[h]\centering
  \begin{subfigure}[t]{\textwidth}\centering
    \includegraphics[trim=0.0in 3.0in 0.0in 0.0in, clip, width=.48\textwidth]{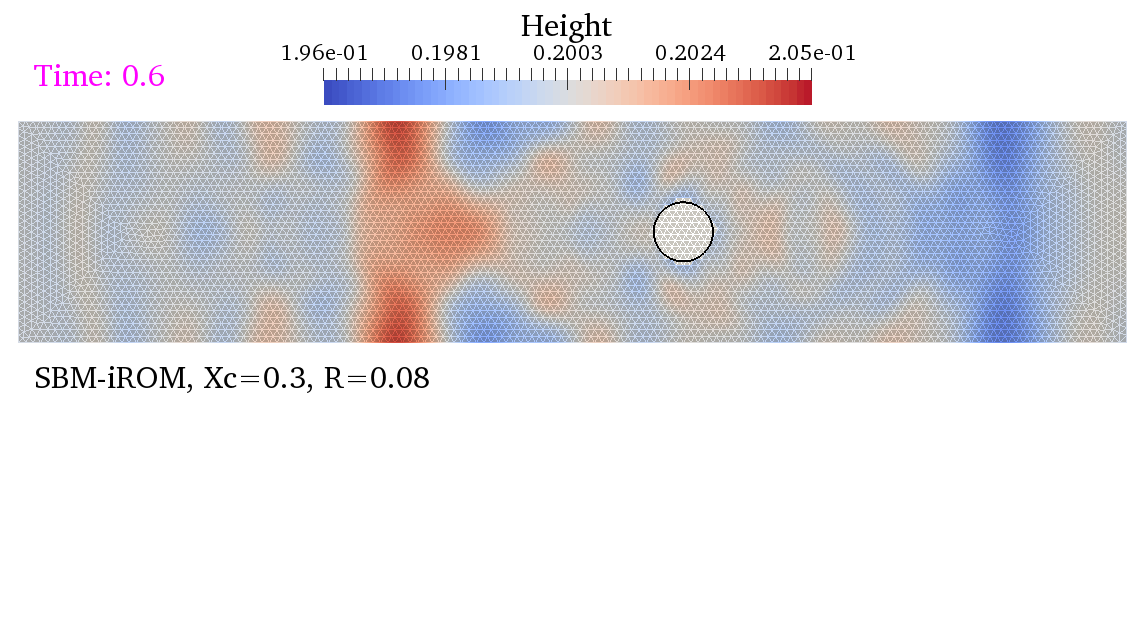}
    \includegraphics[trim=0.0in 3.0in 0.0in 0.0in, clip, width=.48\textwidth]{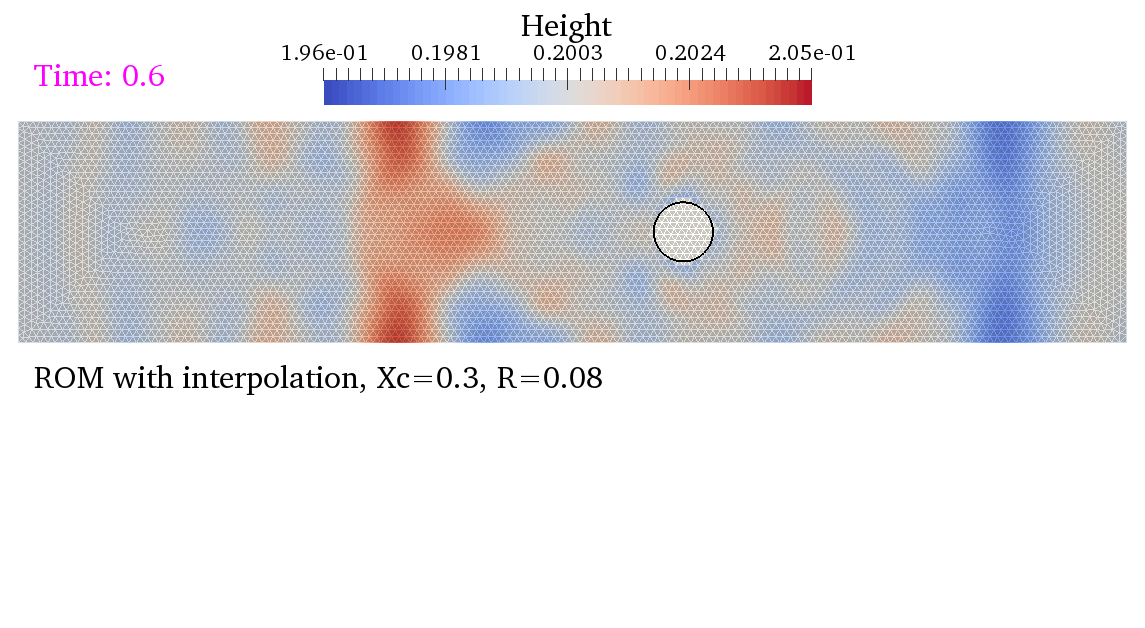}
    \caption{$R=0.08$.}
  \end{subfigure}
  \begin{subfigure}[t]{\textwidth}\centering
    \includegraphics[trim=0.0in 3.0in 0.0in 0.0in, clip, width=.48\textwidth]{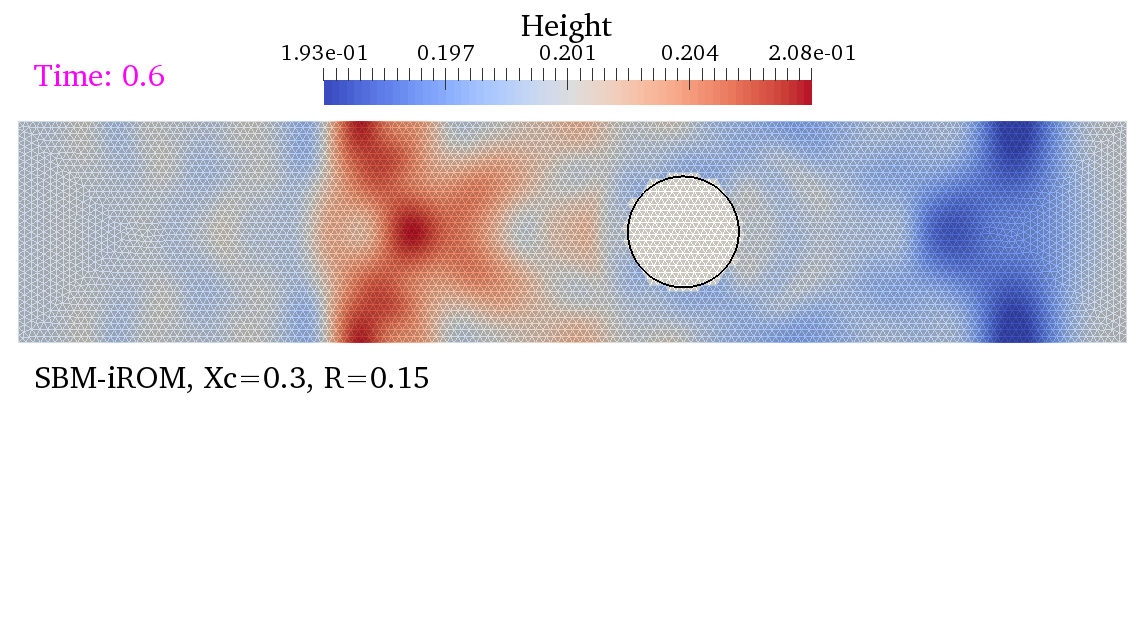}
    \includegraphics[trim=0.0in 3.0in 0.0in 0.0in, clip, width=.48\textwidth]{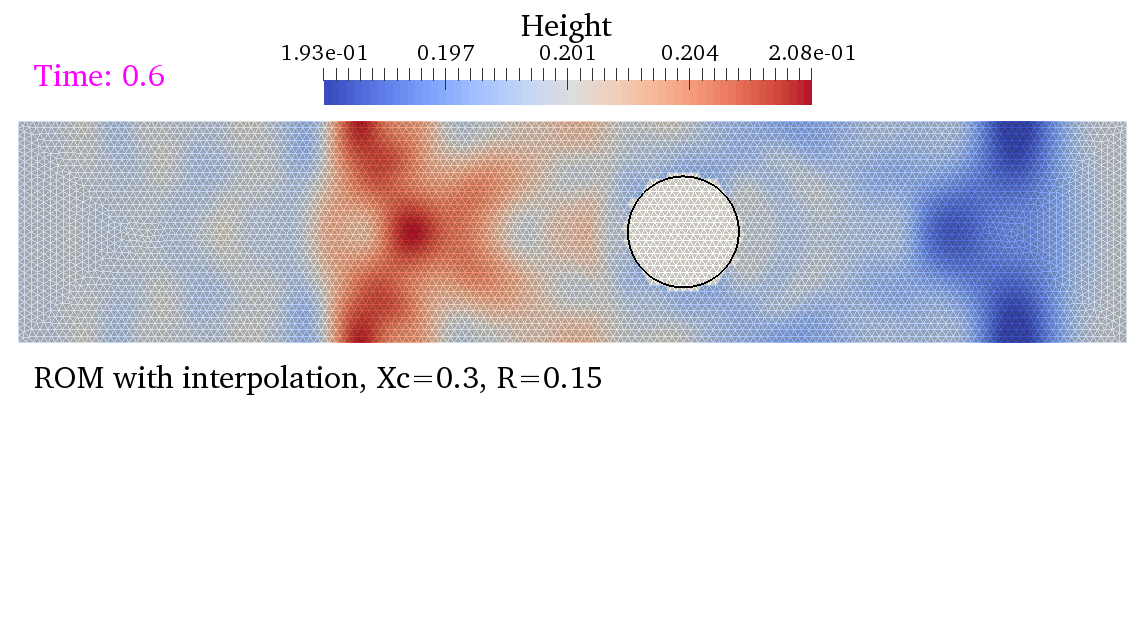}
    \caption{$R=0.15$.}
  \end{subfigure}
  \begin{subfigure}[t]{\textwidth}\centering
    \includegraphics[trim=0.0in 3.0in 0.0in 0.0in, clip, width=.48\textwidth]{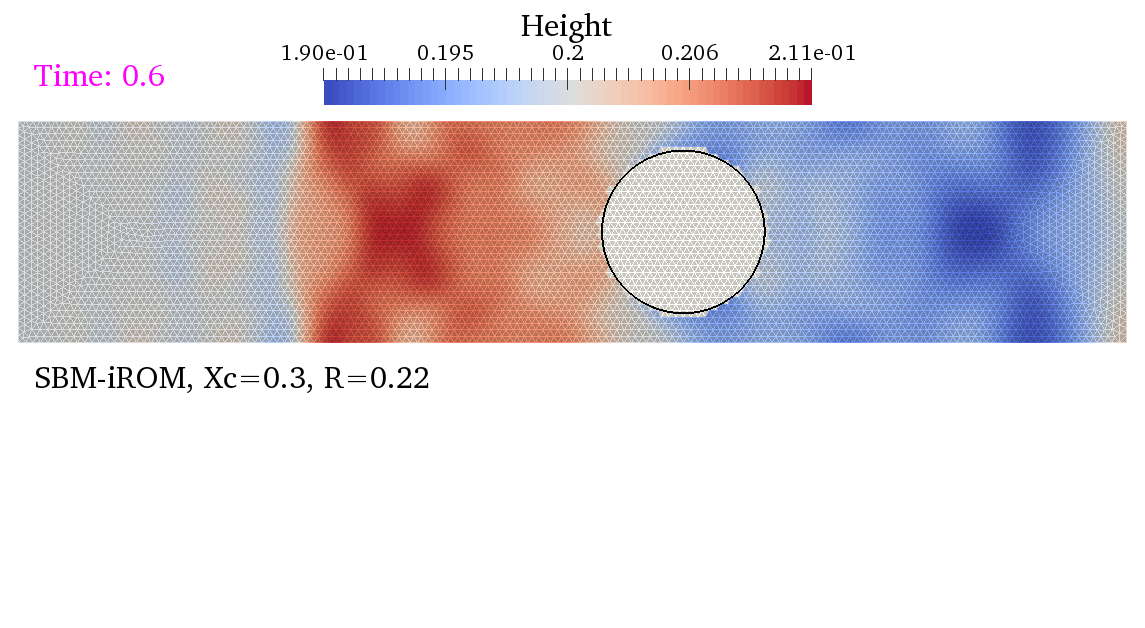}
    \includegraphics[trim=0.0in 3.0in 0.0in 0.0in, clip, width=.48\textwidth]{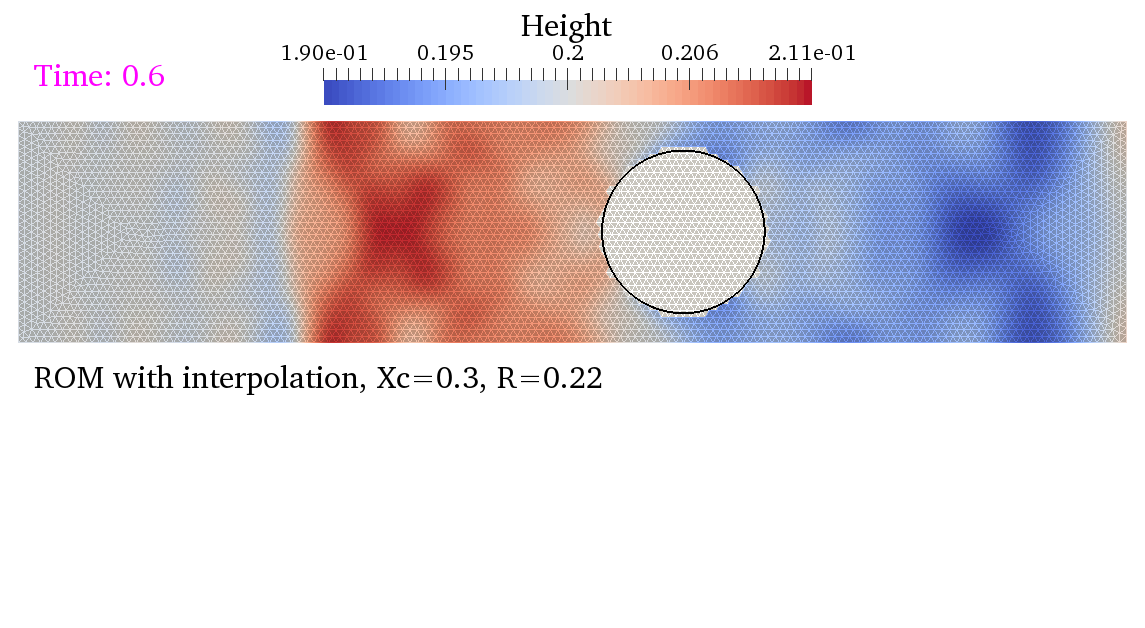}
    \caption{$R=0.22$.}
  \end{subfigure}
  \caption{The height ($h$) in the case of $R=0.08$ (top row), $R=0.15$ (middle row), or $R=0.22$ (bottom row), and $x_c=0.3$ computed by FOM (left panels) and SBM-iROM (right panels), with the legend range set according to the FOM computations.
  The ROM computations are performed using POD modes corresponding to $\mu_{\textrm{pod}}=1-10^{-6}$.}
  \label{fg:num_cyl_xr_x0d3_sol}
\end{figure}

\begin{figure}[h]\centering
  \begin{subfigure}[t]{\textwidth}\centering
    \includegraphics[trim=0.0in 3.0in 0.0in 0.0in, clip, width=.48\textwidth]{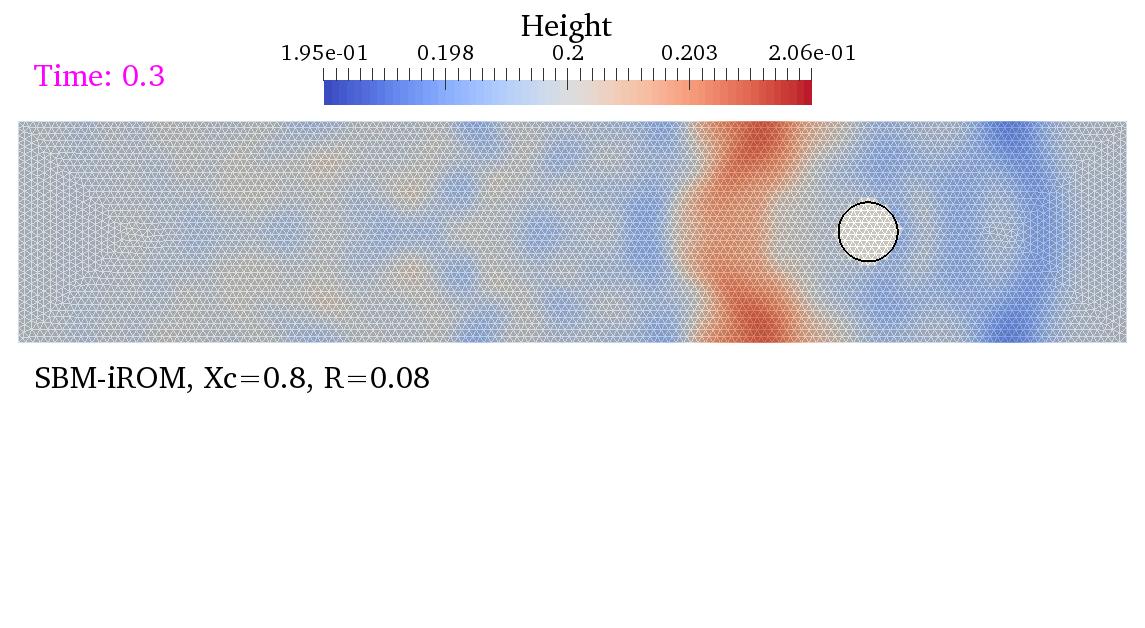}
    \includegraphics[trim=0.0in 3.0in 0.0in 0.0in, clip, width=.48\textwidth]{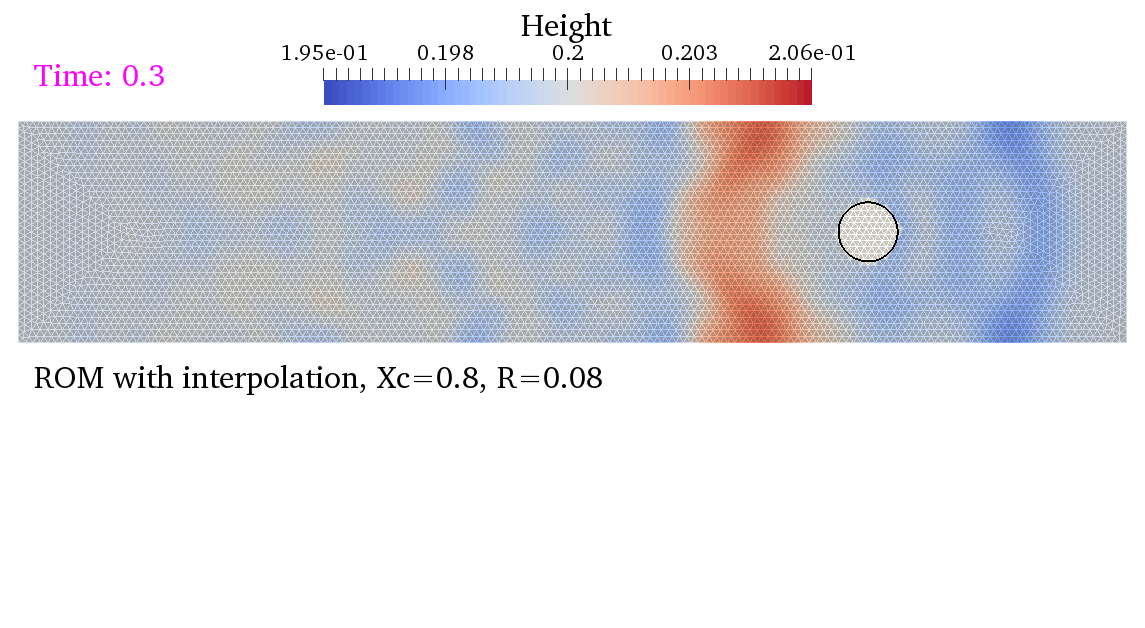}
    \caption{$R=0.08$.}
  \end{subfigure}
  \begin{subfigure}[t]{\textwidth}\centering
    \includegraphics[trim=0.0in 3.0in 0.0in 0.0in, clip, width=.48\textwidth]{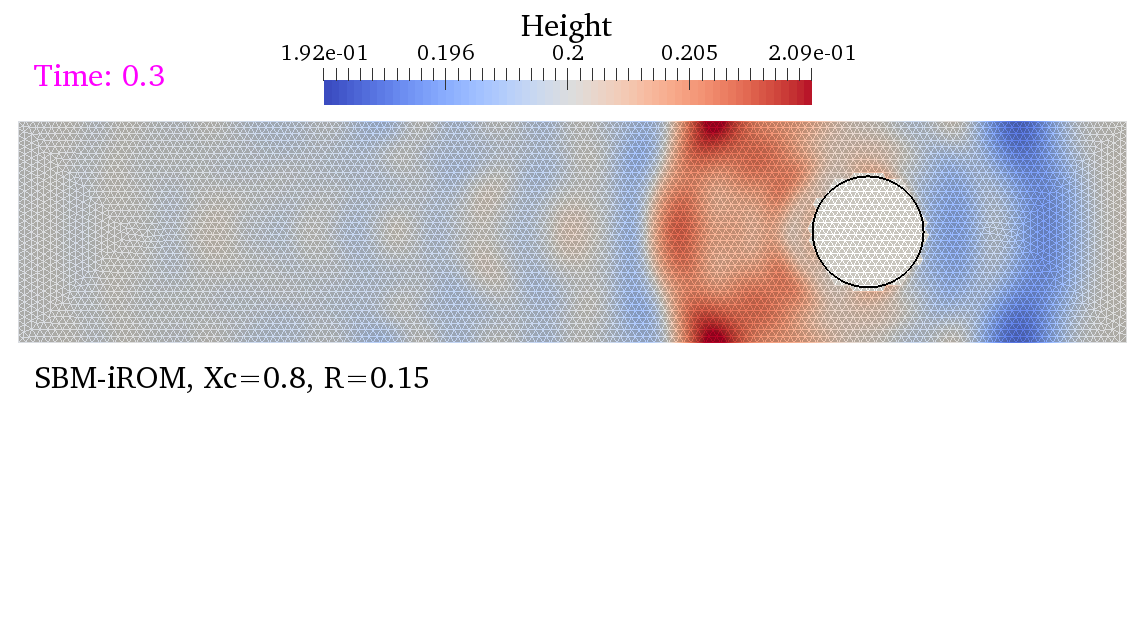}
    \includegraphics[trim=0.0in 3.0in 0.0in 0.0in, clip, width=.48\textwidth]{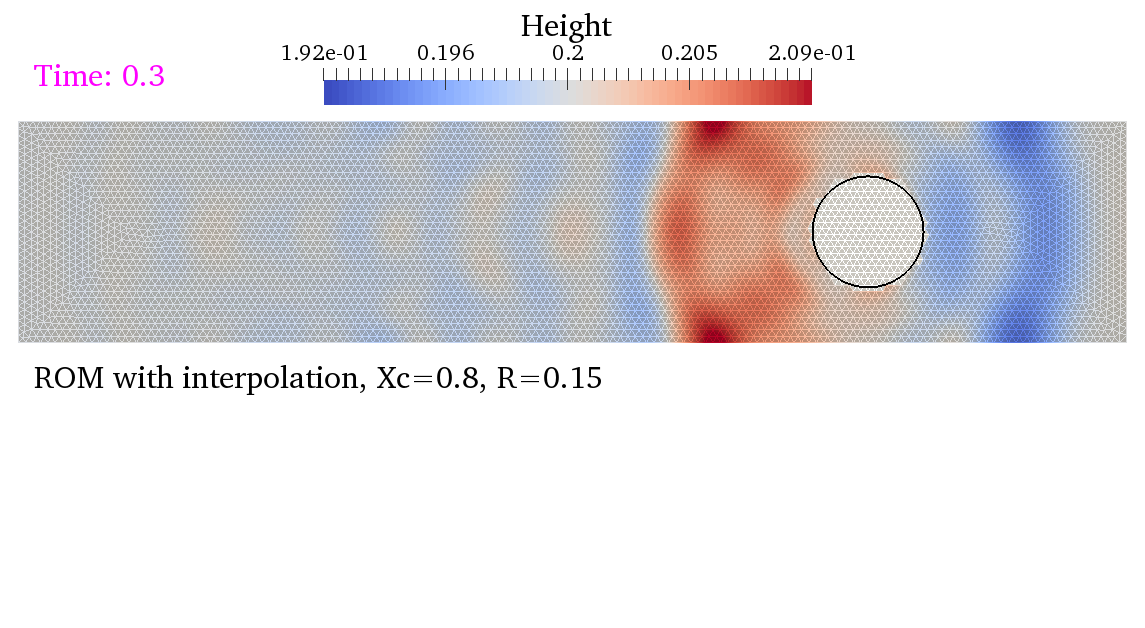}
    \caption{$R=0.15$.}
  \end{subfigure}
  \begin{subfigure}[t]{\textwidth}\centering
    \includegraphics[trim=0.0in 3.0in 0.0in 0.0in, clip, width=.48\textwidth]{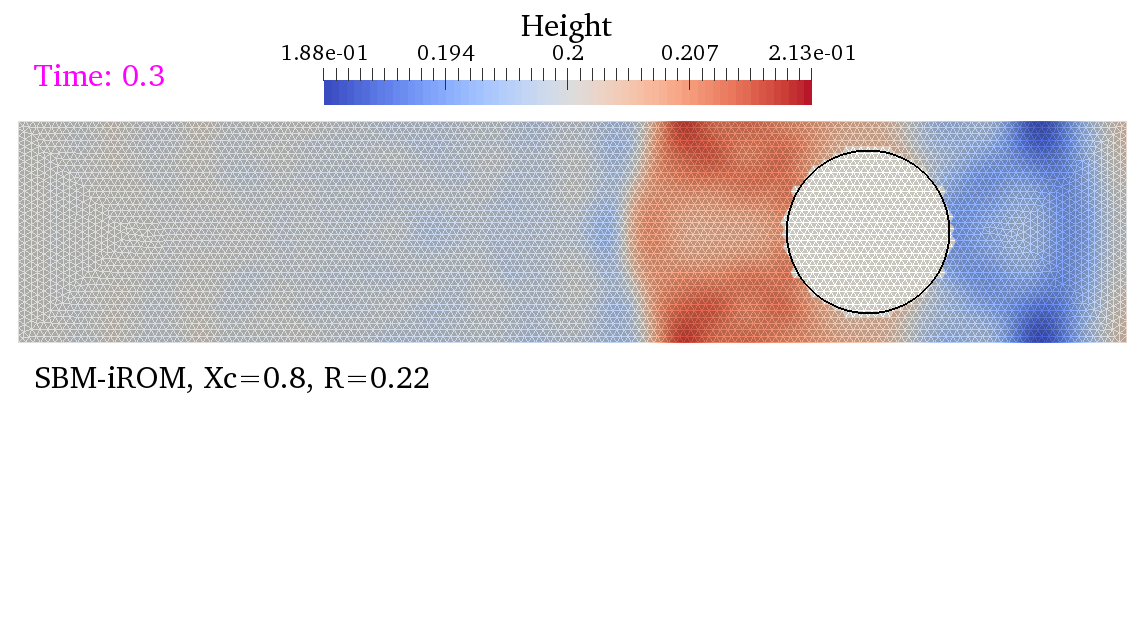}
    \includegraphics[trim=0.0in 3.0in 0.0in 0.0in, clip, width=.48\textwidth]{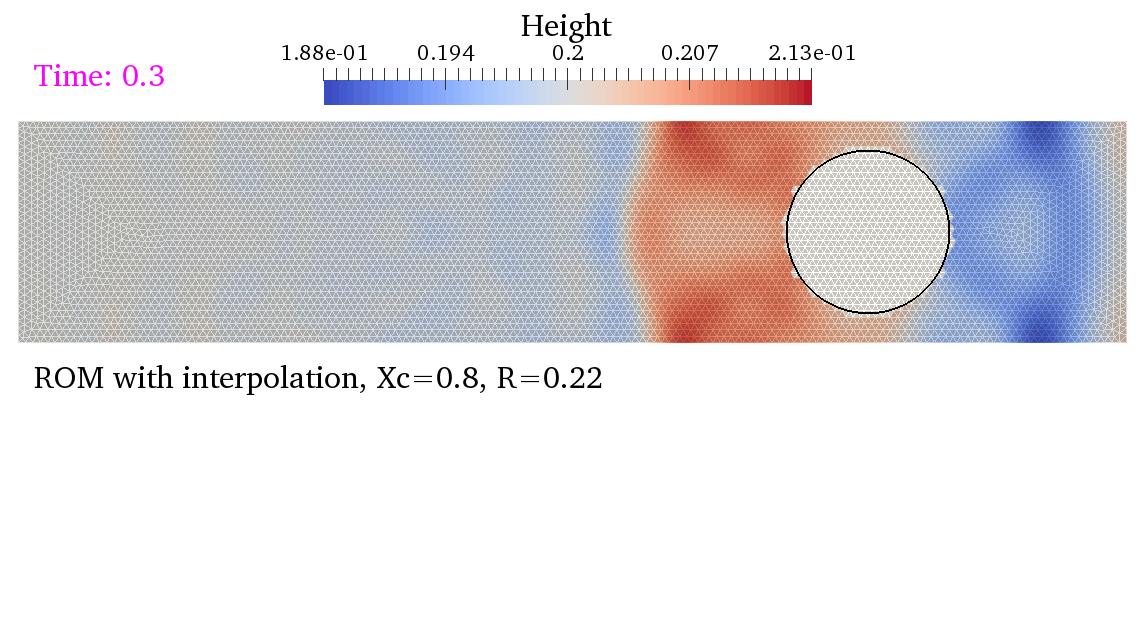}
    \caption{$R=0.22$.}
  \end{subfigure}
  \caption{The height ($h$) in the case of $R=0.08$ (top row), $R=0.15$ (middle row), or $R=0.22$ (bottom row), and $x_c=0.8$ computed by FOM (left panels) and SBM-iROM (right panels), with the legend range set according to the FOM computations.
  The ROM computations are performed using POD modes corresponding to $\mu_{\textrm{pod}}=1-10^{-6}$.}
  \label{fg:num_cyl_xr_x0d8_sol}
\end{figure}

In the last test case, we vary $R$ the same way as in Test 1 and $x_c$ the same way as in Test 2; hence there are 9 FOM computations with a total number of 1032 snapshots sampled with the frequency $n_{\freq}=10$.
The ROM computations also follow the same variation in $R$ and in $x_c$ as in previous two test cases.
In Figure~\ref{fg:num_cyl_xr_warp_r0d08_x-0d15} and Figure~\ref{fg:num_cyl_xr_warp_r0d15_x0d3}, we plot the height solution surfaces for the case $R=0.08$, $x_c=-0.15$ and $R=0.15$, $x_c=0.3$, respectively.
Note that as all nodes are active in some snapshots, SBM-ROM can handle $R=0.08$ (c.f. Test 1), but the solution is highly unstable.
The situation for SBM-ROM slightly improves when $R=0.15$; nevertheless, in both demonstrations ROM with interpolation is significant superior to that without interpolation.
Similar as in Test 2, we omit the solution plots for SBM-ROM, and compare the terminal water height solutions computed by FOM and SBM-iROM in Figures~\ref{fg:num_cyl_xr_x-0d65_sol}--\ref{fg:num_cyl_xr_x0d8_sol}.

\begin{figure}[h]\centering
  \begin{subfigure}[t]{\textwidth}\centering
    \includegraphics[width=.32\textwidth]{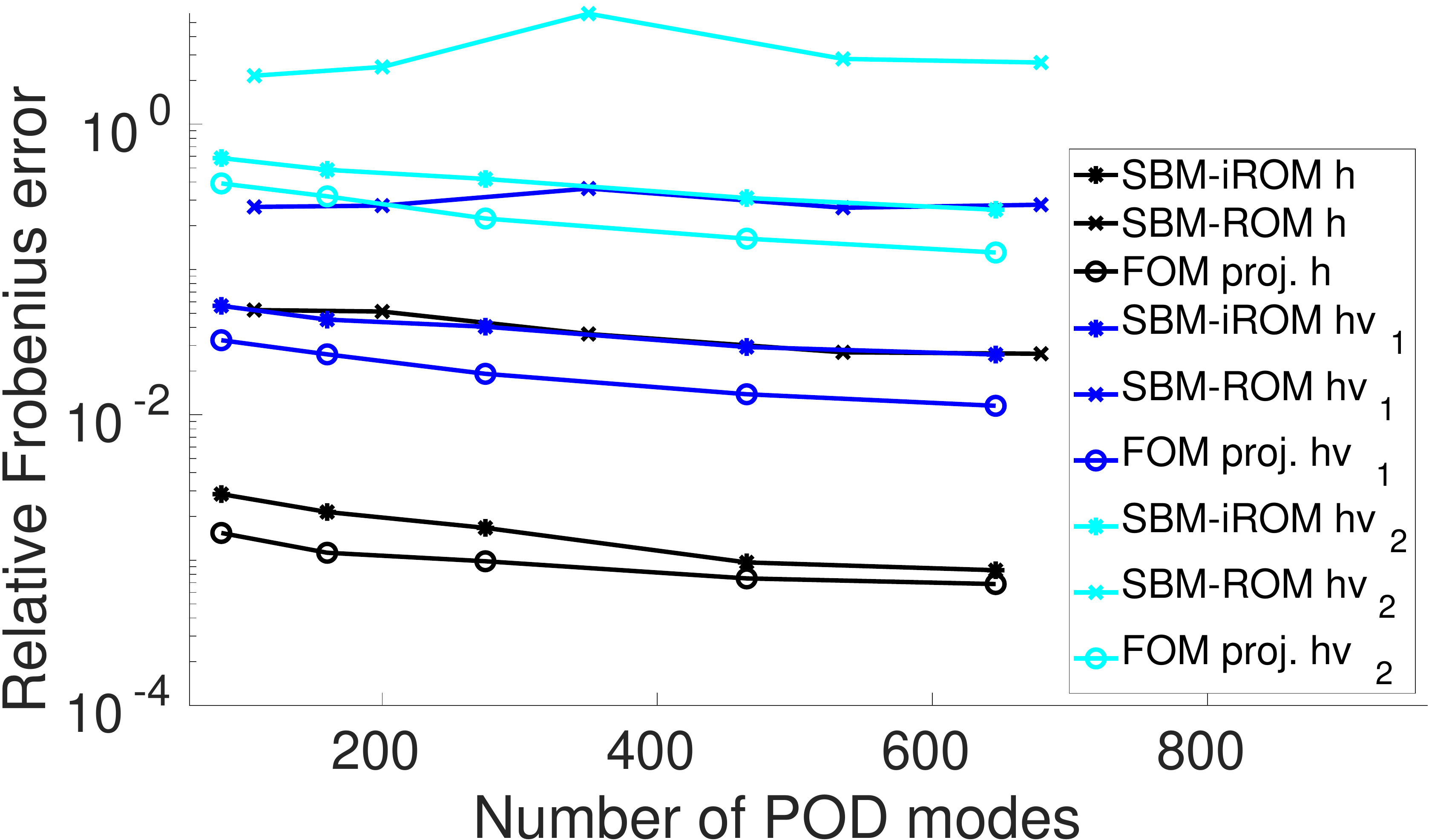}
    \includegraphics[width=.32\textwidth]{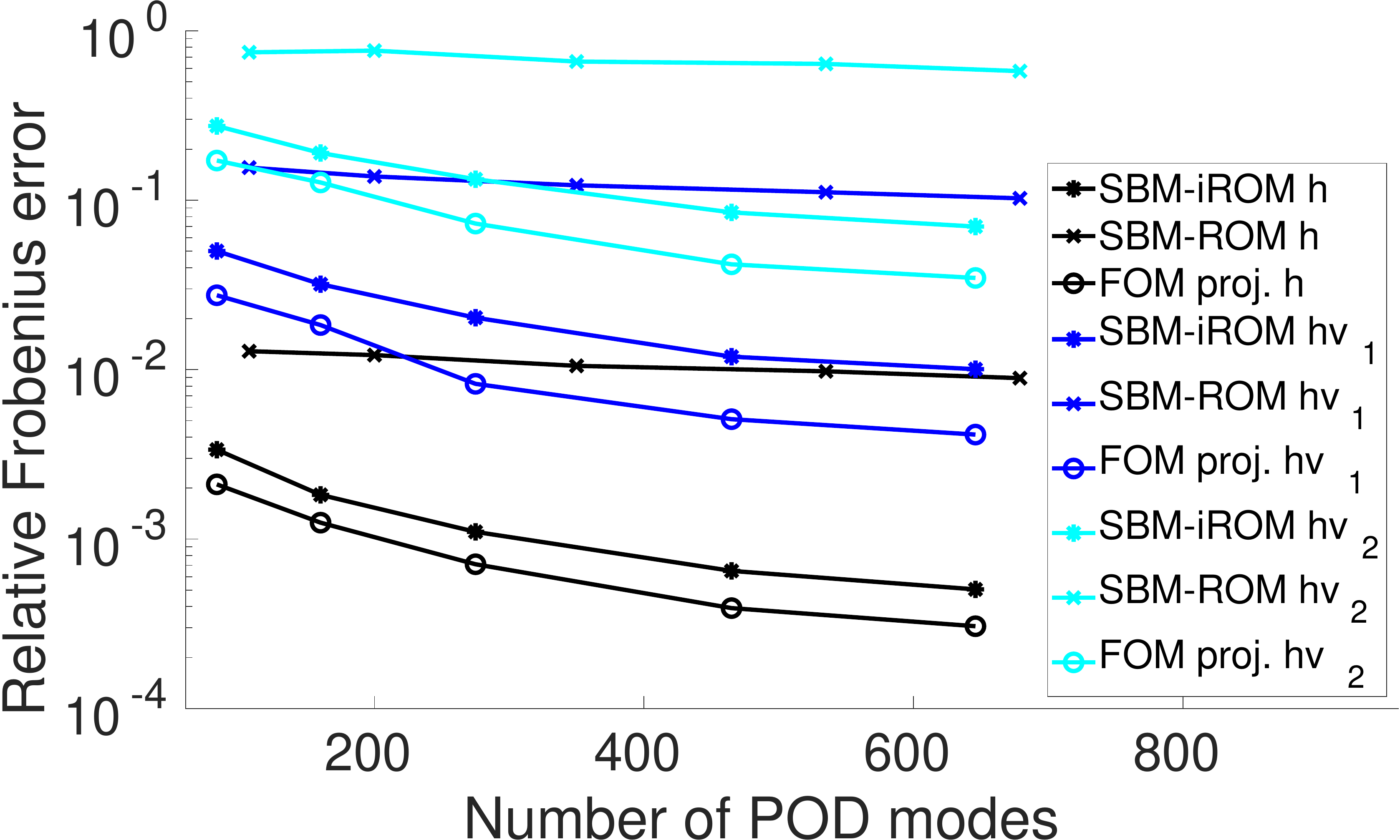}
    \includegraphics[width=.32\textwidth]{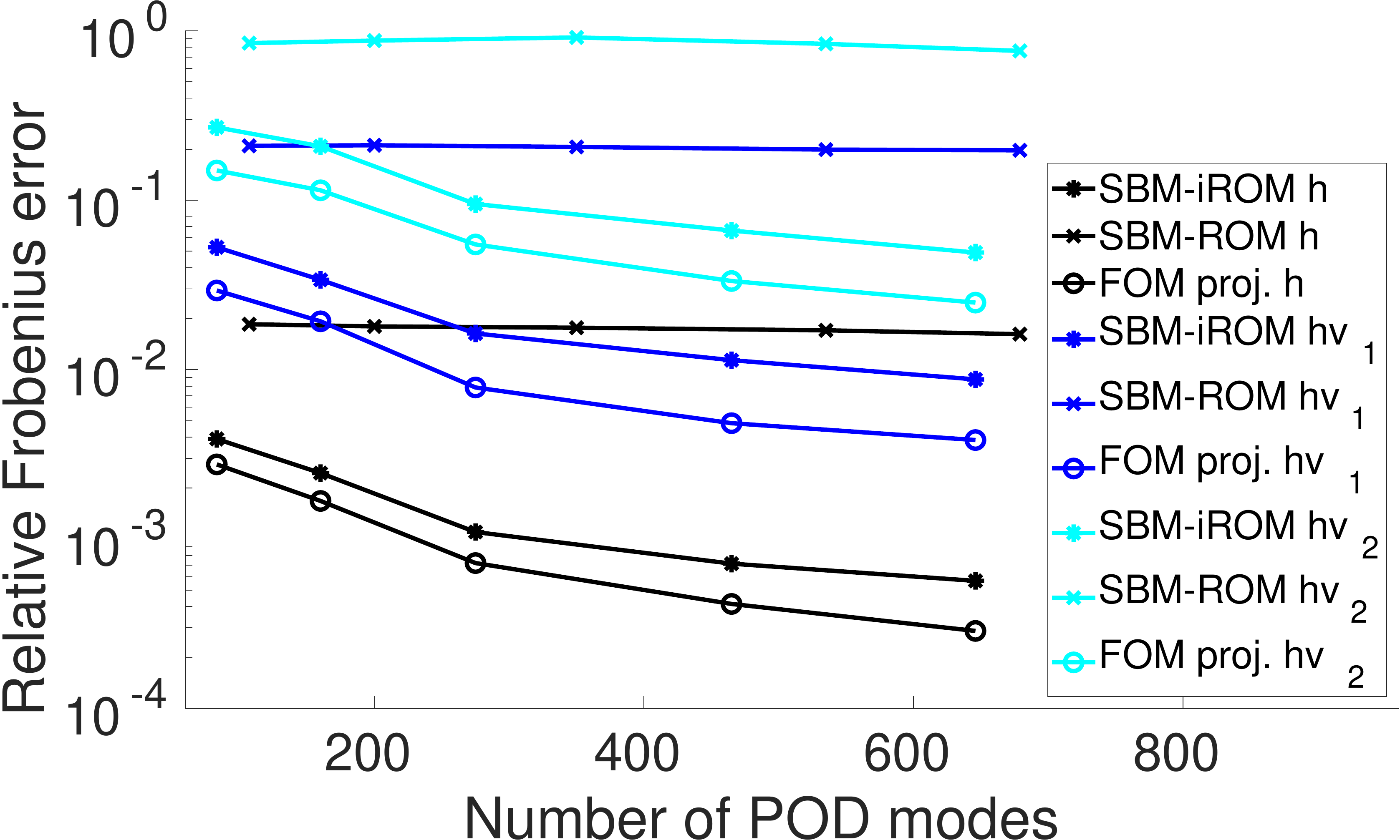}
    \caption{$x_c=-0.65$ with $R=0.08$ (left), $R=0.15$ (middle), and $R=0.22$ (right).}
    \label{fg:num_cyl_xr_froberr_x-0d65}
  \end{subfigure}
  \begin{subfigure}[t]{\textwidth}\centering
    \includegraphics[width=.32\textwidth]{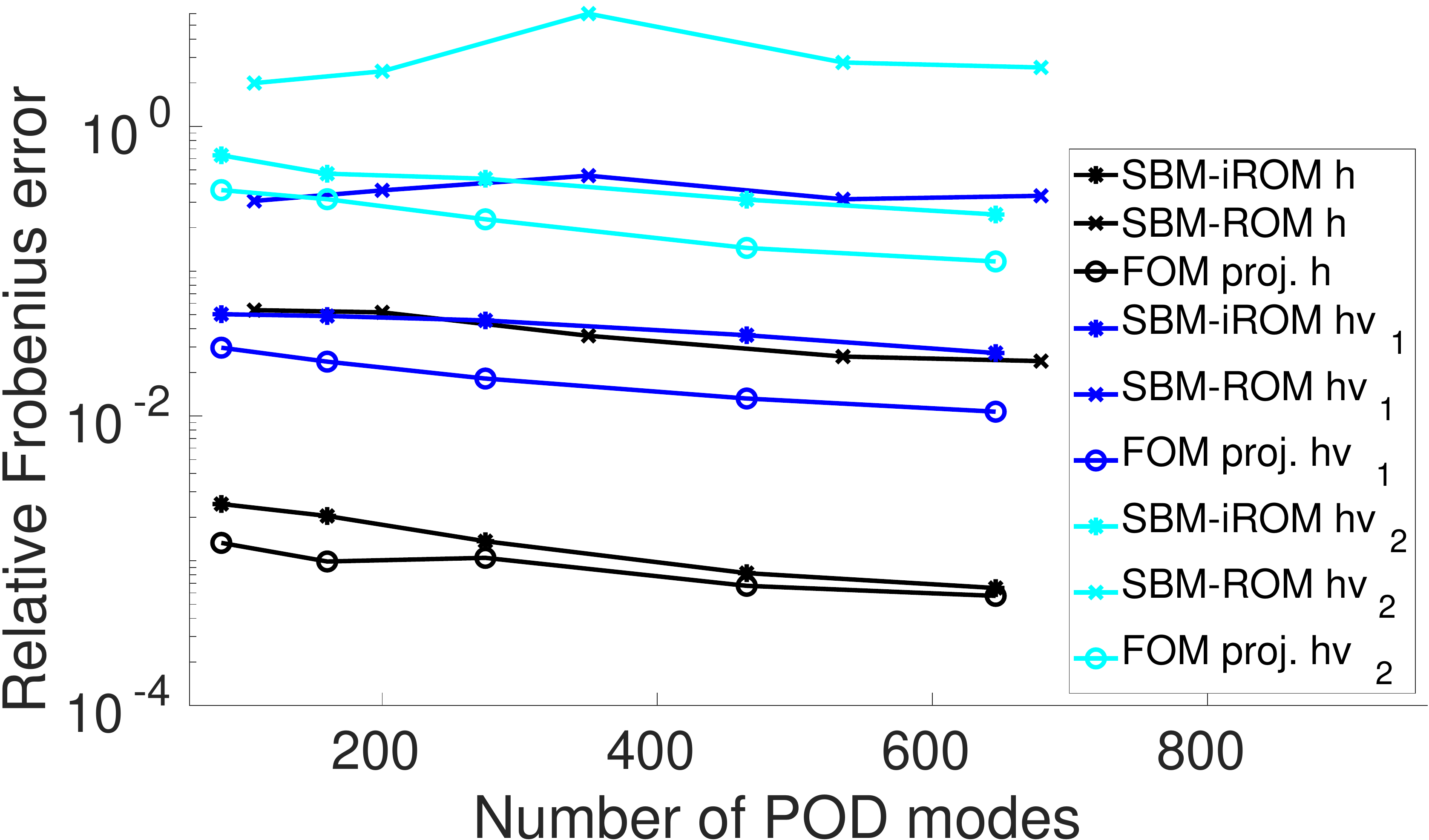}
    \includegraphics[width=.32\textwidth]{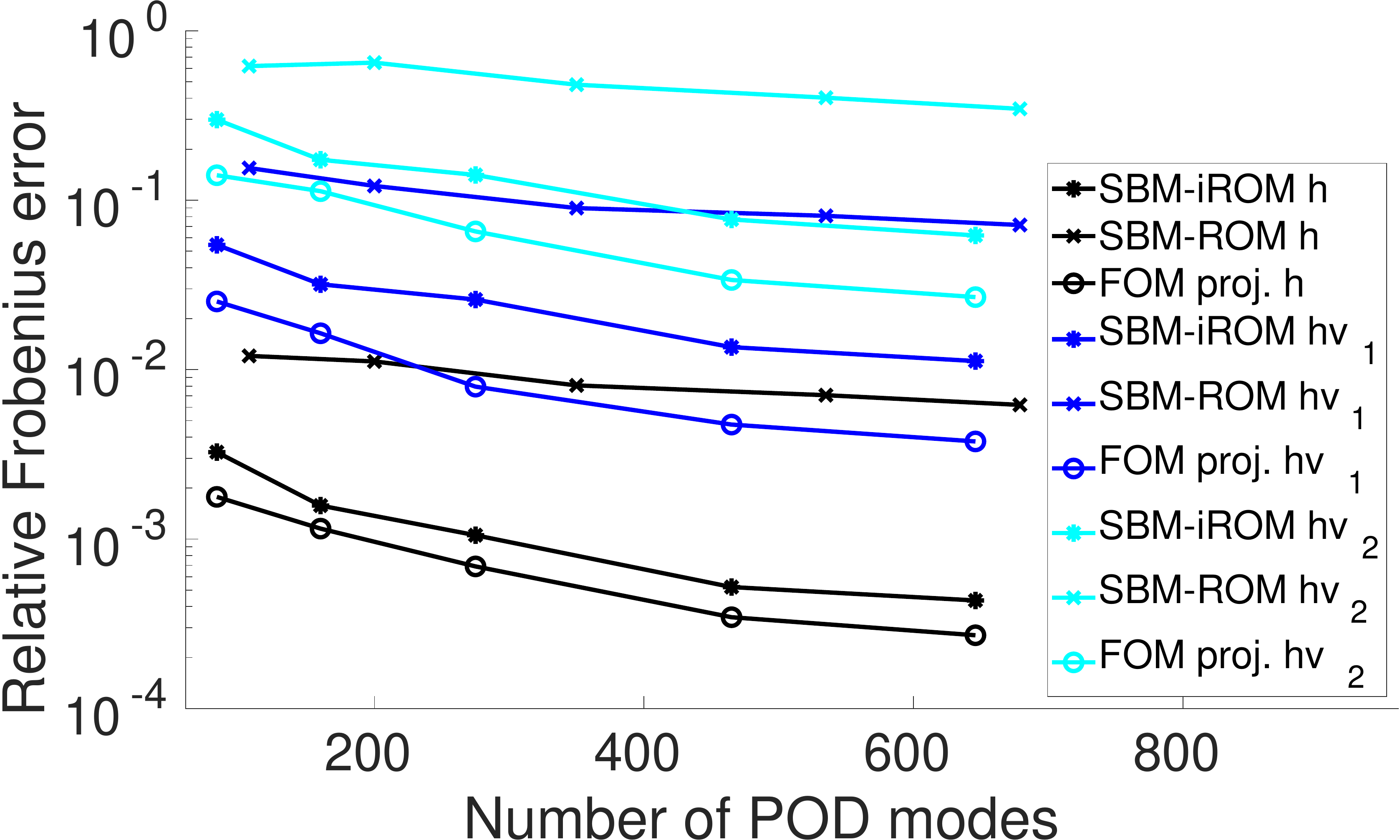}
    \includegraphics[width=.32\textwidth]{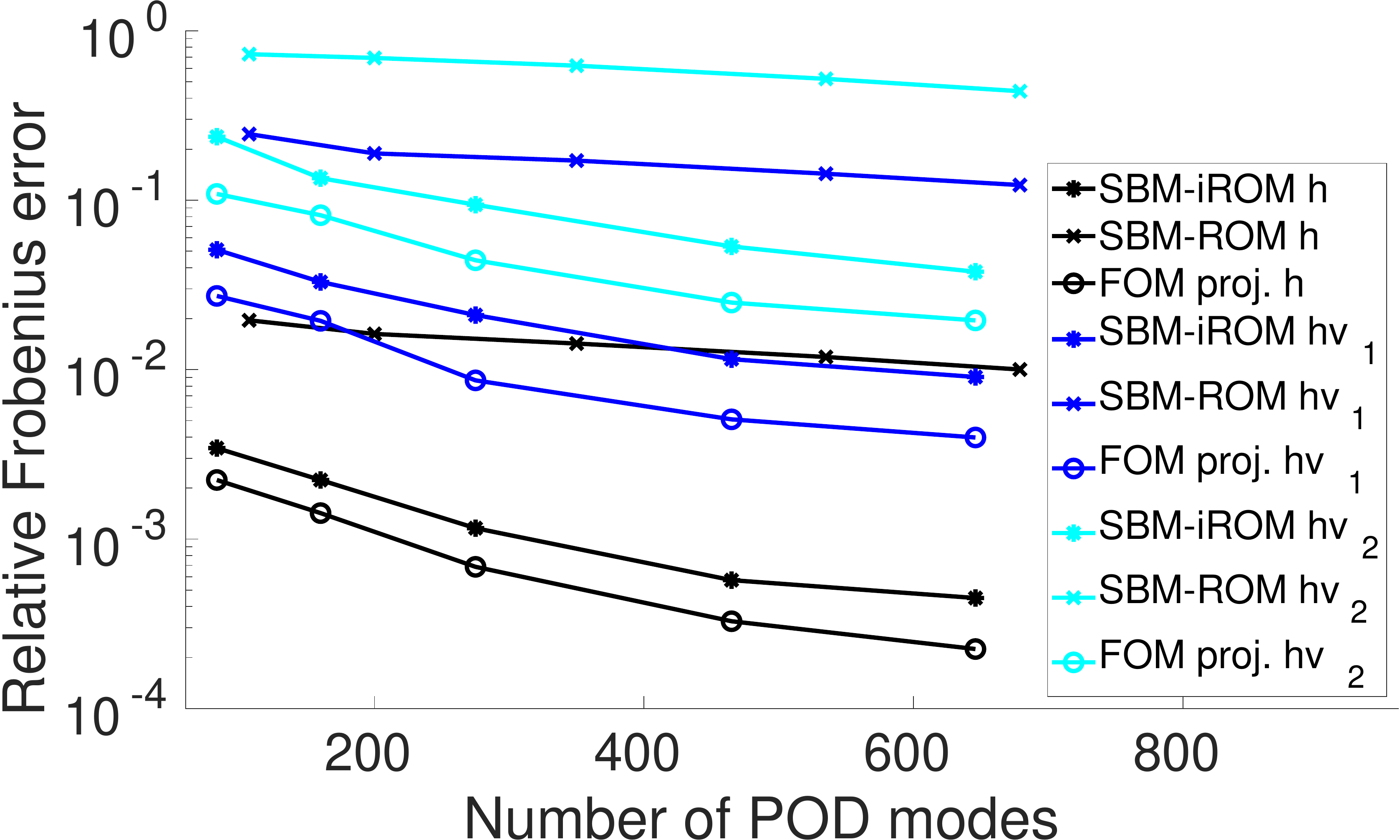}
    \caption{$x_c=-0.15$ with $R=0.08$ (left), $R=0.15$ (middle), and $R=0.22$ (right).}
    \label{fg:num_cyl_xr_froberr_x-0d15}
  \end{subfigure}
  \begin{subfigure}[t]{\textwidth}\centering
    \includegraphics[width=.32\textwidth]{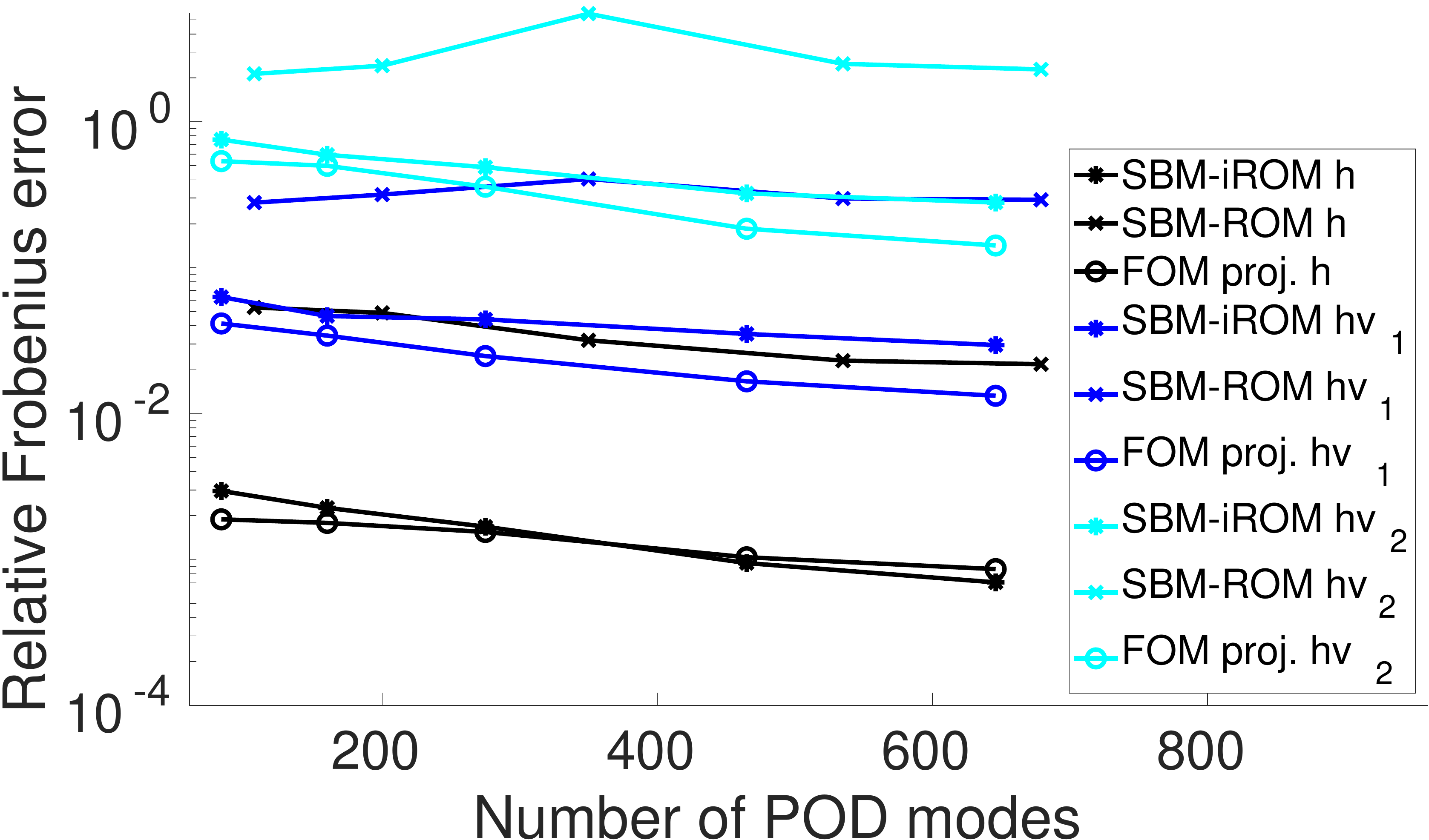}
    \includegraphics[width=.32\textwidth]{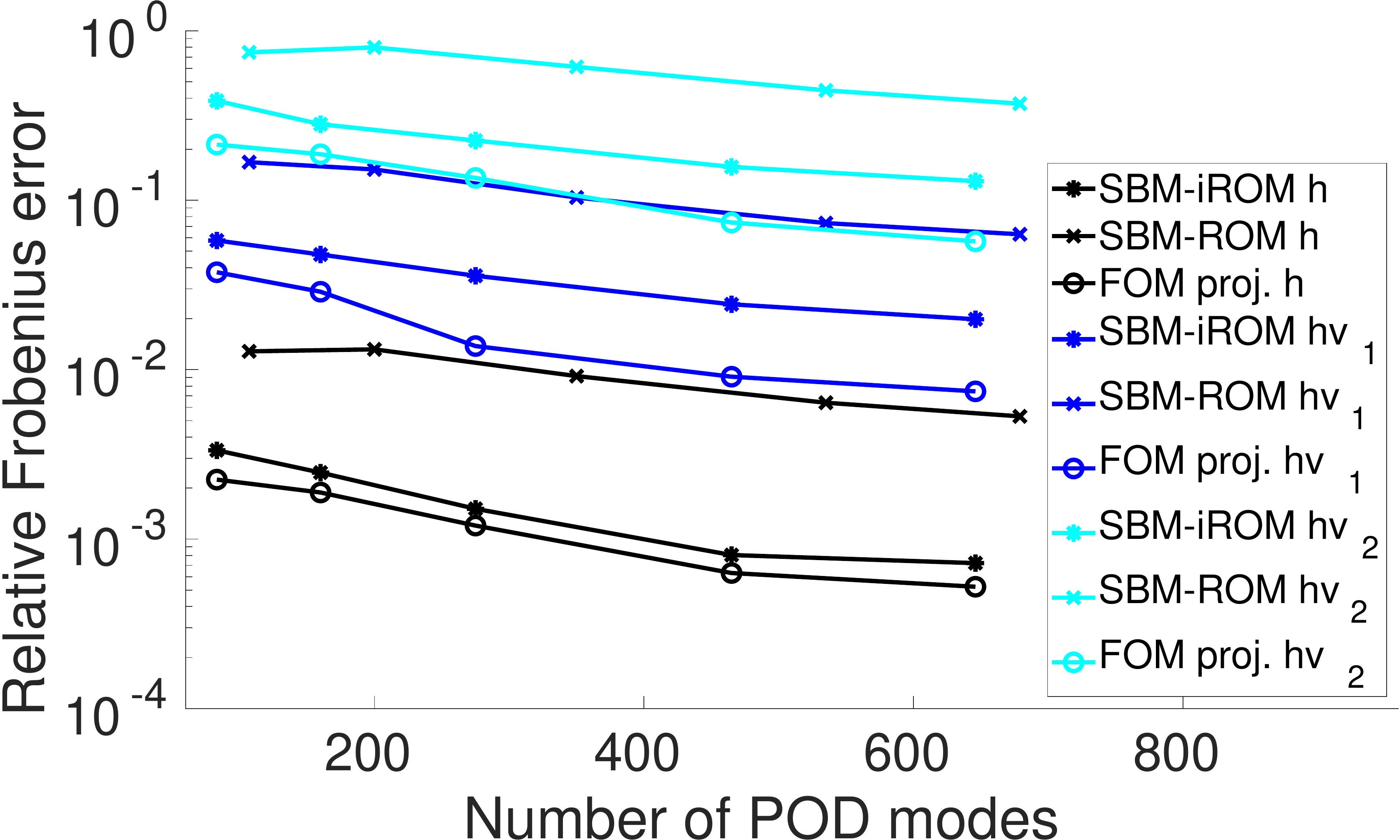}
    \includegraphics[width=.32\textwidth]{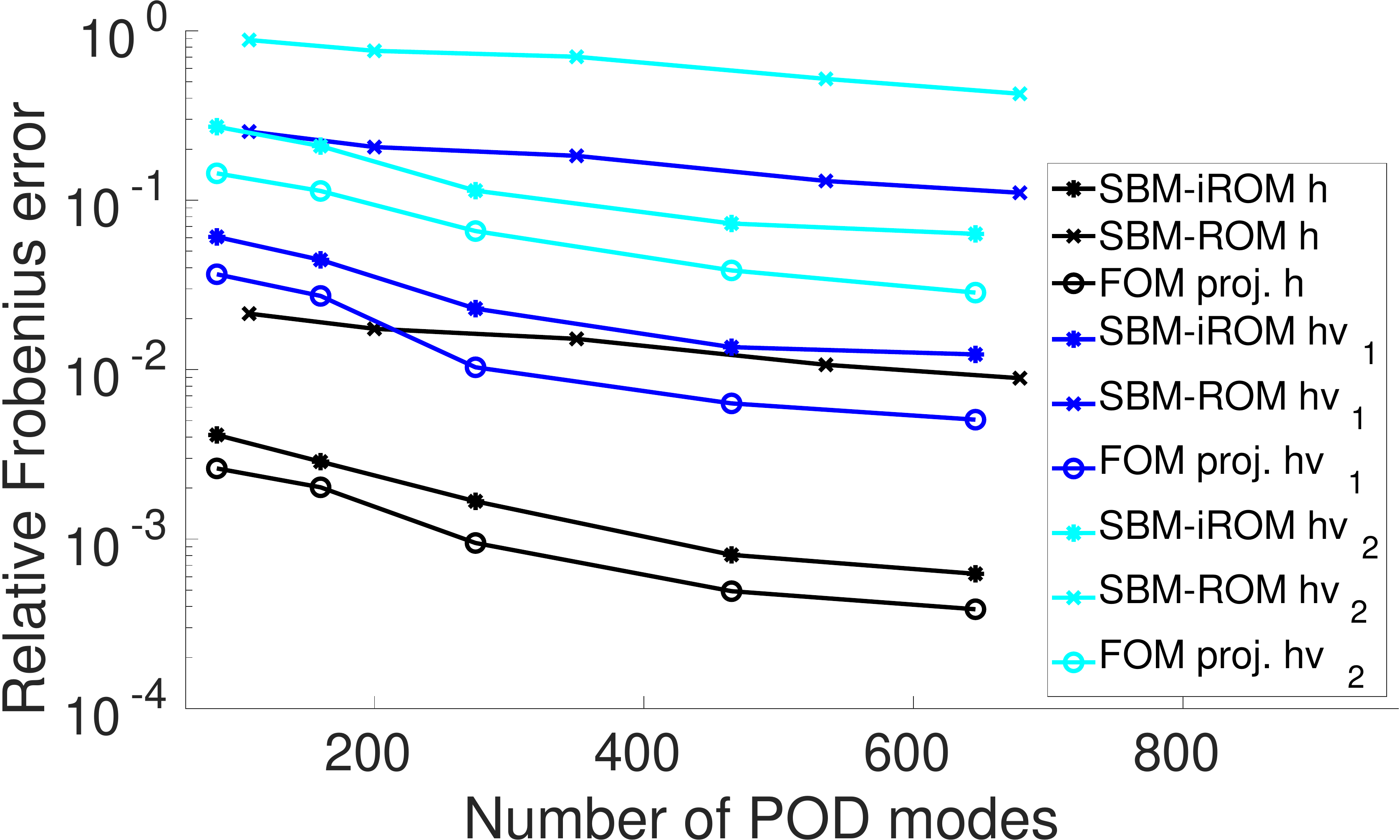}
    \caption{$x_c=0.3$ with $R=0.08$ (left), $R=0.15$ (middle), and $R=0.22$ (right).}
    \label{fg:num_cyl_xr_froberr_x0d3}
  \end{subfigure}
  \begin{subfigure}[t]{\textwidth}\centering
    \includegraphics[width=.32\textwidth]{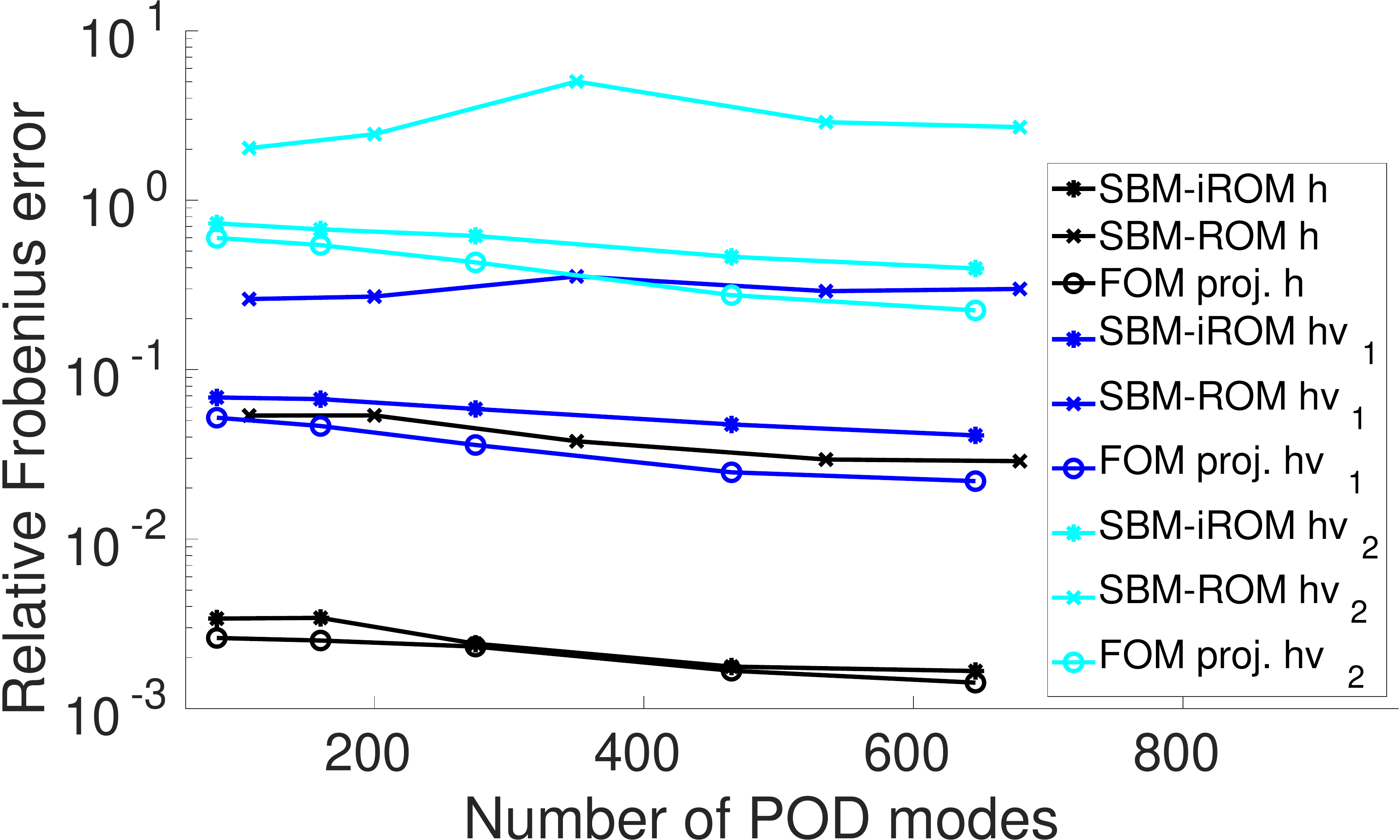}
    \includegraphics[width=.32\textwidth]{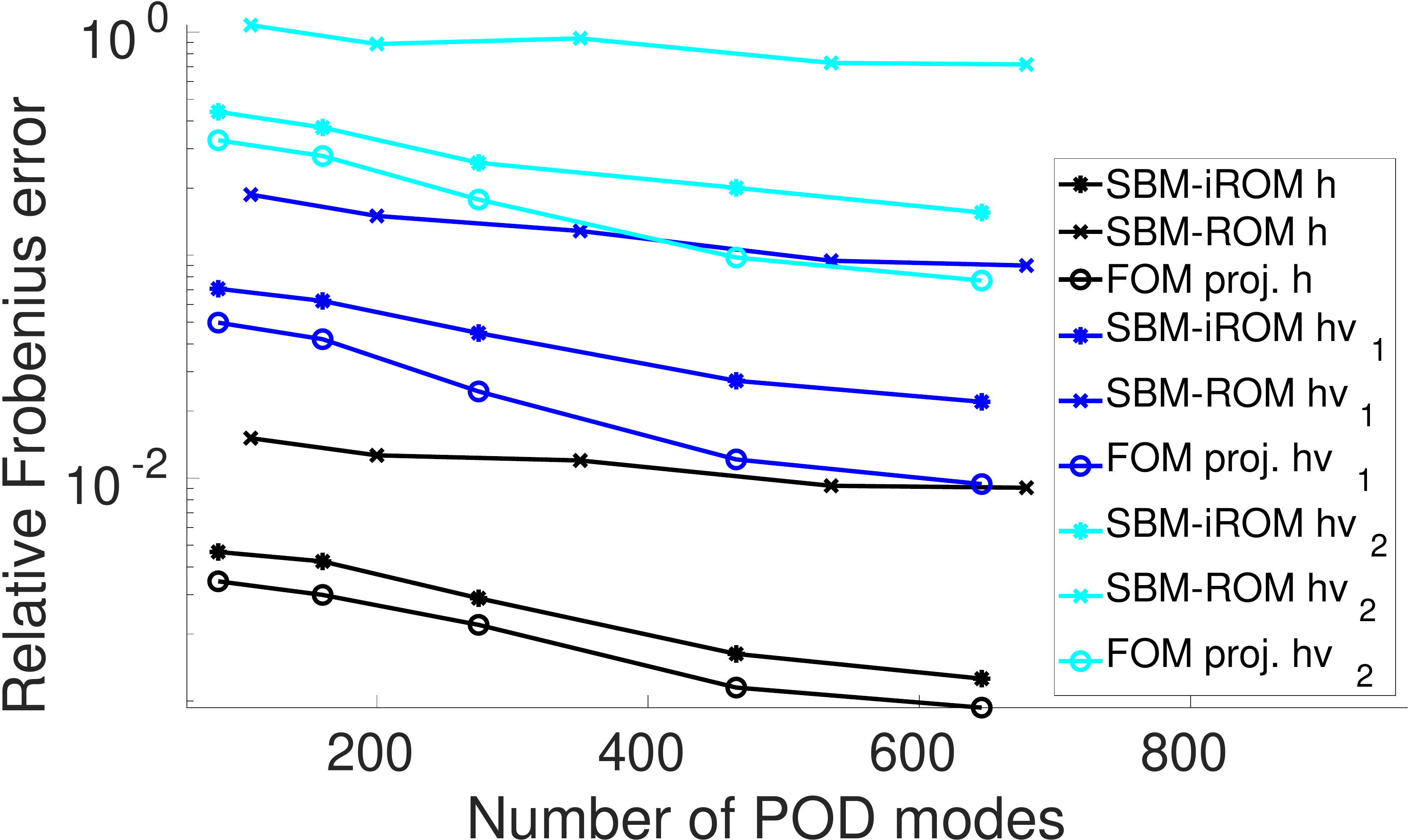}
    \includegraphics[width=.32\textwidth]{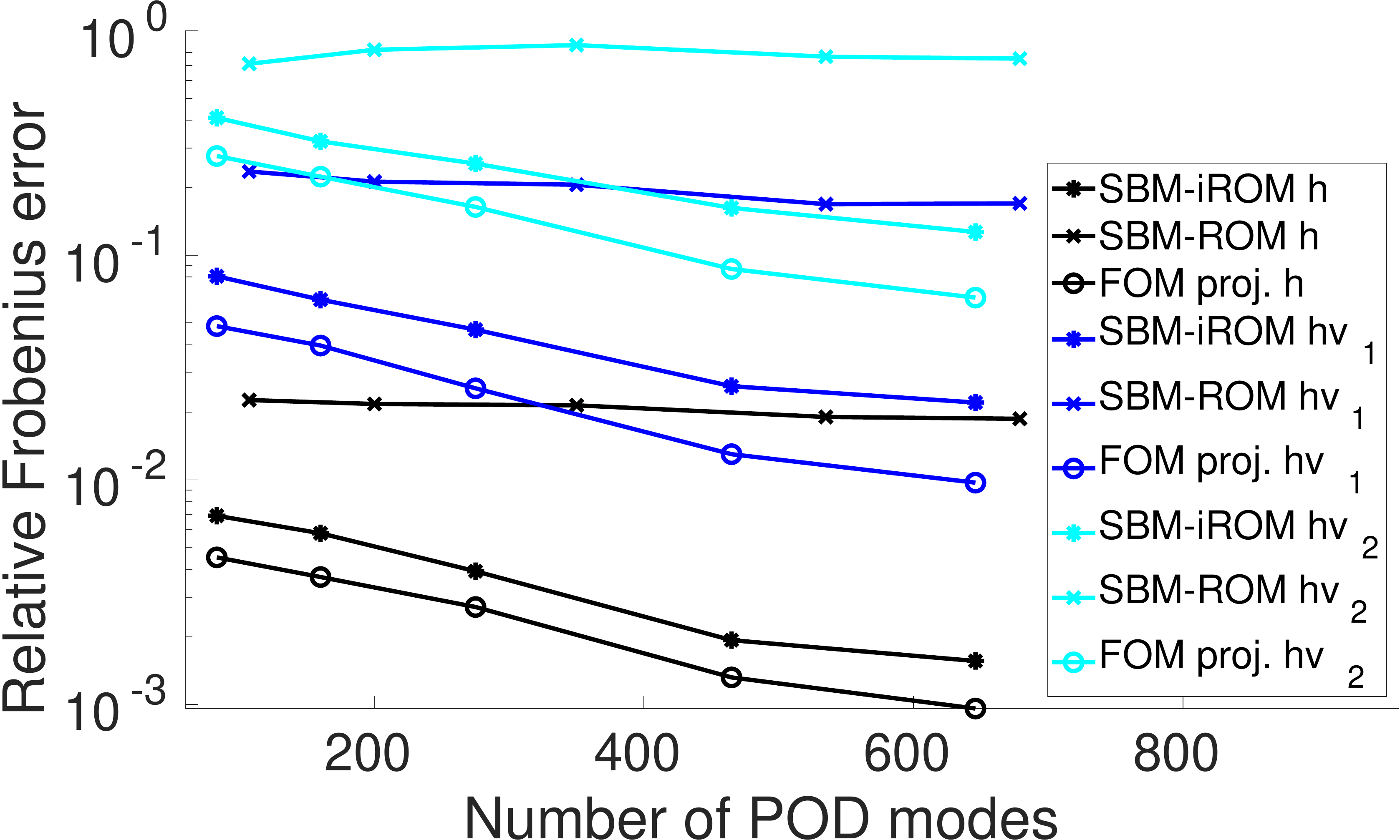}
    \caption{$x_c=0.8$ with $R=0.08$ (left), $R=0.15$ (middle), and $R=0.22$ (right).}
    \label{fg:num_cyl_xr_froberr_x0d8}
  \end{subfigure}
  \caption{The relative Frobenius errors of Test 3 by SBM-iROM (marked by $\ast$), SBM-ROM (marked by x), and FOM projection error (FOM proj., marked by o) in logarithmic scale plotted against the number of POD modes.
  The errors in $h$, $hv_1$, $hv_2$ are plotted in black, blue, and cyan colors, respectively.}
  \label{fg:num_cyl_xr_froberr}
\end{figure}

Lastly, we summarize the relative space-time Frobenius errors in Tables~\ref{tb:num_cyl_xr_froberr_r0d08_x-0d65}--\ref{tb:num_cyl_xr_froberr_r0d22_x0d8} and Figure~\ref{fg:num_cyl_xr_froberr}.
\begin{table}[h]\centering
  \caption{The relative space-time Frobenius errors of the projected FOM solutions and ROM computations in the case $R=0.08$ and $x_c=-0.65$ of Test 3.}
  \label{tb:num_cyl_xr_froberr_r0d08_x-0d65}
  \begin{tabular}{@{}lccccccccccc@{}}
  \toprule[.5mm]
  & \multicolumn{3}{l}{FOM projection} & & \multicolumn{3}{l}{SBM-iROM} & & \multicolumn{3}{l}{SBM-ROM} \\ \cmidrule[.2mm](lr){2-4} \cmidrule[.2mm](lr){6-8} \cmidrule[.2mm](l){10-12}
  $\mu_{\textrm{pod}}$ & $h$ & $hv_1$ & $hv_2$ & & $h$ & $hv_1$ & $hv_2$ & & $h$ & $hv_1$ & $hv_2$ \\ \cmidrule[.3mm](l){2-12}
  $1-10^{-5}$ & 1.54e-3 & 3.26e-2 & 3.91e-1 & & 2.85e-3 & 5.62e-2 & 5.84e-1 & & 5.26e-2 & 2.70e-1 & 2.16e+0 \\
  $1-10^{-6}$ & 1.12e-3 & 2.61e-2 & 3.19e-1 & & 2.14e-3 & 4.52e-2 & 4.85e-1 & & 5.14e-2 & 2.75e-1 & 2.47e+0 \\
  $1-10^{-7}$ & 9.83e-4 & 1.92e-2 & 2.24e-1 & & 1.67e-3 & 4.04e-2 & 4.20e-1 & & 3.60e-2 & 3.61e-1 & 5.76e+0 \\
  $1-10^{-8}$ & 7.49e-4 & 1.39e-2 & 1.63e-1 & & 9.65e-4 & 2.93e-2 & 3.10e-1 & & 2.69e-2 & 2.66e-1 & 2.81e+0 \\
  $1-10^{-9}$ & 6.85e-4 & 1.15e-2 & 1.31e-1 & & 8.53e-4 & 2.60e-2 & 2.58e-1 & & 2.63e-2 & 2.79e-1 & 2.65e+0 \\
  \bottomrule[.4mm]
  \end{tabular}
\end{table}
\begin{table}[h]\centering
  \caption{The relative space-time Frobenius errors of the projected FOM solutions and ROM computations in the case $R=0.15$ and $x_c=-0.65$ of Test 3.}
  \label{tb:num_cyl_xr_froberr_r0d15_x-0d65}
  \begin{tabular}{@{}lccccccccccc@{}}
  \toprule[.5mm]
  & \multicolumn{3}{l}{FOM projection} & & \multicolumn{3}{l}{SBM-iROM} & & \multicolumn{3}{l}{SBM-ROM} \\ \cmidrule[.2mm](lr){2-4} \cmidrule[.2mm](lr){6-8} \cmidrule[.2mm](l){10-12}
  $\mu_{\textrm{pod}}$ & $h$ & $hv_1$ & $hv_2$ & & $h$ & $hv_1$ & $hv_2$ & & $h$ & $hv_1$ & $hv_2$ \\ \cmidrule[.3mm](l){2-12}
  $1-10^{-5}$ & 2.10e-3 & 2.75e-2 & 1.71e-1 & & 3.37e-3 & 5.01e-2 & 2.74e-1 & & 1.28e-2 & 1.56e-1 & 7.46e-1 \\
  $1-10^{-6}$ & 1.25e-3 & 1.83e-2 & 1.27e-1 & & 1.82e-3 & 3.19e-2 & 1.90e-1 & & 1.22e-2 & 1.38e-1 & 7.64e-1 \\
  $1-10^{-7}$ & 7.10e-4 & 8.25e-3 & 7.28e-2 & & 1.10e-3 & 2.03e-2 & 1.33e-1 & & 1.05e-2 & 1.22e-1 & 6.58e-1 \\
  $1-10^{-8}$ & 3.91e-4 & 5.10e-3 & 4.18e-2 & & 6.49e-4 & 1.19e-2 & 8.47e-2 & & 9.78e-3 & 1.11e-1 & 6.37e-1 \\
  $1-10^{-9}$ & 3.06e-4 & 4.14e-3 & 3.48e-2 & & 5.05e-4 & 1.01e-2 & 6.99e-2 & & 8.91e-3 & 1.03e-1 & 5.78e-1 \\
  \bottomrule[.4mm]
  \end{tabular}
\end{table}
\begin{table}[h]\centering
  \caption{The relative space-time Frobenius errors of the projected FOM solutions and ROM computations in the case $R=0.22$ and $x_c=-0.65$ of Test 3.}
  \label{tb:num_cyl_xr_froberr_r0d22_x-0d65}
  \begin{tabular}{@{}lccccccccccc@{}}
  \toprule[.5mm]
  & \multicolumn{3}{l}{FOM projection} & & \multicolumn{3}{l}{SBM-iROM} & & \multicolumn{3}{l}{SBM-ROM} \\ \cmidrule[.2mm](lr){2-4} \cmidrule[.2mm](lr){6-8} \cmidrule[.2mm](l){10-12}
  $\mu_{\textrm{pod}}$ & $h$ & $hv_1$ & $hv_2$ & & $h$ & $hv_1$ & $hv_2$ & & $h$ & $hv_1$ & $hv_2$ \\ \cmidrule[.3mm](l){2-12}
  $1-10^{-5}$ & 2.76e-3 & 2.93e-2 & 1.50e-1 & & 3.90e-3 & 5.27e-2 & 2.69e-1 & & 1.85e-2 & 2.09e-1 & 8.46e-1 \\
  $1-10^{-6}$ & 1.68e-3 & 1.93e-2 & 1.14e-1 & & 2.45e-3 & 3.39e-2 & 2.08e-1 & & 1.80e-2 & 2.11e-1 & 8.75e-1 \\
  $1-10^{-7}$ & 7.22e-4 & 7.85e-3 & 5.48e-2 & & 1.10e-3 & 1.64e-2 & 9.50e-2 & & 1.77e-2 & 2.06e-1 & 9.13e-1 \\
  $1-10^{-8}$ & 4.14e-4 & 4.83e-3 & 3.34e-2 & & 7.16e-4 & 1.14e-2 & 6.62e-2 & & 1.71e-2 & 1.99e-1 & 8.37e-1 \\
  $1-10^{-9}$ & 2.87e-4 & 3.84e-3 & 2.49e-2 & & 5.67e-4 & 8.76e-3 & 4.91e-2 & & 1.62e-2 & 1.97e-1 & 7.60e-1 \\
  \bottomrule[.4mm]
  \end{tabular}
\end{table}
\begin{table}[h]\centering
  \caption{The relative space-time Frobenius errors of the projected FOM solutions and ROM computations in the case $R=0.08$ and $x_c=-0.15$ of Test 3.}
  \label{tb:num_cyl_xr_froberr_r0d08_x-0d15}
  \begin{tabular}{@{}lccccccccccc@{}}
  \toprule[.5mm]
  & \multicolumn{3}{l}{FOM projection} & & \multicolumn{3}{l}{SBM-iROM} & & \multicolumn{3}{l}{SBM-ROM} \\ \cmidrule[.2mm](lr){2-4} \cmidrule[.2mm](lr){6-8} \cmidrule[.2mm](l){10-12}
  $\mu_{\textrm{pod}}$ & $h$ & $hv_1$ & $hv_2$ & & $h$ & $hv_1$ & $hv_2$ & & $h$ & $hv_1$ & $hv_2$ \\ \cmidrule[.3mm](l){2-12}
  $1-10^{-5}$ & 1.33e-3 & 2.96e-2 & 3.63e-1 & & 2.47e-3 & 5.05e-2 & 6.32e-1 & & 5.38e-2 & 3.06e-1 & 1.99e+0 \\
  $1-10^{-6}$ & 9.88e-4 & 2.37e-2 & 3.14e-1 & & 2.04e-3 & 4.90e-2 & 4.72e-1 & & 5.20e-2 & 3.62e-1 & 2.40e+0 \\
  $1-10^{-7}$ & 1.05e-3 & 1.81e-2 & 2.28e-1 & & 1.36e-3 & 4.57e-2 & 4.35e-1 & & 3.58e-2 & 4.57e-1 & 6.01e+0 \\
  $1-10^{-8}$ & 6.71e-4 & 1.32e-2 & 1.45e-1 & & 8.21e-4 & 3.61e-2 & 3.12e-1 & & 2.57e-2 & 3.14e-1 & 2.76e+0 \\
  $1-10^{-9}$ & 5.72e-4 & 1.07e-2 & 1.17e-1 & & 6.48e-4 & 2.72e-2 & 2.46e-1 & & 2.39e-2 & 3.32e-1 & 2.55e+0 \\
  \bottomrule[.4mm]
  \end{tabular}
\end{table}
\begin{table}[h]\centering
  \caption{The relative space-time Frobenius errors of the projected FOM solutions and ROM computations in the case $R=0.15$ and $x_c=-0.15$ of Test 3.}
  \label{tb:num_cyl_xr_froberr_r0d15_x-0d15}
  \begin{tabular}{@{}lccccccccccc@{}}
  \toprule[.5mm]
  & \multicolumn{3}{l}{FOM projection} & & \multicolumn{3}{l}{SBM-iROM} & & \multicolumn{3}{l}{SBM-ROM} \\ \cmidrule[.2mm](lr){2-4} \cmidrule[.2mm](lr){6-8} \cmidrule[.2mm](l){10-12}
  $\mu_{\textrm{pod}}$ & $h$ & $hv_1$ & $hv_2$ & & $h$ & $hv_1$ & $hv_2$ & & $h$ & $hv_1$ & $hv_2$ \\ \cmidrule[.3mm](l){2-12}
  $1-10^{-5}$ & 1.77e-3 & 2.53e-2 & 1.40e-1 & & 3.26e-3 & 5.45e-2 & 2.99e-1 & & 1.21e-2 & 1.55e-1 & 6.19e-1 \\
  $1-10^{-6}$ & 1.15e-3 & 1.64e-2 & 1.13e-1 & & 1.58e-3 & 3.19e-2 & 1.73e-1 & & 1.12e-2 & 1.21e-1 & 6.49e-1 \\
  $1-10^{-7}$ & 6.90e-4 & 7.94e-3 & 6.54e-2 & & 1.05e-3 & 2.59e-2 & 1.41e-1 & & 8.07e-3 & 8.99e-2 & 4.81e-1 \\
  $1-10^{-8}$ & 3.46e-4 & 4.74e-3 & 3.38e-2 & & 5.21e-4 & 1.36e-2 & 7.70e-2 & & 7.07e-3 & 8.09e-2 & 4.03e-1 \\
  $1-10^{-9}$ & 2.71e-4 & 3.77e-3 & 2.68e-2 & & 4.34e-4 & 1.12e-2 & 6.19e-2 & & 6.19e-3 & 7.14e-2 & 3.46e-1 \\
  \bottomrule[.4mm]
  \end{tabular}
\end{table}
\begin{table}[h]\centering
  \caption{The relative space-time Frobenius errors of the projected FOM solutions and ROM computations in the case $R=0.22$ and $x_c=-0.15$ of Test 3.}
  \label{tb:num_cyl_xr_froberr_r0d22_x-0d15}
  \begin{tabular}{@{}lccccccccccc@{}}
  \toprule[.5mm]
  & \multicolumn{3}{l}{FOM projection} & & \multicolumn{3}{l}{SBM-iROM} & & \multicolumn{3}{l}{SBM-ROM} \\ \cmidrule[.2mm](lr){2-4} \cmidrule[.2mm](lr){6-8} \cmidrule[.2mm](l){10-12}
  $\mu_{\textrm{pod}}$ & $h$ & $hv_1$ & $hv_2$ & & $h$ & $hv_1$ & $hv_2$ & & $h$ & $hv_1$ & $hv_2$ \\ \cmidrule[.3mm](l){2-12}
  $1-10^{-5}$ & 2.23e-3 & 2.72e-2 & 1.09e-1 & & 3.44e-3 & 5.10e-2 & 2.37e-1 & & 1.96e-2 & 2.46e-1 & 7.27e-1 \\
  $1-10^{-6}$ & 1.42e-3 & 1.94e-2 & 8.15e-2 & & 2.23e-3 & 3.29e-2 & 1.35e-1 & & 1.62e-2 & 1.89e-1 & 6.91e-1 \\
  $1-10^{-7}$ & 6.86e-4 & 8.64e-3 & 4.42e-2 & & 1.16e-3 & 2.10e-2 & 9.41e-2 & & 1.43e-2 & 1.71e-1 & 6.22e-1 \\
  $1-10^{-8}$ & 3.27e-4 & 5.09e-3 & 2.49e-2 & & 5.72e-4 & 1.15e-2 & 5.33e-2 & & 1.19e-2 & 1.43e-1 & 5.20e-1 \\
  $1-10^{-9}$ & 2.24e-4 & 3.98e-3 & 1.95e-2 & & 4.49e-4 & 9.05e-3 & 3.78e-2 & & 1.00e-2 & 1.23e-1 & 4.40e-1 \\
  \bottomrule[.4mm]
  \end{tabular}
\end{table}
\begin{table}[h]\centering
  \caption{The relative space-time Frobenius errors of the projected FOM solutions and ROM computations in the case $R=0.08$ and $x_c=0.3$ of Test 3.}
  \label{tb:num_cyl_xr_froberr_r0d08_x0d3}
  \begin{tabular}{@{}lccccccccccc@{}}
  \toprule[.5mm]
  & \multicolumn{3}{l}{FOM projection} & & \multicolumn{3}{l}{SBM-iROM} & & \multicolumn{3}{l}{SBM-ROM} \\ \cmidrule[.2mm](lr){2-4} \cmidrule[.2mm](lr){6-8} \cmidrule[.2mm](l){10-12}
  $\mu_{\textrm{pod}}$ & $h$ & $hv_1$ & $hv_2$ & & $h$ & $hv_1$ & $hv_2$ & & $h$ & $hv_1$ & $hv_2$ \\ \cmidrule[.3mm](l){2-12}
  $1-10^{-5}$ & 1.89e-3 & 4.15e-2 & 5.38e-1 & & 2.96e-3 & 6.30e-2 & 7.55e-1 & & 5.34e-2 & 2.79e-1 & 2.13e+0 \\
  $1-10^{-6}$ & 1.78e-3 & 3.43e-2 & 4.99e-1 & & 2.26e-3 & 4.67e-2 & 5.94e-1 & & 4.90e-2 & 3.17e-1 & 2.42e+0 \\
  $1-10^{-7}$ & 1.55e-3 & 2.48e-2 & 3.59e-1 & & 1.68e-3 & 4.43e-2 & 4.89e-1 & & 3.18e-2 & 4.06e-1 & 5.51e+0 \\
  $1-10^{-8}$ & 1.04e-3 & 1.67e-2 & 1.85e-1 & & 9.46e-4 & 3.52e-2 & 3.23e-1 & & 2.31e-2 & 2.98e-1 & 2.50e+0 \\
  $1-10^{-9}$ & 8.60e-4 & 1.33e-2 & 1.42e-1 & & 6.98e-4 & 2.95e-2 & 2.80e-1 & & 2.18e-2 & 2.92e-1 & 2.29e+0 \\
  \bottomrule[.4mm]
  \end{tabular}
\end{table}
\begin{table}[h]\centering
  \caption{The relative space-time Frobenius errors of the projected FOM solutions and ROM computations in the case $R=0.15$ and $x_c=0.3$ of Test 3.}
  \label{tb:num_cyl_xr_froberr_r0d15_x0d3}
  \begin{tabular}{@{}lccccccccccc@{}}
  \toprule[.5mm]
  & \multicolumn{3}{l}{FOM projection} & & \multicolumn{3}{l}{SBM-iROM} & & \multicolumn{3}{l}{SBM-ROM} \\ \cmidrule[.2mm](lr){2-4} \cmidrule[.2mm](lr){6-8} \cmidrule[.2mm](l){10-12}
  $\mu_{\textrm{pod}}$ & $h$ & $hv_1$ & $hv_2$ & & $h$ & $hv_1$ & $hv_2$ & & $h$ & $hv_1$ & $hv_2$ \\ \cmidrule[.3mm](l){2-12}
  $1-10^{-5}$ & 2.24e-3 & 3.76e-2 & 2.13e-1 & & 3.34e-3 & 5.79e-2 & 3.86e-1 & & 1.28e-2 & 1.68e-1 & 7.46e-1 \\
  $1-10^{-6}$ & 1.88e-3 & 2.88e-2 & 1.87e-1 & & 2.47e-3 & 4.78e-2 & 2.81e-1 & & 1.32e-2 & 1.52e-1 & 7.98e-1 \\
  $1-10^{-7}$ & 1.20e-3 & 1.38e-2 & 1.35e-1 & & 1.51e-3 & 3.58e-2 & 2.25e-1 & & 9.16e-3 & 1.04e-1 & 6.12e-1 \\
  $1-10^{-8}$ & 6.31e-4 & 9.09e-3 & 7.39e-2 & & 8.05e-4 & 2.43e-2 & 1.57e-1 & & 6.38e-3 & 7.33e-2 & 4.45e-1 \\
  $1-10^{-9}$ & 5.23e-4 & 7.44e-3 & 5.71e-2 & & 7.21e-4 & 1.98e-2 & 1.30e-1 & & 5.30e-3 & 6.29e-2 & 3.71e-1 \\
  \bottomrule[.4mm]
  \end{tabular}
\end{table}
\begin{table}[h]\centering
  \caption{The relative space-time Frobenius errors of the projected FOM solutions and ROM computations in the case $R=0.22$ and $x_c=0.3$ of Test 3.}
  \label{tb:num_cyl_xr_froberr_r0d22_x0d3}
  \begin{tabular}{@{}lccccccccccc@{}}
  \toprule[.5mm]
  & \multicolumn{3}{l}{FOM projection} & & \multicolumn{3}{l}{SBM-iROM} & & \multicolumn{3}{l}{SBM-ROM} \\ \cmidrule[.2mm](lr){2-4} \cmidrule[.2mm](lr){6-8} \cmidrule[.2mm](l){10-12}
  $\mu_{\textrm{pod}}$ & $h$ & $hv_1$ & $hv_2$ & & $h$ & $hv_1$ & $hv_2$ & & $h$ & $hv_1$ & $hv_2$ \\ \cmidrule[.3mm](l){2-12}
  $1-10^{-5}$  & 2.61e-3 & 3.66e-2 & 1.44e-1 & & 4.12e-3 & 6.08e-2 & 2.72e-1 & & 2.14e-2 & 2.54e-1 & 8.82e-1 \\
  $1-10^{-6}$ & 2.02e-3 & 2.72e-2 & 1.13e-1 & & 2.86e-3 & 4.45e-2 & 2.09e-1 & & 1.74e-2 & 2.06e-1 & 7.61e-1 \\
  $1-10^{-7}$ & 9.48e-4 & 1.03e-2 & 6.58e-2 & & 1.67e-3 & 2.29e-2 & 1.14e-1 & & 1.52e-2 & 1.83e-1 & 7.02e-1 \\
  $1-10^{-8}$ & 4.93e-4 & 6.32e-3 & 3.85e-2 & & 8.07e-4 & 1.36e-2 & 7.29e-2 & & 1.07e-2 & 1.30e-1 & 5.20e-1 \\
  $1-10^{-9}$ & 3.85e-4 & 5.07e-3 & 2.85e-2 & & 6.24e-4 & 1.23e-2 & 6.33e-2 & & 8.90e-3 & 1.11e-1 & 4.25e-1 \\
  \bottomrule[.4mm]
  \end{tabular}
\end{table}
\begin{table}[h]\centering
  \caption{The relative space-time Frobenius errors of the projected FOM solutions and ROM computations in the case $R=0.08$ and $x_c=0.8$ of Test 3.}
  \label{tb:num_cyl_xr_froberr_r0d08_x0d8}
  \begin{tabular}{@{}lccccccccccc@{}}
  \toprule[.5mm]
  & \multicolumn{3}{l}{FOM projection} & & \multicolumn{3}{l}{SBM-iROM} & & \multicolumn{3}{l}{SBM-ROM} \\ \cmidrule[.2mm](lr){2-4} \cmidrule[.2mm](lr){6-8} \cmidrule[.2mm](l){10-12}
  $\mu_{\textrm{pod}}$ & $h$ & $hv_1$ & $hv_2$ & & $h$ & $hv_1$ & $hv_2$ & & $h$ & $hv_1$ & $hv_2$ \\ \cmidrule[.3mm](l){2-12}
  $1-10^{-5}$ & 2.61e-3 & 5.21e-2 & 5.99e-1 & & 3.40e-3 & 6.85e-2 & 7.29e-1 & & 5.36e-2 & 2.61e-1 & 2.03e+0 \\
  $1-10^{-6}$ & 2.52e-3 & 4.65e-2 & 5.43e-1 & & 3.43e-3 & 6.70e-2 & 6.73e-1 & & 5.37e-2 & 2.70e-1 & 2.45e+0 \\
  $1-10^{-7}$ & 2.32e-3 & 3.59e-2 & 4.30e-1 & & 2.41e-3 & 5.86e-2 & 6.16e-1 & & 3.78e-2 & 3.55e-1 & 5.02e+0 \\
  $1-10^{-8}$ & 1.67e-3 & 2.48e-2 & 2.75e-1 & & 1.77e-3 & 4.75e-2 & 4.64e-1 & & 2.95e-2 & 2.91e-1 & 2.89e+0 \\
  $1-10^{-9}$ & 1.42e-3 & 2.20e-2 & 2.24e-1 & & 1.66e-3 & 4.09e-2 & 3.96e-1 & & 2.88e-2 & 2.99e-1 & 2.70e+0 \\
  \bottomrule[.4mm]
  \end{tabular}
\end{table}
\begin{table}[h]\centering
  \caption{The relative space-time Frobenius errors of the projected FOM solutions and ROM computations in the case $R=0.15$ and $x_c=0.8$ of Test 3.}
  \label{tb:num_cyl_xr_froberr_r0d15_x0d8}
  \begin{tabular}{@{}lccccccccccc@{}}
  \toprule[.5mm]
  & \multicolumn{3}{l}{FOM projection} & & \multicolumn{3}{l}{SBM-iROM} & & \multicolumn{3}{l}{SBM-ROM} \\ \cmidrule[.2mm](lr){2-4} \cmidrule[.2mm](lr){6-8} \cmidrule[.2mm](l){10-12}
  $\mu_{\textrm{pod}}$ & $h$ & $hv_1$ & $hv_2$ & & $h$ & $hv_1$ & $hv_2$ & & $h$ & $hv_1$ & $hv_2$ \\ \cmidrule[.3mm](l){2-12}
  $1-10^{-5}$ & 3.45e-3 & 4.98e-2 & 3.28e-1 & & 4.67e-3 & 7.07e-2 & 4.39e-1 & & 1.51e-2 & 1.87e-1 & 1.07e+0 \\
  $1-10^{-6}$ & 2.99e-3 & 4.20e-2 & 2.78e-1 & & 4.23e-3 & 6.23e-2 & 3.74e-1 & & 1.27e-2 & 1.50e-1 & 8.85e-1 \\
  $1-10^{-7}$ & 2.20e-3 & 2.45e-2 & 1.78e-1 & & 2.89e-3 & 4.47e-2 & 2.60e-1 & & 1.20e-2 & 1.28e-1 & 9.38e-1 \\
  $1-10^{-8}$ & 1.15e-3 & 1.21e-2 & 9.76e-2 & & 1.63e-3 & 2.73e-2 & 2.00e-1 & & 9.24e-3 & 9.44e-2 & 7.28e-1 \\
  $1-10^{-9}$ & 9.35e-4 & 9.40e-3 & 7.69e-2 & & 1.26e-3 & 2.20e-2 & 1.55e-1 & & 9.06e-3 & 8.98e-2 & 7.16e-1 \\
  \bottomrule[.4mm]
  \end{tabular}
\end{table}
\begin{table}[h]\centering
  \caption{The relative space-time Frobenius errors of the projected FOM solutions and ROM computations in the case $R=0.22$ and $x_c=0.8$ of Test 3.}
  \label{tb:num_cyl_xr_froberr_r0d22_x0d8}
  \begin{tabular}{@{}lccccccccccc@{}}
  \toprule[.5mm]
  & \multicolumn{3}{l}{FOM projection} & & \multicolumn{3}{l}{SBM-iROM} & & \multicolumn{3}{l}{SBM-ROM} \\ \cmidrule[.2mm](lr){2-4} \cmidrule[.2mm](lr){6-8} \cmidrule[.2mm](l){10-12}
  $\mu_{\textrm{pod}}$ & $h$ & $hv_1$ & $hv_2$ & & $h$ & $hv_1$ & $hv_2$ & & $h$ & $hv_1$ & $hv_2$ \\ \cmidrule[.3mm](l){2-12}
  $1-10^{-5}$ & 4.52e-3 & 4.84e-2 & 2.77e-1 & & 6.92e-3 & 8.07e-2 & 4.10e-1 & & 2.27e-2 & 2.36e-1 & 7.13e-1 \\
  $1-10^{-6}$ & 3.69e-3 & 3.96e-2 & 2.24e-1 & & 5.78e-3 & 6.34e-2 & 3.22e-1 & & 2.18e-2 & 2.13e-1 & 8.23e-1 \\
  $1-10^{-7}$ & 2.71e-3 & 2.56e-2 & 1.64e-1 & & 3.92e-3 & 4.67e-2 & 2.56e-1 & & 2.15e-2 & 2.06e-1 & 8.63e-1 \\
  $1-10^{-8}$ & 1.32e-3 & 1.30e-2 & 8.69e-2 & & 1.93e-3 & 2.61e-2 & 1.63e-1 & & 1.91e-2 & 1.69e-1 & 7.66e-1 \\
  $1-10^{-9}$ & 9.58e-4 & 9.72e-3 & 6.48e-2 & & 1.56e-3 & 2.21e-2 & 1.27e-1 & & 1.87e-2 & 1.70e-1 & 7.52e-1 \\
  \bottomrule[.4mm]
  \end{tabular}
\end{table}

%% file: section_concl.tex
\section{Concluding remarks and future developments}\label{sec:conclusions}
In this article we have proposed and analyzed the coupling between the Shifted Boundary Method and POD-Galerkin methods for reduced order modelling in presence of geometrical parameters considering a case of hyperbolic systems. The methodology has been applied to shallow water equations discretized using an explicit time integration scheme and tested on three numerical benchmarks of increasing complexity.

In order to tackle one of the issues arising with the coupling of immersed methods and projection-based reduced-order models we proposed a new preprocessing technique applied on the full-order snapshots prior to the POD basis computation. This interpolation techniques demonstrates to increase considerably the accuracy of the results with respect to the unpreprocessed approach without any increase of the online computational costs.

The proposed preprocessing approach produces accurate results for cases that would otherwise yield completely unreliable results. This fact was particularly evident when the reduced-order model required to activate grid nodes that remained inactive for all training snapshots. 

From a general remark, we observe that the coupling of immersed methods and projection based reduced models suffer from the necessity to use a large number of training points in order to produce accurate results during the online computations.

In future works we aim to study possible approaches in order to circumvent this
limitation applying suitable transformations to the POD basis. In order to
speed-up the online computations, other interesting aspects are related to the
hyper-reduction of the problem using residual evaluations only in a limited
number of nodes. Despite the fact that hyper-reduction techniques are now well
developed, their application to immersed methods might be challenging. 

We also plan to extend the proposed methodologies to fluid structure interaction problems.

%% file: section_aknw.tex
\section*{Acknowledgments}
This research has been supported by the Army Research Office (ARO) under Grant W911NF-18-1-0308 (GS), the U.S. National Science Foundation under Grant 2137934, European Union Funding for Research and Innovation -Horizon 2020 Program- in the framework of European Research Council Executive Agency:
Consolidator Grant H2020 ERC CoG 2015 AROMA-CFD project 681447 ``Advanced Reduced Order Methods with Applications in Computational Fluid Dynamics" (PI Prof. Gianluigi Rozza). We also acknowledge the INDAM-GNCS project ``Tecniche Numeriche Avanzate per Applicazioni Industriali", and by project FSE-European Social Fund-HEaD ``Higher Education and Development" SISSA operazione 1, Regione Autonoma Friuli-Venezia Giulia, the HFRI and GSRT under grant agreement No 1115,  the ``First Call for H.F.R.I. Research Projects to support Faculty members and Researchers and the procurement of high-cost research equipment" grant 3270, and the support  of the National Infrastructures for Research and Technology S.A. (GRNET S.A.) under project ID pa190902.